\documentclass[10pt,fleqn]{article}

\usepackage{url}

\usepackage[T1]{fontenc}
\usepackage{textcomp}

\usepackage{graphicx}
\usepackage{mathpazo}
\usepackage{amsmath}
\usepackage{amssymb}
\usepackage{amsfonts}
\usepackage{subfigure}
\usepackage{epsfig}
\usepackage{epstopdf}

\newcommand{\thefont}[2]{\fontsize{#1}{#2}\fontshape{n}\selectfont}
\newcommand{\1}{\rlap{\thefont{10pt}{12pt}1}\kern.16em\rlap{\thefont{11pt}{13.2pt}1}\kern.4em}

\renewcommand{\baselinestretch} {1.1}

\makeatletter
\def\singlespace{\def\baselinestretch{1}\@normalsize}

\usepackage{color}

\title{{\sc On the usefulness of Meyer wavelets for deconvolution and density estimation}}

\author{
  {\em J\'er\'emie Bigot}}

\date{March 2009}

\usepackage{apacite}
	
	\renewcommand{\BBAA}{and}

\newcommand{\opO}{\ensuremath{\mathcal O}}

\newcommand{\reals}{\ensuremath{{\mathbb R}}}
\newcommand{\CC}{\ensuremath{{\mathbb C}}}
\newcommand{\RR}{\reals}
\newcommand{\ZZ}{\ensuremath{{\mathbb Z}}}

\newcommand{\PP}{\ensuremath{{\mathbb P}}}
\newcommand{\EE}{\ensuremath{{\mathbb E}}}
\newcommand{\var}{\mbox{Var}}
\newcommand{\pen}{\mbox{pen}}

\newtheorem{ass}{Assumption}
\newtheorem{theo}{Theorem}
\newtheorem{prop}{Proposition}
\newtheorem{lemma}{Lemma}

\numberwithin{equation}{section}
\numberwithin{ass}{section}
\numberwithin{theo}{section}
\numberwithin{prop}{section}
\numberwithin{lemma}{section}
\numberwithin{definition}{section}
\numberwithin{rmq}{section}

\begin{document}

\maketitle

\thispagestyle{empty}

\begin{abstract} 
The aim of this paper is to show the usefulness of Meyer wavelets for the classical problem of  density estimation and for density deconvolution from noisy observations. By using such wavelets, the computation of the empirical  wavelet coefficients relies on the fast Fourier transform  and on the fact that Meyer wavelets are band-limited functions. This makes such estimators very simple to compute and this avoids the problem of evaluating wavelets at non-dyadic points which is the main drawback of classical wavelet-based density estimators.  Our approach is based on term-by-term thresholding of the empirical wavelet coefficients with random thresholds depending on an estimation of the variance of each  coefficient. Such estimators are shown to achieve the same performances of an oracle estimator up to a logarithmic term. These estimators also achieve near-minimax rates of convergence over a large class of Besov spaces. A simulation study is proposed to show the good finite sample performances of the estimator for both problems of direct density estimation and density deconvolution.
\end{abstract}

{\footnotesize

%\vspace{0.5cm}



\subsubsection*{Affiliations}
Institut de Math\'ematiques de Toulouse, Universit\'e de Toulouse et CNRS (UMR 5219), \url{Jeremie.Bigot@math.univ-toulouse.fr}

\subsubsection*{Acknowledgements}

We gratefully acknowledge Yves Rozenholc for providing the Matlab code to compute the model selection estimator.

\setcounter{page}{0}
%\newpage
}

\section{Introduction}

Density estimation is a well-known problem in statistics that has been thoroughly studied. It consists in estimating an unknown probability density function $f$ from an independent and identically distributed (iid) sample  of random variables $X_{i}$ for $i=1,\ldots,n$, with $n$ representing the sample size. Wavelet decomposition is known to be a powerful method for nonparametric estimation, see e.g.  \citeA{DJKP95jrssb}. The advantages of wavelet methods is their ability in estimating spatially inhomogeneous functions. They can be used to estimate functions in Besov spaces with optimal rates of convergence, and have therefore received special attention in the literature over the last two decades, in particular  for density estimation, see e.g. \citeA{DJKP96aos} and \citeA{vida} for a detailed review on the subject. 

For a given scaling function $\phi$ and mother wavelet $\psi$, the scaling and wavelet coefficients are usually estimated as \citeA{DJKP96aos}
\begin{equation}
\hat{c}_{j_{0},k} = \frac{1}{n} \sum_{i=1}^{n} \phi_{j_{0},k}(X_{i}) \mbox{ and } \hat{\beta}_{j,k} =  \frac{1}{n} \sum_{i=1}^{n}  \psi_{j,k}(X_{i}) \label{eq:hatcoefintro}
\end{equation}
where $ \phi_{j_{0},k}(x)  =  \phi(2^{j_{0}}x-k)$,  $ \psi_{j,k}(x)  =  \psi(2^{j}x-k)$, and $j_{0}$ denotes the usual coarse level of resolution. A hard-thresholding estimator of $f$ then takes the form
\begin{equation}
\hat{f}_{n}(x) = \sum_{k} \hat{c}_{j_{0},k}   \phi_{j_{0},k}(x) + \sum_{j=j_{0}}^{j_{1}} \sum_{k}  \hat{\beta}_{j,k}  \1_{\{|  \hat{\beta}_{j,k}  | \geqslant \tau_{j,k}\}}   \psi_{j,k}(x), \label{eq:hatf}
\end{equation}
where  $j_{1}$ is an appropriate frequency cut-off, and where the $\tau_{j,k}$'s are appropriate thresholds (positive numbers) that are possibly level-dependent. However, in practice the computation of the coefficients (\ref{eq:hatcoefintro}) requires some numerical approximation by interpolation  as one typically only knows the values of the functions $\phi_{j_{0},k}$ and $ \psi_{j,k}$ at dyadic points. Various approaches have been used to approximate numerically the coefficients in (\ref{eq:hatcoefintro}). For instance \citeA{kookoo}, \citeA{antgreg} use a binning method followed by the the standard discrete wavelet transform, while the algorithm proposed in \citeA{nason} is based on an approximation of the scaling coefficients at a sufficiently small level of resolution.

%{\bf A completer avec references sur approximation et calcul des coefficients d'ondelettes en densite } \\
 In this paper, we propose to avoid the use of such numerical approximation schemes. This is achieved by using Meyer wavelets which are band-limited functions. Indeed, using such wavelets and thanks to the Plancherel's identity, the empirical coefficients can be easily computed from the Fourier transform of the data $X_{i}, i=1,\ldots,n$. Such an approach therefore takes full advantages of the fast Fourier transform and of existing fast algorithms for Meyer wavelet decomposition developed by \citeA{K94}. 

As Meyer wavelets are band-limited functions, they have been recently used for deconvolution problems in nonparametric regression by \citeA{JKPR04jrssb}, \citeA{sapatinaspensky}, and for density deconvolution by \citeA{PV99aos,FK02ieee,BVB07}. Moreover, the use of such wavelets leads to fast algorithms, see   \citeA{K94} and \citeA{rast07}.

Density deconvolution is the problem of estimating the function $f$ when the observation of the random variable $X$ is contaminated by an independent additive noise. In this case, the observations at hand are a sample of variables $Y_{i},i=1,\ldots,n$ such that
\begin{equation}
Y_{i} = X_{i} + \epsilon_{i}, \quad i=1,\ldots,n, \label{eq:densdec}
\end{equation}
where $X_{i}$ are iid variables with unknown density $f$, and $\epsilon_{i}$ are iid variables with {\it known} density $h$ which represents some additive noise independent of the $X_{i}$'s. Density estimation from a noisy sample is of fundamental importance in many practical situations, and applications can be found in communication theory \citeA{M03ieee}, experimental physics \shortcite<e.g.>{KSTL03cpc} or econometrics \citeA{PVR02}. In this setting, the density of the observed variables $Y_{i}$ is the convolution of the density $f$ with the density $h$ of the additive noise. Hence, the problem of estimating $f$ relates to nonparametric methods of deconvolution which is a widely studied inverse problem in statistics and signal processing. However the indirect observation of the data leads to different optimality properties, for instance in terms of rate of convergence, than the direct problem of density estimation without an additive error. Standard techniques recently studied for density deconvolution include model selection \citeA{CRT06cjs},  kernel smoothing \citeA{CH88jasa}, spline deconvolution \citeA{K99sjs}, spectral cut-off \citeA{J08aos}  and wavelet thresholding \citeA{PV99aos,FK02ieee,BVB07}, to name but a few. Again, Meyer wavelets can be used to easily compute estimators of the wavelet coefficients of $f$ by using the Fourier coefficients of the noise density $h$ without using any numerical approximation scheme.

The second contribution of this paper is the use of random thresholds $\tau_{j,k}$. Classically, the thresholds used in wavelet density estimation are deterministic, but such thresholds may be too large in practice as they do not take into account the fact that the variance of each empirical wavelet coefficient $\hat{\beta}_{j,k}$ depends on its location $(j,k)$.  The use of random thresholds has been recently proposed in \citeA{judlam} in the context of density estimation, and in \citeA{reyriv} for Poisson intensity estimation.  However, the estimation procedure in  \citeA{judlam} and \citeA{reyriv} is different from the one we propose, since it is based on biorthogonal bases and on the use of the Haar basis to compute the wavelet coefficients. %These data-driven thresholds are based on an estimation of the variance of the empirical wavelet coefficients  $\hat{\beta}_{j,k}$.
The use of data based threshold exploiting the variance structure of the empirical wavelet coefficients is also proposed in \citeA{nason}, but the theoretical properties of the resulting algorithms are not studied. In this paper, we show that using Meyer wavelets allows one to compute easily an estimation of an upper bound of the variance of the $\hat{\beta}_{j,k}$'s that is then used to compute random thresholds $\tau_{j,k}$. 

Then, a third contribution of this paper is that the resulting hard-thresholding estimators are shown to attain the same performances (up to  logarithmic terms) of an ideal estimator, called oracle as its computation depends on unknown quantities such as the variance of the $\hat{\beta}_{j,k}$'s or the magnitude of the true wavelet coefficients. Oracle inequalities is an active research area in nonparametric statistics (see e.g. \citeA{johnstone}, \citeA{candes} for detailed expositions) which has recently gained popularity.  Deriving an oracle inequality is the problem of bounding the risk of a statistical procedure by the performances of an ideal estimator which represents the best model for the function to recover.  Oracle inequalities are currently used in many different contexts in statistics. They have been introduced in \citeA{donojohn} for nonparametric regression with wavelets, then used by \citeA{cgpt} for inverse problems in the white noise model, by \citeA{rigo}, \citeA{caste}, \citeA{efro}, \citeA{bunea} for density estimation problems, and by \citeA{reyriv} for estimating the intensity of a Poisson process, to name but a few. 

The rest of this paper is organized as follows. In Section \ref{sec:wav}, we provide some background on Meyer wavelets and we define the corresponding wavelet estimators with random thresholds for both direct density estimation and density deconvolution. In Section \ref{sec:oracle}, it is explained how the performances of such estimates can be compared to those of an oracle estimate, which leads to non-asymptotic oracle inequalities. Asymptotic properties of the estimators are then studied in Section \ref{sec:minimax}. Depending on the problem at hand.% (either direct density estimation or density deconvolution)
, these estimators are shown to achieve near-minimax rates of convergence over a large class of Besov spaces. Finally, a simulation study is proposed in Section \ref{sec:simu} to evaluate the numerical performances of our estimators and to compare them with other procedures existing in the literature. The proofs of the main theorems are gathered in a  technical Appendix.

\section{Density estimation with Meyer wavelets} \label{sec:wav}

In what follows, it will be assumed that the density function $f$ of the $X_{i}$'s has a compact support included in $[0,1]$.  Of course, assuming that the support of $f$ is included in $[0,1]$ would not hold in many practical applications and this is mainly made for mathematical convenience to simplify the presentation of the estimator. If the range of the data is outside $[0,1]$, one can simply rescale and center them such that they fall into $[0,1]$, and then apply the inverse transformation to the estimated density.% In this section, we recall some general properties of Meyer wavelet bases on the interval $[0,1]$, and we explain how they can be used to quickly compute an estimator of $f$ for both the direct problem of density estimation and for density deconvolution. 

\subsection{Wavelet decomposition and the periodized Meyer wavelet basis} 

Let $(\phi^{\ast},\psi^{\ast})$ be the Meyer scaling and wavelet function respectively  (see \citeA{M92} for further details). It is constructed from a scaling function $\phi^{\ast}$ with Fourier transform
	\begin{equation*}
	\tilde{\phi}(\omega) = \int_{\RR} \phi^{\ast}(x) e^{-i \omega x} dx = \left\{
	\begin{array}{cll}
	\frac{\tilde{h}(\omega/2)}{\sqrt{2}}& \mbox{ if } & |\omega| \leqslant 4\pi/3,\\
	0 & \mbox{ if } & |\omega| > 4\pi/3,
	\end{array}
	\right.
	\end{equation*}
where $\tilde{h} : \CC \to \RR$ is a smooth function chosen as a polynomial of degree 3 in our simulations. Scaling and wavelet function at scale $j \geq 0$ are defined by
$$
\phi_{j,k}^{\ast}(x) = 2^{j/2} \phi^{\ast}(2^{j}x-k) \mbox{ and } \psi_{j,k}^{\ast}(x) = 2^{j/2} \psi^{\ast}(2^{j}x-k), \; k=0,\ldots,2^{j}-1.
$$
As in  \citeA{JKPR04jrssb} , one can then define the periodized Meyer wavelet basis of $L^{2}([0,1])$ (the space of squared integrable functions on $[0,1]$) by periodizing the functions  $(\phi^{\ast},\psi^{\ast})$ i.e.
$$
\phi_{j,k}(x) = 2^{j/2} \sum_{i \in \ZZ} \phi^{\ast}(2^{j}(x+i)-k) \mbox{ and } \psi_{j,k}(x) = 2^{j/2} \sum_{i \in \ZZ} \psi^{\ast}(2^{j}(x+i)-k), \; k=0,\ldots,2^{j}-1.
$$
For any function $f$ of $L^{2}([0,1])$, its wavelet decomposition can be written as:
	\begin{equation*}
	f = \sum_{k =0}^{2^{j_{0}}-1} c_{j_{0},k} \phi_{j_{0},k} + \sum_{j = j_{0}}^{+ \infty} \sum_{k=0}^{2^{j}-1} \beta_{j,k} \psi_{j,k},
	\end{equation*}
where $c_{j_{0},k} = \langle f , \phi_{j_{0},k} \rangle = \int_{0}^{1} f(u) \phi_{j_{0},k}(u) du$, $\beta_{j,k} = \langle f,\psi_{j,k} \rangle =  \int_{0}^{1} f(u) ,\psi_{j,k}(u) du$ and $j_{0} \geq 0$ denotes the usual coarse level of resolution. Moreover, the $L^{2}$ norm of $f$ is given by
$$
\|f\|^{2} = \sum_{k =0}^{2^{j_{0}}-1} c_{j_{0},k}^{2} + \sum_{j = j_{0}}^{+ \infty} \sum_{k=0}^{2^{j}-1} \beta_{j,k}^{2}.
$$
Meyer wavelets can be used to efficiently compute the coefficients $c_{j,k}$ and $\beta_{j,k}$ by using the Fourier transform. Indeed,  let $e_{\ell}(x) = e^{2\pi i \ell x }, \ell \in \ZZ$ and denote by $f_{\ell} = \langle f,e_{\ell}\rangle = \int_{0}^{1}f(u)e^{-2\pi i \ell u }du$ the Fourier coefficients of $f$ supposed to be a function in $L^{2}([0,1])$. By the Plancherel's identity, we obtain that
	\begin{equation}
	\beta_{j,k} = \langle f,\psi_{j,k} \rangle = \sum_{\ell \in C_{j}}  \psi^{j,k}_{\ell}  f_{\ell}, \label{eq:plancherel}
	\end{equation}
where $\psi^{j,k}_{\ell} =  \langle \psi_{j,k},e_{\ell}\rangle $ denote the Fourier coefficients of $\psi_{j,k}$ and $C_{j} = \{\ell \in \ZZ;  \psi^{j,k}_{\ell} \neq 0 \}$. As Meyer wavelets $\psi_{j,k}$ are band-limited $C_{j}$ is a finite subset set of $[-2^{j+2}c_{0},-2^{j}c_{0}] \cup [2^{j}c_{0},2^{j+2}c_{0}]$ with $c_{0} = 2\pi/3$ (see \citeA{JKPR04jrssb}).

\subsection{The case of direct density estimation}

Based on a sample $X_{1},\ldots,X_{n}$, an unbiased estimator of $f_{\ell}$ is given by $ \frac{1}{n} \sum_{m=1}^{n} \exp (-2\pi i \ell X_{m}) $ which yields an unbiased estimator of  $\beta_{j,k}$ given by
\begin{equation}
\label{wav:coefdirect}
\hat{\beta}_{j,k} = \sum_{\ell  \in C_{j}}   \psi^{j,k}_{\ell}  \left( \frac{1}{n} \sum_{m=1}^{n} e^{-2\pi i \ell X_{m}} \right).
\end{equation}
%Remark that $ \frac{1}{n} \sum_{m=1}^{n}e^{-2\pi i \ell X_{m}}$ simply corresponds to the discrete Fourier transform of the observations. 
The coefficients (\ref{wav:coefdirect})  can therefore be easily calculated by combining the fast Fourier transformation of the data with the fast algorithm of  \citeA{K94} which relies on the fact that the first sum in  (\ref{wav:coefdirect}) only involves a finite number of terms.  Equation (\ref{wav:coefdirect}) also shows that using Meyer wavelets, which are band-limited functions, avoids the use of numerical schemes to approximate the computation of the coefficients (\ref{eq:hatcoefintro})  as one typically only knows the values of the functions $\phi_{j_{0},k}$ and $ \psi_{j,k}$ at dyadic points. We define the estimators of the scaling coefficients $c_{j_{0},k}$ analogously, with $\phi$ instead of $\psi$ and $C_{j_{0}} = \{\ell \in \ZZ;  \phi^{j_{0},k}_{\ell} \neq 0 \}$ instead of $C_{j}$. 

\subsection{The case of density deconvolution}

Consider now the problem (\ref{eq:densdec}) of density deconvolution. Denote by $f_{\ell}=  \langle f,e_{\ell}\rangle $ the Fourier coefficients of $f$, and by $h_{\ell}=  \int_{\RR} h(u)e_{\ell}(u) du $ the Fourier coefficients of the error density $h$.  Since $f_{\ell} = \EE(e^{-2\pi i \ell X_{1} })$ and $h_{\ell} = \EE(e^{-2\pi i \ell \epsilon_{1} })$, it follows by independence that $\EE(e^{-2\pi i \ell Y_{1} }) = f_{\ell} h_{\ell}$. This equality and the Plancherel's identity (\ref{eq:plancherel}) therefore implies that an unbiased estimator of $\beta_{j,k}$ is given by (assuming that $h_{\ell}  \neq 0$ for all $\ell \in \ZZ$)
\begin{equation}
\label{wav:coefdec}
\hat{\beta}_{j,k} = \sum_{\ell \in C_{j}}   \tilde{\psi}^{j,k}_{\ell}  \left( \frac{1}{n} \sum_{m=1}^{n} e^{-2\pi i \ell Y_{m}} \right) \mbox{ where }   \tilde{\psi}^{j,k}_{\ell}  = \frac{  \psi^{j,k}_{\ell}}{h_{\ell}}.
\end{equation}
%Again, remark that $ \frac{1}{n} \sum_{m=1}^{n}e^{-2\pi i \ell Y_{m}}$ simply corresponds to the discrete Fourier transform of the observations, and thus the coefficients (\ref{wav:coefdec})  can be easily computed. 
%Note that a similar estimation for Meyer wavelets coefficients has been proposed in  \citeA{BVB07}. However, \citeA{BVB07} focus on log-density deconvolution to ensure the positivity of the estimators, and their results are only asymptotic with deterministic thresholds which makes the procedure proposed in this paper significantly different.

It is well-known that the difficulty of the deconvolution problem is quantified by the smoothness of the error density $h$. The so-called ill-posedness of such inverse problems depends on how fast the Fourier coefficients $h_\ell$ tend to zero. Depending on the decay of these coefficients, the estimation of $f_\ell$ will be more or less accurate. In this paper, we consider the case where the  $h_\ell$'s have a polynomial decay which is usually referred to as ordinary smooth convolution (see e.g. \citeA{F91aos}):

\begin{ass}\label{assordi}
	The Fourier coefficients of $h$ have a polynomial decay which means that  there exists a real $\nu \geqslant 0$  and two positive constants $C_{\min}, C_{\max}$ such that for all $\ell \in \ZZ$, $C_{\min} |\ell|^{-\nu} \leq  |h_{\ell}| \leq C_{\max} |\ell|^{-\nu}$. 
\end{ass}
The rates of convergence that can be expected from a wavelet estimator depend on such smoothness assumptions and are well-studied in the literature, and we refer to \citeA{PV99aos,FK02ieee} for further details.

Finally note that to simplify the presentation, we prefer to use the same notation $\hat{\beta}_{j,k}$ and  $\hat{c}_{j_{0},k}$ for both problems of  direct density estimation and density deconvolution.

\subsection{Thresholding of the empirical wavelet coefficients}

Based on an estimation of the scaling and wavelet coefficients, a linear wavelet estimator of $f$ is of the form
$
\hat{f}_{L} = \sum_{k =0}^{2^{j_{0}}-1} \hat{c}_{j_{0},k} \phi_{j_{0},k} 
	+ \sum_{j = j_{0}}^{j_{1}}  \sum_{k=0}^{2^{j}-1} \hat{\beta}_{j,k} \psi_{j,k}.
$
For an appropriate choice of $j_{1}$ one can show that $\hat{f}_{L}$ achieves optimal rates of convergence among the class of linear estimators. Typically, if $f$ belongs to a Sobolev space $H^{s}$ with smoothness order $s$, then  for direct density estimation the choice $2^{j_{1}} \approx n^{\frac{1}{2s+1}}$ yields optimal rates of convergence for the quadratic risk, see \citeA{DJKP96aos},  \citeA{judlam} for further details . However, this choice is not adaptive because it depends on the unknown smoothness $s$ of $f$. It is well known that adaptivity can be obtained  by using nonlinear estimators based on appropriate thresholding of the estimated wavelet coefficients. A non-linear estimator by hard-thresholding is defined by 
\begin{equation} \label{eq:hatfnonlin}
\hat{f}_{n} = \sum_{k =0}^{2^{j_{0}}-1} \hat{c}_{j_{0},k} \phi_{j_{0},k} + \sum_{j = j_{0}}^{j_{1}}  \sum_{k=0}^{2^{j}-1}   \hat{\beta}_{j,k}  \1_{\{|  \hat{\beta}_{j,k}  | \geqslant \tau_{j,k}\}} \psi_{j,k}
\end{equation}
where the $\tau_{j,k}$'s are appropriate thresholds (positive numbers). Various choices for $j_{1}$ and the threshold $\tau_{j,k}$ have been proposed. In the case of direct density estimation, \citeA{DJKP96aos} recommended to take level-dependent threshold $\tau_{j,k} \sim \sqrt{j/n}$ and $2^{j_{1}} \sim \frac{n}{\log(n)}$. For density deconvolution, one possible calibration in ordinary smooth deconvolution is $2^{j_1} \sim n^{\frac{1}{2\nu+1}}$ and $\tau_{j,k} \sim \frac{2^{\nu j}}{\sqrt{n}}$, see \citeA{PV99aos}. The choices $\tau_{j,k} \sim 2^{\nu j}\sqrt{j/n}$ and  $\tau_{j,k} \sim 2^{\nu j}\sqrt{ 2 \frac{\log(n)}{n} }$  have also been considered in \citeA{FK02ieee} and \citeA{BVB07} respectively.

However, such thresholds may be too large in practice as they do not take into account the fact that the variance of each empirical wavelet coefficient $\hat{\beta}_{j,k}$ depends on its location $(j,k)$. Consider the problem of density deconvolution and let us denote the variance of $\hat{\beta}_{j,k}$ by
$$
\sigma^{2}_{j,k} = \EE(\hat{\beta}_{j,k} - \beta_{j,k})^{2}.
$$
Let us also denote by $\tilde{\psi}_{j,k}$ the function defined for $y \in \RR$ by
$$
\tilde{\psi}_{j,k}(y) = \sum_{\ell \in C_{j}} \tilde{\psi}^{j,k}_{\ell} e^{-2\pi i \ell y }.
$$
By definition, it follows that $\hat{\beta}_{j,k} = \frac{1}{n}  \sum_{m=1}^{n}  \tilde{\psi}_{j,k}(Y_{m})$ and thus $\sigma^{2}_{j,k} =  \frac{1}{n} \var(  \tilde{\psi}_{j,k}(Y_{1}))$. Hence, a simple upper bound for $\sigma^{2}_{j,k}$ is $V_{jk} = \frac{1}{n}  \EE(\tilde{\psi}_{j,k}(Y_{1})^{2})$. Then, simple algebra shows that in the case of density deconvolution
\begin{equation}
V_{j,k} =  \frac{1}{n} \int_{\RR} |\tilde{\psi}_{j,k}(y)|^{2} f^{Y}(y)dy =  \frac{1}{n} \sum_{\ell, \ell' \in C_{j}}  \tilde{\psi}^{j,k}_{\ell} \overline{ \tilde{\psi}^{j,k}_{\ell'} } f^{Y}_{\ell-\ell'}, \quad \label{eq:Vdens}
\end{equation}
with  $f^{Y}(y) = \int_{0}^{1}f(u)h(y-u)du$ and  $f^{Y}_{\ell-\ell'} = \EE  e^{-2\pi i (\ell-\ell') Y_{1} }$. An unbiased estimator of $V_{j,k}$ is thus given by
\begin{equation}
\hat{V}_{j,k} =  \frac{1}{n} \sum_{\ell, \ell' \in C_{j}}  \tilde{\psi}^{j,k}_{\ell} \overline{ \tilde{\psi}^{j,k}_{\ell'} } \left( \frac{1}{n} \sum_{m=1}^{n} e^{-2\pi i (\ell-\ell') Y_{m}}\right) =  \frac{1}{n^{2}}  \sum_{m=1}^{n} \left|\sum_{\ell  \in C_{j}} \tilde{\psi}^{j,k}_{\ell}  e^{-2\pi i \ell Y_{m}}  \right|^{2}   . \label{eq:hatVdens}
\end{equation}
Similar computations can be made for the case of direct density estimation with $\psi^{j,k}_{\ell}$ instead $\tilde{\psi}^{j,k}_{\ell}$, $f_{\ell-\ell'}$ instead of $f^{Y}_{\ell-\ell'}$, and $X_{m}$ instead of $Y_{m}$ in the above equations (\ref{eq:Vdens}) and (\ref{eq:hatVdens}). Note that $\hat{V}_{j,k}$ can also be written as $\hat{V}_{j,k} =    \frac{1}{n^{2}}  \sum_{m=1}^{n} \left| \tilde{\psi}_{j,k}(Y_{m}) \right|^{2}$, but its calculation is obtained from the Fourier coefficients  $(\tilde{\psi}^{j,k}_{\ell} )_{\ell \in C_{j}}$ and not from $ \tilde{\psi}_{j,k}$ whose computation at non dyadic points requires numerical approximation. Alternatively, one could also use an estimation of the variance $\sigma^{2}_{j,k}$  given by
$$
\hat{\sigma}^{2}_{j,k} = \hat{V}_{j,k} - \frac{1}{n} |\hat{\beta}_{j,k}|^{2},
$$
instead of the upper bound $V_{j,k}$. However, this does not change significantly our results since $\hat{\sigma}^{2}_{j,k}$ and $\hat{V}_{j,k}$ are very close  for $n$ sufficiently large. Moreover, oracle inequalities are simpler to derive using an estimated upper bound for  $\sigma^{2}_{j,k}$ .

A thresholding rule is usually chosen by controlling the probability of deviation of $\hat{\beta}_{j,k}$ from the true wavelet coefficient  $\beta_{j,k}$. From Lemma \ref{lemma:ineg1} (see the Appendix) one has that for any positive $x$,
$
\PP \left( |\hat{\beta}_{j,k} - \beta_{j,k} | \geq \sqrt{2 V_{jk} x} + \frac{ \eta_{j} }{3n} x \right) \leq 2 \exp(-x)
$
where
\begin{equation} \label{eq:eta}
\eta_{j} = \sum_{\ell \in C_{j}}  |\tilde{\psi}^{j,k}_{\ell}|,
\end{equation}
which would suggest to take a threshold of the form
$$
\tau_{j,k}^{\ast} = \sqrt{2 \delta \log(n) V_{j,k}  } + \frac{\delta \log(n)}{3n} \eta_{j}, 
$$
where $\delta > 0$ is a tuning parameter. Thinking of the classical universal threshold, one would take $\delta = 1$. However, this choice can be too conservative and the results of this paper shows that it is possible to take $\delta$ smaller that 1. Moreover, it is shown that the choice of this tuning parameter depends on the highest resolution level $j_{1}$ and the degree of ill-posedness $\nu$ in the case of density deconvolution. Throughout the paper, we discuss the choice of $\delta$ and finally propose data-based values for its calibration.
 
Obviously $\tau_{j,k}^{\ast}$ is an ideal threshold as $ V_{j,k} $ is unknown. Based on Lemma \ref{lemma:ineg2}  (see the Appendix) which gives a control on the probability of deviation of $\hat{V}_{j,k}$ from  $ V_{j,k} $, we propose to use the following random thresholds
\begin{equation} \label{eq:thr}
\tau_{j,k} = \sqrt{2 \delta \log(n)\left( \hat{V}_{j,k}  + \sqrt{2 \delta \log(n) \hat{V}_{j,k}  \frac{\eta_{j}^{2}}{n^{2}}   } + \delta \log(n) \kappa \frac{\eta_{j}^{2}}{n^{2}} \right) } + \frac{\delta \log(n)}{3n} \eta_{j} ,
\end{equation}
where $\kappa = \frac{4}{3} + \sqrt{\frac{5}{3}}$

Again, for direct density estimation, we take the same thresholds with $\psi^{j,k}_{\ell}$ instead of $\tilde{\psi}^{j,k}_{\ell}$ to compute $\eta_{j}$ in (\ref{eq:eta}). The above choice for $\tau_{j,k}$ resembles to the universal threshold $\sigma \sqrt{2 \log(n) }$ proposed by \citeA{donojohn} in the context of nonparametric regression with homoscedastic variance $\sigma^{2}$. Here we exploit the fact that in the context of density estimation the variance of a wavelet coefficient depends on its location $(j,k)$ and has to be estimated.  

The additive terms  $\frac{\delta \log(n)}{3n} \eta_{j}$ and  $ \sqrt{2 \delta \log(n) \hat{V}_{j,k}  \frac{\eta_{j}^{2}}{n^{2}}   } + \delta \log(n) \kappa \frac{\eta_{j}^{2}}{n^{2}}$ will allow us to derive oracle inequalities. Indeed, in Section \ref{sec:oracle}, we compare the quadratic risk of $\hat{f}_{n}$ to the risk of the following oracle estimator
\begin{equation} \label{eq:oracle}
\tilde{f}_{n} = \sum_{k =0}^{2^{j_{0}}-1} \hat{c}_{j_{0},k} \phi_{j_{0},k} + \sum_{j = j_{0}}^{j_{1}}  \sum_{k=0}^{2^{j}-1}   \hat{\beta}_{j,k}  \1_{\{|\beta_{j,k}|^{2}  \geqslant \sigma^{2}_{j,k} \}} \psi_{j,k}
\end{equation}
Note that $\tilde{f}_{n}$ is an ideal estimator that can not be computed in practice as its depends on the unknown coefficients $\beta_{j,k}$ of $f$ and the unknown variance terms $\sigma^{2}_{jk}$. However, we shall use it as a benchmark to assess the quality of our estimator.  The quadratic risk of $\tilde{f}_{n}$ is
\begin{equation} \label{eq:riskoracle}
\EE \|\tilde{f}_{n}-f\|^{2} = \sum_{k =0}^{2^{j_{0}}-1} \sigma^{2}_{j_{0},k} + \sum_{j = j_{0}}^{j_{1}}  \sum_{k=0}^{2^{j}-1}   \min(\beta_{j,k}^{2} ,\sigma^{2}_{j,k})  + \sum_{j = j_{1} +1}^{+ \infty}   \sum_{k=0}^{2^{j}-1}\beta_{j,k}^{2}  ,
\end{equation}
where $\sigma^{2}_{j_{0},k} =  \EE(\hat{c}_{j_{0},k} - c_{j_{0},k})^{2}$. Equation (\ref{eq:riskoracle}) shows that we retrieve the classical formula for the quadratic risk of an oracle estimator given in \citeA{donojohn} except that the variance term $\sigma^{2}_{j,k}$ is not constant as in standard nonparametric regression with homoscedastic variance.

Data-driven thresholds based on an estimation of the variance of the empirical wavelet coefficients have already been proposed in \citeA{judlam} in the context of density estimation. However, the estimation procedure in  \citeA{judlam} is different from ours since it is  based on biorthogonal bases and on the use of the Haar basis to compute the wavelet coefficients. Note that our choice for $\tau_{j,k}$ is similar to the random threshold in  \citeA{judlam}, but the addition of the terms $ \sqrt{2 \delta \log(n) \hat{V}_{j,k}  \frac{\eta_{j}^{2}}{n^{2}}   } + \delta \log(n) \kappa \frac{\eta_{j}^{2}}{n^{2}}$ and  $\frac{\delta \log(n)}{3n} \eta_{j}$ will allow us to compare the performances of $\hat{f}_{n}$ with those of the oracle $\tilde{f}_{n}$. The addition of similar deterministic and stochastic terms is also proposed in \citeA{reyriv} to derive  oracle inequalities in the context of Poisson intensity estimation.

\section{Oracle inequalities} \label{sec:oracle}

To derive oracle inequalities, we need a further smoothness assumption on the error density $h$:
\begin{ass} \label{ass:h}
There exists a constant $C > 0$ and a real $\rho > 1$ such that the density $h$ satisfies $h(x) \leq \frac{C}{1 + |x|^{\rho}}$ for all $x \in \RR$.
\end{ass}
Obviously, the above condition for $h$ is not very restrictive as $h$ is by definition an integrable function on $\RR$. For a bounded function $f \in L^{2}([0,1])$ we denote by $\| f \|_{\infty} = \sup_{x \in [0,1]}\{ |f(x)|\}$ its supremum norm. Let us also define the following class of densities
$$
D^{2}([0,1]) = \{ f \in  L^{2}([0,1]), \mbox{ with } f \geq 0 \mbox{ and } \int_{0}^{1} f(x) dx = 1 \}.
$$

\subsection{The case of direct density estimation}

The following theorem states that for appropriate choices of $j_{1}, j_{0}$ and the tuning parameter $\delta$, then the estimator $\hat{f}_{n}$ behaves essentially as the oracle $\tilde{f}_{n}$ up to logarithmic terms.

\begin{theo} \label{theo:oracledirect}
Assume that $f \in D^{2}([0,1])$ with $\| f \|_{\infty} < + \infty$ . Let $\alpha \geq 0$ and $1/2 \geq \eta > 0$ be some fixed constants. For any $n \geq \exp(1)$, define $j_{1} = j_{1}(n)$ to be the integer such that $2^{j_{1}} > n^{\eta} (\log n)^{\alpha} \geq 2^{j_{1}-1}$, and $j_{0} = j_{0}(n)$ to be the integer such that $2^{j_{0}} > \log(n) \geq 2^{j_{0}-1}$ and  suppose that $\eta$ and $\alpha$ are such that $j_{1} \geq j_{0}$. Assume that $\delta > \eta$, and take the random thresholds $\tau_{j,k}$ given by equation (\ref{eq:thr}). Then, the estimator $\hat{f}_{n}$ satisfies the following oracle inequality
\begin{equation} \label{eq:ineqoracle}
\EE \|\hat{f}_{n}-f\|^{2} \leq C_{1}(\delta) \left[ \sum_{k =0}^{2^{j_{0}}-1} \sigma^{2}_{j_{0},k} + \sum_{j = j_{0}}^{j_{1}}  \sum_{k=0}^{2^{j}-1}   \min(\beta_{j,k}^{2} ,\log(n) \sigma^{2}_{j,k})  + \sum_{j = j_{1} +1}^{+ \infty}   \sum_{k=0}^{2^{j}-1}\beta_{j,k}^{2} \right] + C_{2}(\delta) \Gamma_{n,1},
\end{equation}
where
$$
\Gamma_{n,1} =  \frac{\log(n)}{n} \sum_{j = j_{0}}^{j_{1}}  \sum_{k=0}^{2^{j}-1} \beta_{j,k}^{2}   + \max(\|f\|_{\infty},1)  \max((\log n)^{\alpha},1)  \frac{(\log n)^{\alpha}}{n} +  \frac{(\log n)^{2+ 2\alpha}}{n^{2(1-\eta)}}
$$
and $C_{1}(\delta)$ and $C_{2}(\delta)$ are two positive constants not depending on $n$ and $f$, and such that $\lim_{\delta \to \eta} C_{1}(\delta) = \lim_{\delta \to \eta} C_{2}(\delta) = + \infty$.
\end{theo}

The above inequality (\ref{eq:ineqoracle}) shows that the performances of  $\hat{f}_{n}$ mimic those of the oracle $\tilde{f}_{n}$ in term of quadratic risk, see equation (\ref{eq:riskoracle}), up to a logarithmic term. The additive term $\Gamma_{n,1}$ depends on two hyperparameters $\alpha$ and $\eta$ which are used to control the effect of the choice of the highest resolution level $j_{1}$ and the tuning parameter $\delta$ on the performances of the estimator. Moreover, the above inequality tends to show that the performances of $\hat{f}_{n}$ deteriorates as $\delta$ tends to $\eta$. Thinking of the classical universal threshold, one would like to take $\delta = 1$. However, if one sets $\eta = 1/2$ and $\alpha = 0$, the additive term $\Gamma_{n,1}$ is bounded by $\frac{(\log n)^{2}}{n}$ which is rate typically faster that the decay of the oracle risk (\ref{eq:riskoracle}) when $f$ belongs to a Sobolev or a Besov space. Hence, Theorem \ref{theo:oracledirect} shows that if one chooses $j_{1} = \lfloor \eta \log_{2}(n) \rfloor + 1$ with $\eta \leq 1/2$, then it is possible to take $\delta$ smaller than 1. Choosing carefully such hyperparameters is of fundamental importance, and a detailed study is therefore proposed in Section \ref{sec:simu} to validate the results of Theorem \ref{theo:oracledirect}, and to analyze the risk of $\hat{f}_{n}$ as a function of the resolution level $j_{1}$ and the tuning constant $\delta$.

Classically, the level $j_{1}$ is chosen as the integer $j_{1}$ such that $2^{j_{1}} \geq \frac{n}{\log(n)} \geq 2^{j_{1}-1}$ (see e.g. \citeA{DJKP96aos}). The results of this paper shows that one can use smaller level $2^{j_{1}}$ of the order $n^{\eta} (\log n)^{\alpha}$.  For $\eta \leq 1/2$, the price to pay is a slightly lower rate for the additive term $\Gamma_{n,1}$ in the oracle inequality (\ref{eq:ineqoracle}) which is of the order $\frac{(\log n)^{2}}{n}$ instead of the rate $\frac{1}{n}$ as classically obtained for the additive term when deriving oracle inequalities. Note that a similar result in the context of Poisson intensity estimation is also given in \citeA{reyriv}.   In particular \citeA{reyriv} obtain an additive term of the order $\frac{1}{n}$ but to derive such a result their proof relies heavily on the fact that the wavelet basis they use (the Haar basis) is such that $\inf_{x \in [0,1]} |\psi(x)| > 0$ which is not the case for Meyer wavelets.

\subsection{The case of density deconvolution}

Consider now the problem of density deconvolution under Assumption \ref{assordi} of ordinary smooth deconvolution.

\begin{theo} \label{theo:oracledec}
Assume that $f \in D^{2}([0,1])$ with $\|f\|_{\infty} < + \infty$, and that $h$ satisfies  Assumption \ref{assordi}  and Assumption \ref{ass:h}. Let $\alpha \geq 0$ and $1/2 \geq \eta > 0$ be some fixed constants. For any $n > \exp(1)$, define $j_{1} = j_{1}(n)$ to be the integer such that $2^{j_{1}} > n^{\eta/(\nu +1)} (\log n)^{\alpha} \geq 2^{j_{1}-1}$, and $j_{0} = j_{0}(n)$ to be the integer such that $2^{j_{0}} > \log(n) \geq 2^{j_{0}-1}$,  and  suppose that $\eta$ and $\alpha$ are such that $j_{1} \geq j_{0}$.  Assume that $\delta > \eta \left(1 + \frac{\nu}{\nu +1}\right)$, and take the random thresholds $\tau_{j,k}$ given by equation (\ref{eq:thr}). Then, the estimator $\hat{f}_{n}$ satisfies the following oracle inequality
\begin{equation} \label{eq:ineqoracledec}
\EE \|\hat{f}_{n}-f\|^{2} \leq C_{3}(\delta) \left[ \sum_{k =0}^{2^{j_{0}}-1} \sigma^{2}_{j_{0},k} + \sum_{j = j_{0}}^{j_{1}}  \sum_{k=0}^{2^{j}-1}   \min(\beta_{j,k}^{2} ,\log(n) \sigma^{2}_{j,k})  + \sum_{j = j_{1} +1}^{+ \infty}   \sum_{k=0}^{2^{j}-1}\beta_{j,k}^{2} \right] + C_{4}(\delta) \Gamma_{n,2},
\end{equation}
where
$$
\Gamma_{n,2} =  \frac{\log(n)}{n} \sum_{j = j_{0}}^{j_{1}}  \sum_{k=0}^{2^{j}-1} \beta_{j,k}^{2}   + \max(\|f\|_{\infty},1) \max((\log n)^{\alpha},1) \frac{(\log n)^{\alpha (2 \nu +1)}}{n} +  \frac{(\log n)^{2+ 2\alpha +  2 \alpha \nu}}{n^{2(1-\eta)}}
$$
and $C_{3}(\delta)$ and $C_{4}(\delta)$ are two positive constants not depending on $n$ and $f$, and such that $\lim_{\delta \to \eta \left(1 + \frac{\nu}{\nu +1}\right)} C_{3}(\delta) = \lim_{\delta \to \eta \left(1 + \frac{\nu}{\nu +1}\right)} C_{4}(\delta) = + \infty$.
\end{theo}

Hence, the above theorem shows that in the case of density deconvolution then the estimator also behaves as an oracle estimate  up to logarithmic terms, and that the performances of the estimator tend to deteriorate as $\delta$  tends to $\eta  \left(1 + \frac{\nu}{\nu +1}\right)$. Similar comments to those given for Theorem  \ref{theo:oracledirect} can be made. If one chooses $\eta = 1/2$, then the additive term $\Gamma_{n,2} $ is of the order $\frac{(\log n)^{2+ 2\alpha +  2 \alpha \nu}}{n}$, and $\delta$ has to be greater  than $\left(1/2 + \frac{\nu}{2\nu +2}\right)$. Hence, this again shows that one can take a value for $\delta$ smaller than 1. However the choice of $\delta$ is typically larger for deconvolution than in the direct case, as it is controlled by the degree $\nu$ of ill-posedness.

 In deconvolution problems, the high-frequency cut-off  $j_{1}$ is usually related to the ill-possedness $\nu$ of the inverse problem and $2^{j_{1}} = \opO\left( (\frac{n}{\log(n)})^{1/(2\nu+1)}\right)$ is a typical choice for various estimators proposed in the literature (see e.g. \citeA{JKPR04jrssb,sapatinaspensky,BVB07}).  This is a standard fact that a smaller $j_{1}$ should be used for ill-posed inverse problems than in the direct case. Again we have introduced  hyperparameters  $\alpha$ and $\eta$ to control the effect of the choice of the highest resolution level $j_{1}$. Understanding the influence of the choice of these hyperparameters on the quality of the estimator is a fundamental issue, and a detailed simulation study is thus proposed in Section \ref{sec:simu} to validate the results of Theorem \ref{theo:oracledec}.

\section{Asymptotic properties and near-minimax optimality} \label{sec:minimax}

It is well known that Besov spaces for periodic functions in $L^{2}([0,1])$  can be characterized in terms of wavelet coefficients (see e.g. \citeA{JKPR04jrssb}). Let $s > 0$ denote the usual smoothness parameter, then  for the Meyer wavelet basis and for a Besov ball $B^{s}_{p,q}(A)$ of radius $A > 0$ with $1 \leq p,q \leq \infty$, one has that  for $s + 1/2-1/p \geq 0$
$$
B^{s}_{p,q}(A) = \left\{ f  \in L^{2}([0,1]) :  \left(\sum_{k =0}^{2^{j_{0}}-1} |c_{j_{0},k}|^{p} \right)^{\frac{1}{p}} + \left( \sum_{j = j_{0}}^{+ \infty} 2^{j(s + 1/2-1/p)q} \left( \sum_{k=0}^{2^{j}-1} |\beta_{j,k}|^{p}\right)^{\frac{q}{p}} \right)^{\frac{1}{q}} \leq A \right\}
$$ 
with the respective above sums replaced by maximum if $p=\infty$ or $q=\infty$. 

The condition that $s + 1/2-1/p \geq 0$ is imposed to ensure that $B^s_{p,q}(A)$ is a subspace of $L^2([0,1])$, and we shall restrict ourselves to this case in this paper. Besov spaces allow for more local variability in local smoothness than is typical for functions in the usual H\"older or Sobolev spaces. For instance, a real function $f$ on $[0,1]$ that is piecewise continuous, but for which each piece is locally in $C^s$, can be an element of $B^s_{p,p}(A)$ with $1 \leq p < 2$, despite the possibility of discontinuities at the transition from one piece to the next. Note that if $s \geq 1$ is not an integer, then $B_{2,2}^{s}(A)$ is equivalent to a Sobolev ball of order $s$, and that the space $B_{p,q}^{s}(A)$ with $1 \leq p < 2$ contains piecewise smooth functions with local irregularities such as discontinuities.Finally let us introduce the following space of densities
$$
D^{s}_{p,q}(A) = D^{2}([0,1]) \cap B^{s}_{p,q}(A). 
$$
The following theorems show that for either the problem of direct density estimation or density deconvolution, the estimator $\hat{f}_{n}$ is asymptotically near-optimal (in the minimax sense) up to a logarithmic factor over a wide range of Besov balls.

To simplify the presentation, all the asymptotic properties of $\hat{f}_{n}$ are given in the case $\eta = 1/2$ and $\alpha = 0$ where $\eta,\alpha$ are the hyperparameters introduced in Theorem  \ref{theo:oracledirect} and  Theorem \ref{theo:oracledec}.

\subsection{The case of direct density estimation}

For $s > 0$ and $1 \leq p \leq \infty$ and let us define
$$
p' = \min(2,p) \mbox{ and } s^{\ast} = s + 1/2 - 1/p'.
$$
Then the following result holds.

\begin{theo} \label{theo:minimaxdirect}
Assume that the conditions of Theorem \ref{theo:oracledirect} hold with $\eta = 1/2$ and $\alpha = 0$. Assume that $f \in D^{s}_{p,q}(A)$ with $s > \frac{1}{2} + \frac{1}{p'}$,   $1 \leq p \leq 2$ and $1 \leq q \leq 2$. Then,  as $n \to + \infty$
$$
\sup_{f \in D^{s}_{p,q}(A)} \EE  \|\hat{f}_{n}-f\|^{2} \leq \opO \left( \frac{n}{\log(n)} \right) ^{-\frac{2s}{2s+1}} .
$$
\end{theo}

Minimax rates of convergence for density estimation over Besov spaces has been studied in detail by \citeA{DJKP96aos}. Hence, Theorem \ref{theo:minimaxdirect} shows that $\hat{f}_{n}$ is an adaptive estimator which converges to $f$ with the minimax rate up to logarithmic factor for the problem of direct density estimation over $D^{s}_{p,q}(A)$. The extra logarithmic factor is usually called the price to pay for adaptivity to the unknown smoothness $s$.

%Note that it can be checked from the proofs in the Appendix that the constant $C_{3}$ depends only on $A,s,p,q$, the parameters $\alpha$ and $\delta$, and the wavelet basis.

\subsection{The case of density deconvolution}

Consider now the problem of density deconvolution under the Assumption \ref{assordi} of ordinary smooth deconvolution.

\begin{theo} \label{theo:minimaxdec}
Assume that the conditions of Theorem \ref{theo:oracledec} hold $\eta = 1/2$ and $\alpha = 0$. Assume that $f \in D^{s}_{p,q}(A)$ with $s > \frac{1}{2} + \frac{1}{p'}$,   $1 \leq p \leq 2$ and $1 \leq q \leq 2$. If $\nu(2-p) < p s^{\ast}$ then as $n \to + \infty$
$$
\sup_{f \in D^{s}_{p,q}(A)}  \EE  \|\hat{f}_{n}-f\|^{2} \leq \opO \left(   \frac{n}{\log(n)}  \right)^{-\frac{2s}{2s+2\nu+1}} ,
$$
 and if  $\nu(2-p) \geq p s^{\ast}$ then as $n \to + \infty$
$$
\sup_{f \in D^{s}_{p,q}(A)}  \EE  \|\hat{f}_{n}-f\|^{2} \leq  \opO \left(   \frac{n}{\log(n)} \right)^{-\frac{2s^{\ast}}{2s^{\ast}+2\nu}} 
$$
\end{theo}

Theorem \ref{theo:minimaxdec} show that there is two different rates of convergence depending on whether $\nu(2-p) < p s^{\ast}$ or $\nu(2-p) \geq p s^{\ast}$. These two conditions are respectively referred to as the dense case when the worst functions $f$ to estimate are spread uniformly over $[0,1]$,  or the sparse case when the hardest functions to estimate have only one non-vanishing wavelet coefficient. This change in the rate of convergence, usually referred to as an elbow effect, has been studied in detail in nonparametric deconvolution problems in the white noise model by \citeA{JKPR04jrssb} and also   
\citeA{sapatinaspensky}. Theorem \ref{theo:minimaxdec} shows that the rate of convergence of $\hat{f}_{n}$ corresponds to minimax rates (up to a logarithmic term) that have been obtained in related deconvolution problems either for density estimation or nonparametric regression in the white noise model (see e.g. \citeA{PV99aos,FK02ieee,JKPR04jrssb,sapatinaspensky}, and references therein).

%%%% Simulations %%%%%%

\section{Simulations} \label{sec:simu}

\begin{figure}[h!]
\centering
\subfigure[]
{ \includegraphics[width=3.5cm]{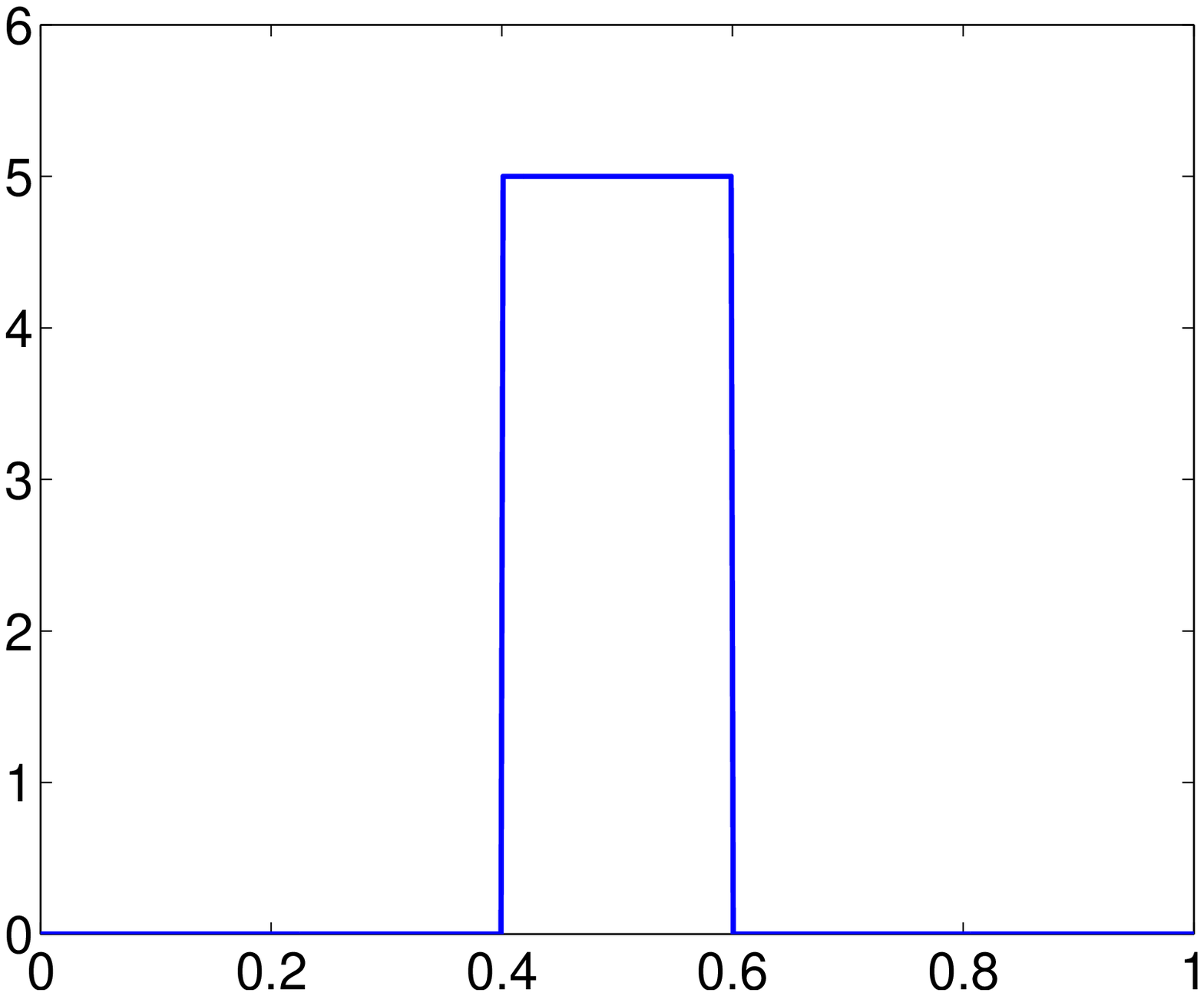} }
\subfigure[]
{ \includegraphics[width= 3.5cm]{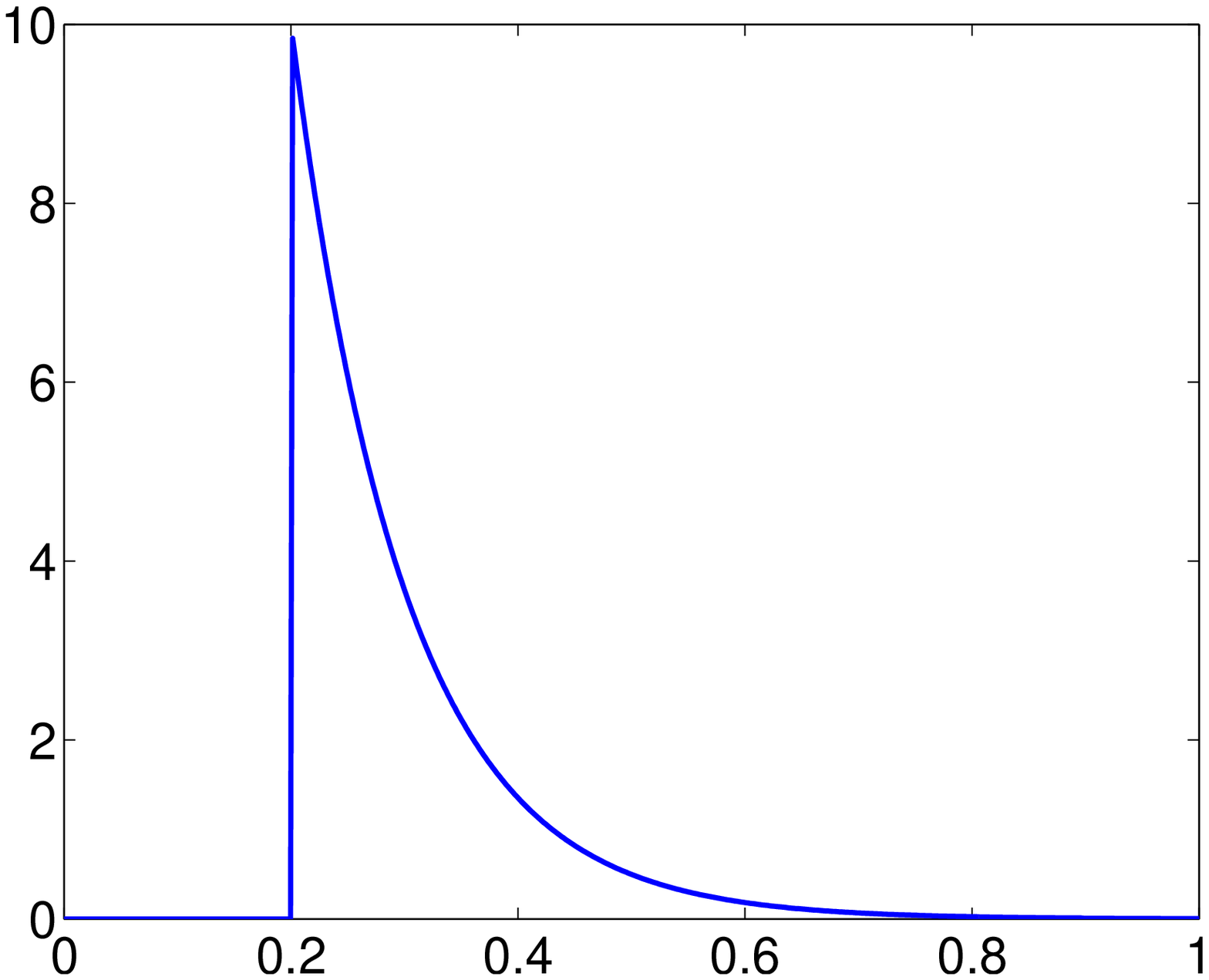} }
\subfigure[]
{ \includegraphics[width= 3.5cm]{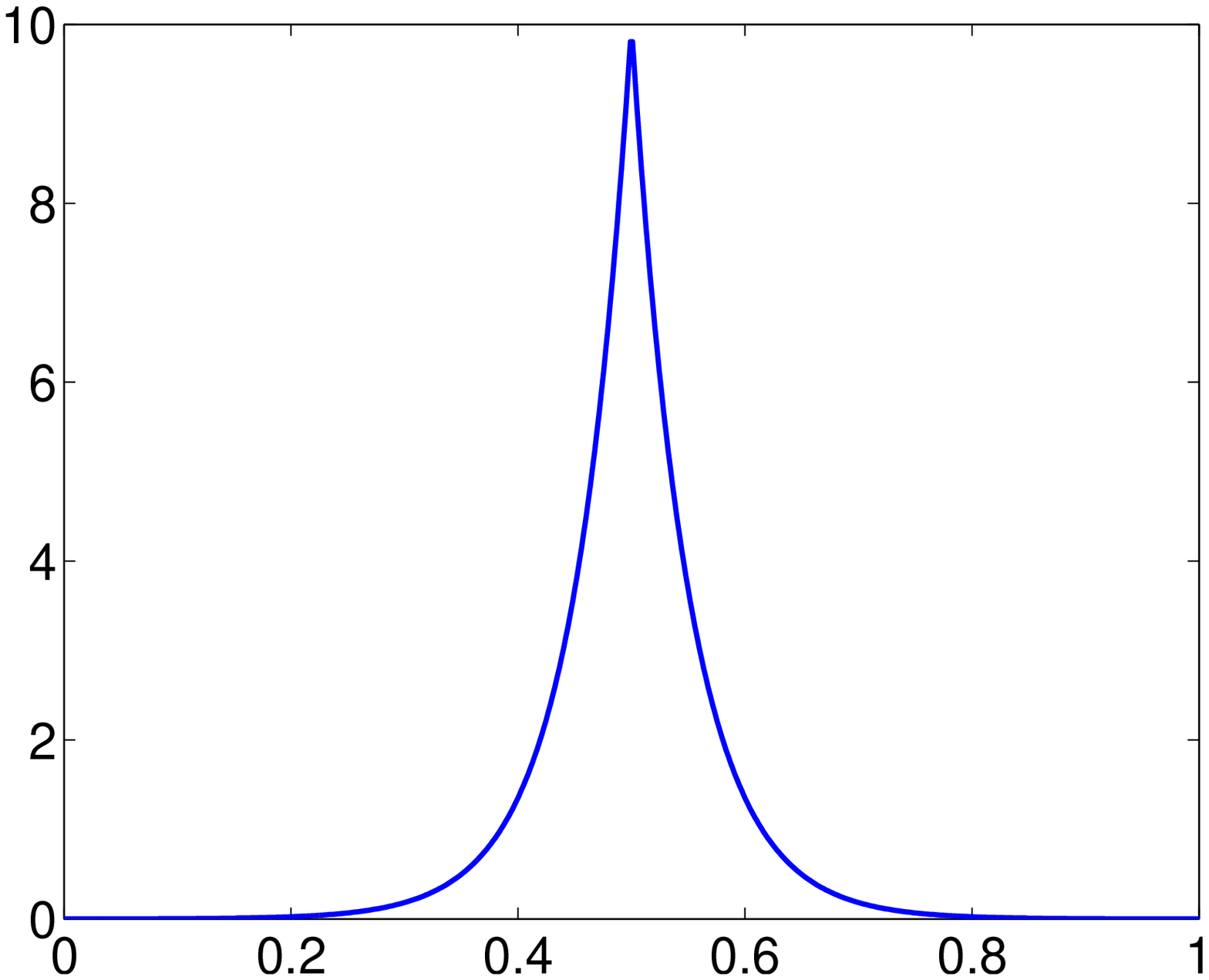} }
\subfigure[]
{ \includegraphics[width= 3.5cm]{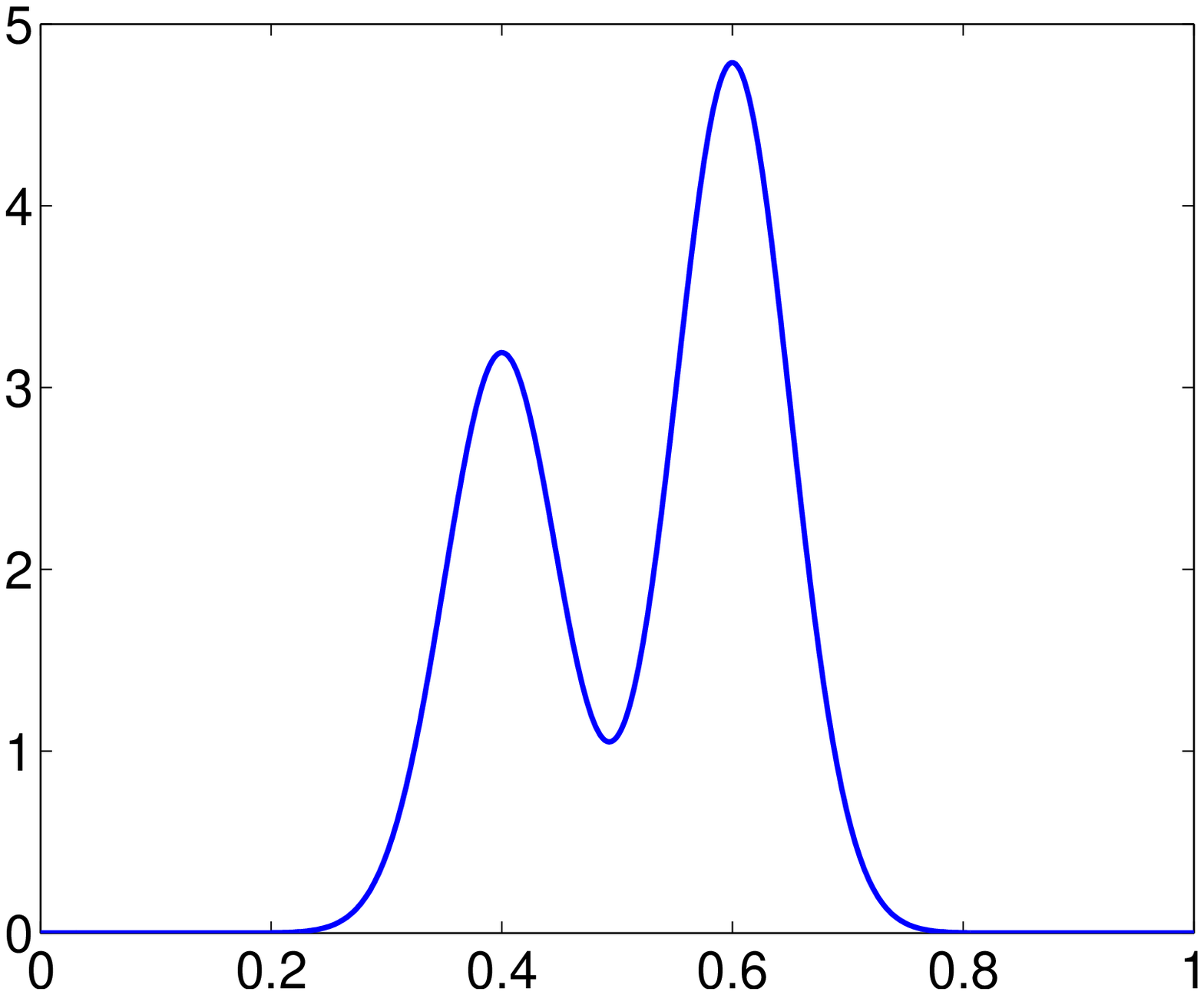} }
\caption{Test densities: (a) Uniform, (b) Exponential, (c) Laplace, (d) MixtGauss. } \label{figdens}
\end{figure}

Simulations use the wavelet toolbox \emph{Wavelab} of Matlab \shortcite{BCDJ95} and the fast algorithm for Meyer wavelet decomposition developed by  \citeA{K94}. As in the simulation study of \citeA{BVB07}, four test densities are considered:
	{\it Uniform distribution}: $f(x) = 5 \1_{[0.4,0.6]}(x)$,
	{\it Exponential distribution}: $f(x) = 10e^{-10(x-0.2)} \1_{[0.2,+\infty[}(x)$,
	{\it Laplace distribution}: $f(x) = 10e^{-20|x-0.5|}$, and
	{\it MixtGauss distribution (mixture of two Gaussian variables)}: $X \sim  \pi_{1}N(\mu_{1},\sigma_{1}^{2})+ \pi_{2}N(\mu_{2},\sigma_{2}^{2})$ with $\pi_{1} = 0.4, \pi_{1} = 0.6$, $\mu_{1} = 0.4, \mu_{2} = 0.6$ and $\sigma_{1} = \sigma_{2} = 0.05$. These four densities, displayed in Figure \ref{figdens}, exhibit different types of smoothness: the Uniform density is a piecewise constant function with two jumps, the Exponential distribution is a piecewise smooth function with a single jump, the Laplace density is a continuous function with a cusp at $x =0.5$, whereas the MixtGauss density is infinitely differentiable. 

\subsection{Direct density estimation}

Assume that an i.i.d sample of variables $X_{1},\ldots,X_{n}$ is drawn from one the test densities shown in  Figure \ref{figdens}. The empirical Fourier coefficients $\hat{f}_{\ell} = \frac{1}{n} \sum_{m=1}^{n} \exp (-2\pi i \ell X_{m})$ are computed for $\ell = -N/2+1,\ldots,N/2$, where $N = 2^J$ is a dyadic integer with $J$ chosen such that $N \geq n$.
The  $\hat{f}_{\ell}$'s are then used as an input of the efficient algorithm of \citeA{K94} in order to compute empirical scaling and wavelet coefficients $(\hat{c}_{j_{0},k})_{k = 0,\ldots,2^{j_{0}-1}},(\hat{\beta}_{j,k})_{k=0,\ldots,2^{j}-1,j=j_{0},\ldots,J-1}$. This algorithm only requires $\opO(N(\log(N))^{2})$ operations to compute the empirical wavelet coefficients from a sample of Fourier coefficients of size $N$. According to Theorem \ref{theo:oracledirect}, one can take $j_{0} = \lfloor \log_{2}(\log(n)) \rfloor + 1$. Then, two hyperparameters essentially  controls the quality of the estimator: the high-frequency cut-off parameter $j_{1}$ and the constant $\delta$ used in the definition of the thresholds  $\tau_{j,k}$ given by equation (\ref{eq:thr}). If the  level $j_{1}$ is such that  $2^{j_{1}} > n^{\eta} (\log n)^{\alpha} \geq 2^{j_{1}-1}$ for some $\eta > 0$, $\alpha \geq 0$, then one must take $\delta > \eta$. Increasing $\alpha$ deteriorates the rate of convergence of the additive term in the oracle inequality (\ref{eq:ineqoracle}), so one can choose to set $\alpha = 0$. Following the choice made for the asymptotic study of $\hat{f}_{n}$ in Section \ref{sec:minimax}, one can take $\eta = 1/2$ and according to Theorem  \ref{theo:oracledirect} one should then take $\delta > 1/2$. However, it is not clear if a value for $\delta$ lower than 1/2 would deteriorate or improve the quality of the estimator $\hat{f}_{n}$.

Alternatively, assume that $j_{1}$ is given. Then, another possibility for choosing $\delta$ is to take $\eta^{\ast} = (j_{1}-1)/\log_{2}(n)$ which is the smallest constant $\eta$ that satisfies $2^{j_{1}} > n^{\eta} \geq 2^{j_{1}-1}$, and then take $\delta >  \eta^{\ast}$. Combining the above arguments, we finally suggest to take
\begin{equation}
j_{1} = j_{1}^{\ast} = \lfloor \frac{1}{2} \log_{2}(n) \rfloor + 1 \mbox{ and } \delta^{\ast} = (j_{1}-1)/\log_{2}(n). \label{eq:choice}
\end{equation}

\begin{figure}[h!]
\centering
\subfigure[]
{ \includegraphics[width=3.5cm]{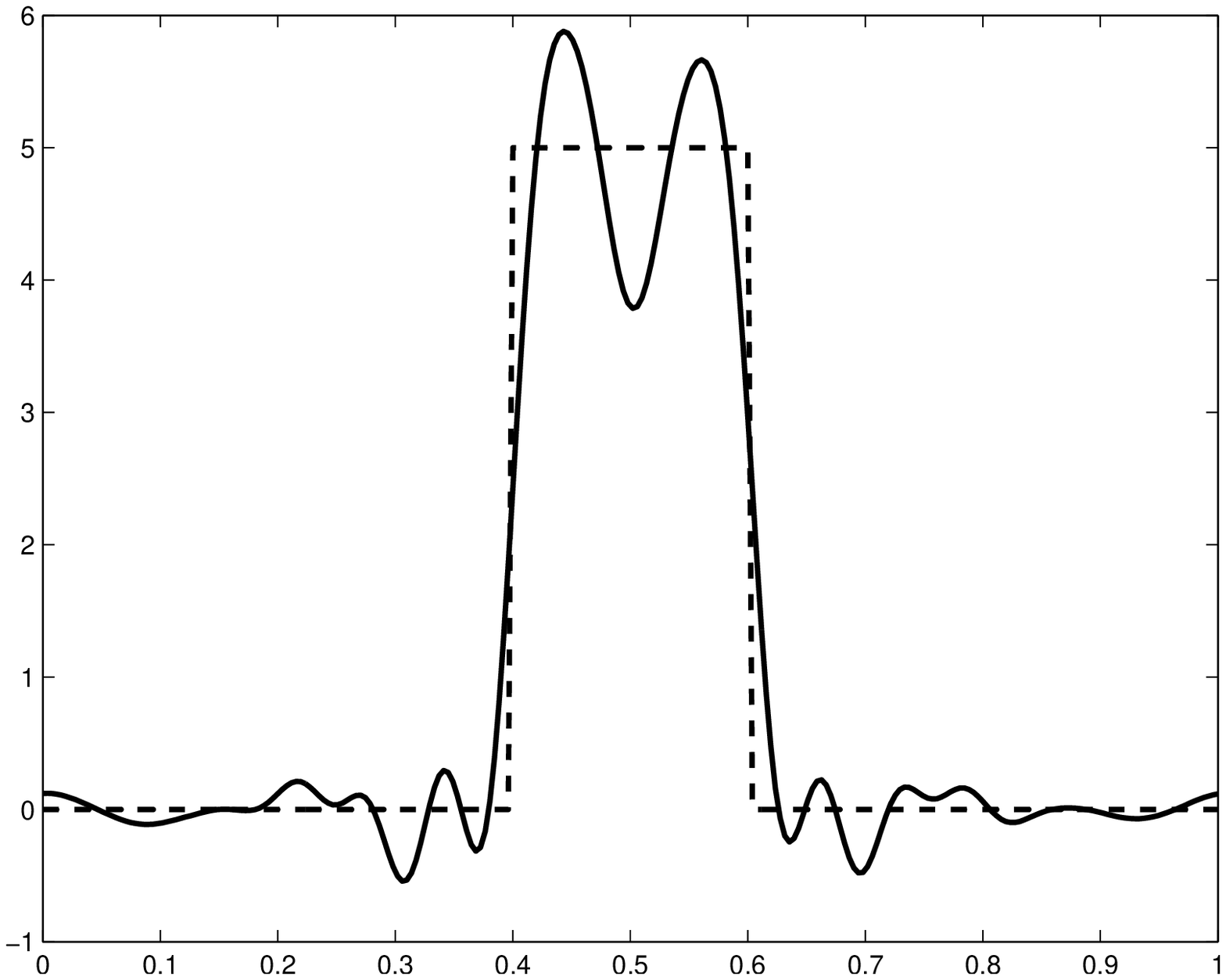} }
\subfigure[]
{ \includegraphics[width= 3.5cm]{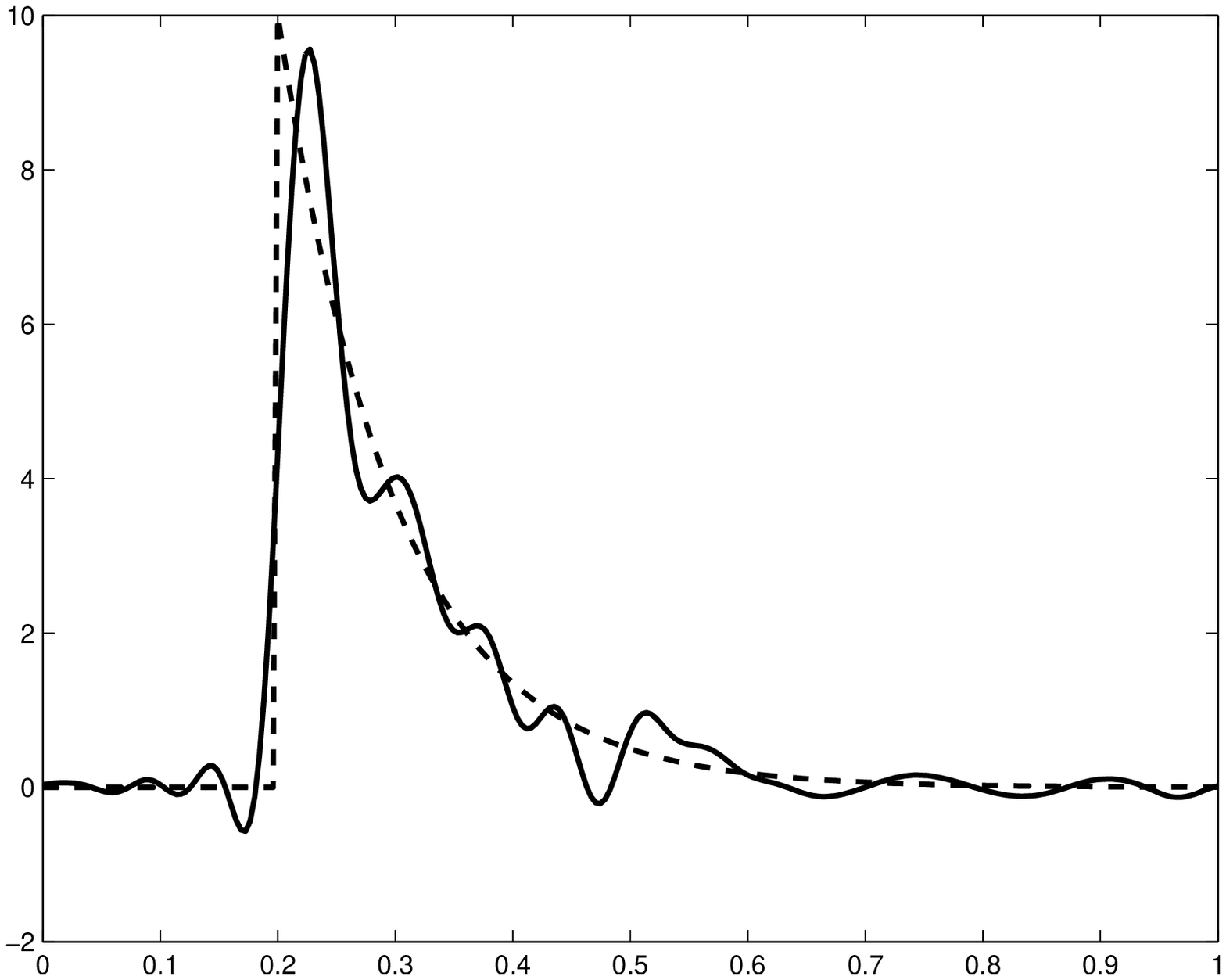} }
\subfigure[]
{ \includegraphics[width= 3.5cm]{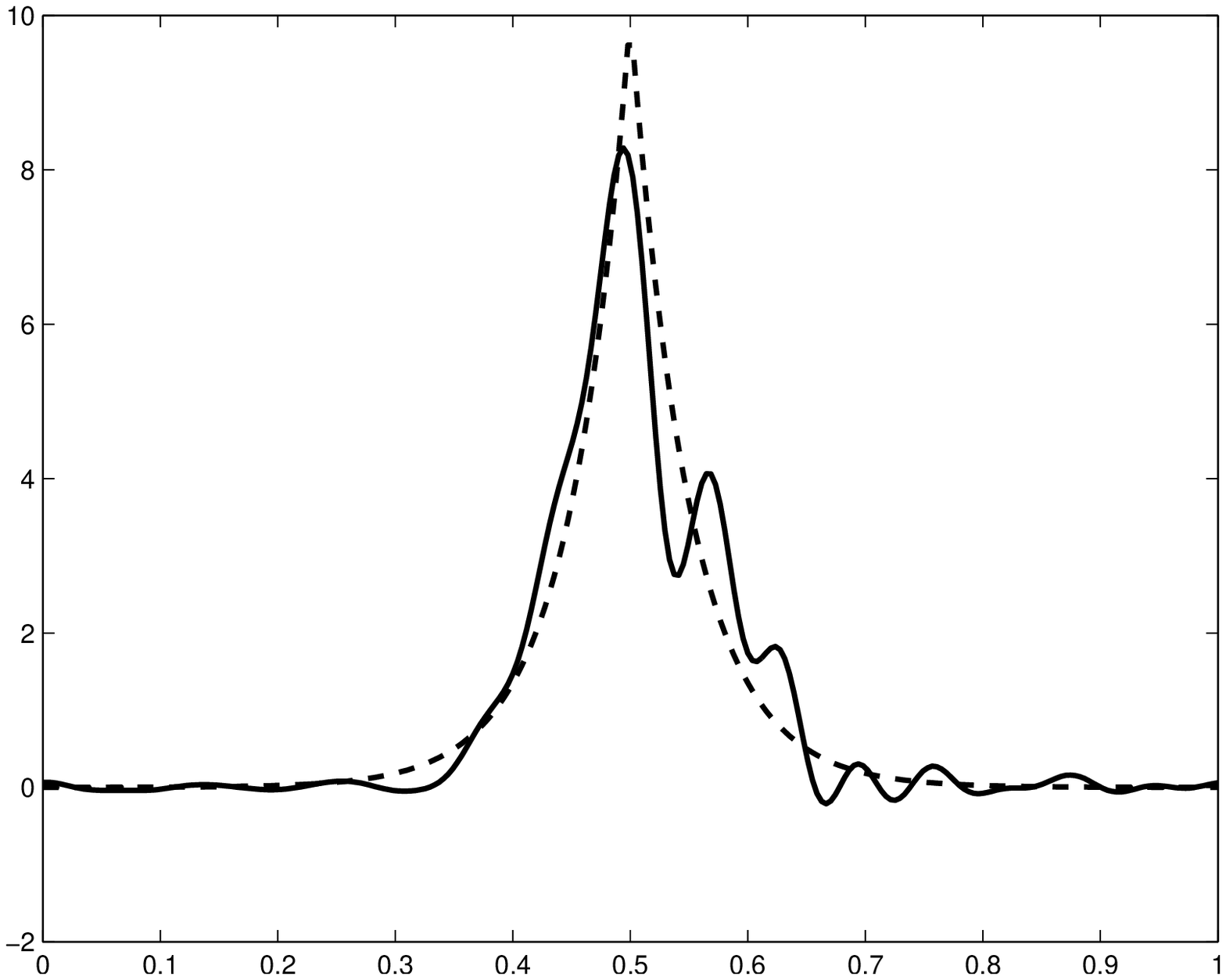} }
\subfigure[]
{ \includegraphics[width= 3.5cm]{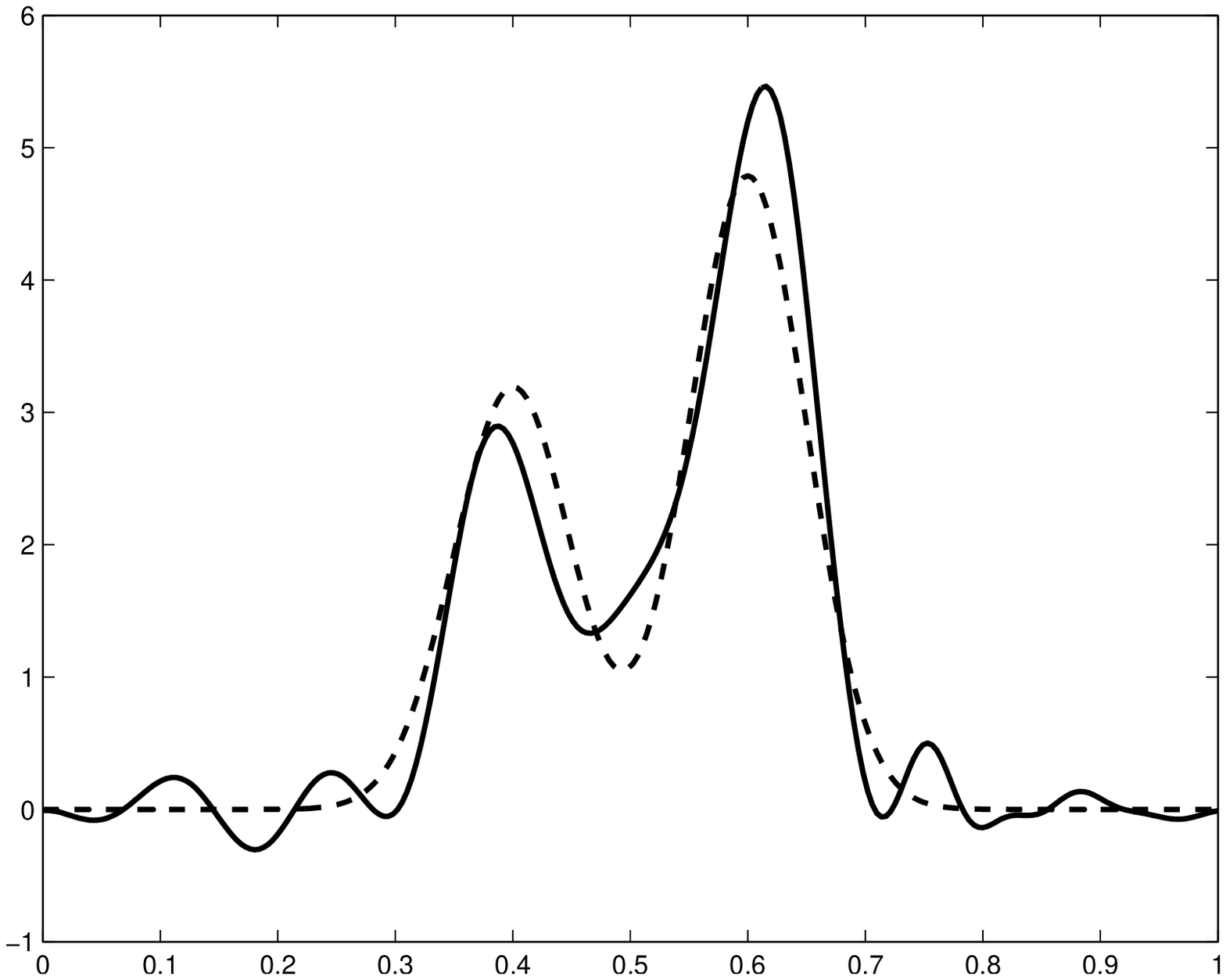} }
\caption{Typical reconstructions from a single simulation with  $n=200$, $j_{1} = j_{1}^{\ast}  = 4$ and $\delta =  (j_{1}-1)/\log_{2}(n) \approx 0.3925$: (a) Uniform, (b) Exponential, (c) Laplace, (d) MixtGauss. } \label{fig:example}
\end{figure}

To give an idea of the quality of $\hat{f}_{n}$, a typical example of estimation with $n=200$, $j_{1}^{\ast}  = 4$ and $\delta^{\ast} \approx 0.3925$  is given in  Figure \ref{fig:example}. Another possibility for choosing a smaller value for $\delta$ is to take $\alpha \neq 0$, and then to choose $\delta^{\ast} = \eta^{\ast} = (j_{1}-1-\alpha\log_{2}(\log(n)) )/\log_{2}(n)$ since $ \eta^{\ast} $  is the smallest constant $\eta$ that satisfies $2^{j_{1}} > n^{\eta} \log(n)^{\alpha} \geq 2^{j_{1}-1}$. For $n = 200$, $\alpha =0.5$ and $j_{1} = 4$ this yields to the choice $\delta \approx 0.2351$. However, as already remarked, the oracle inequality (\ref{eq:ineqoracle}) shows that taking $\alpha \neq 0$ may deteriorate the risk of $\hat{f}_{n}$.

The goal of this numerical section is thus to validate the above choices (\ref{eq:choice}) for $j_{1}$ and $\delta$, by studying the risk of $\hat{f}_{n}$ (compared to the risk of the oracle $\tilde{f}_{n}$) as a function of these two hyperparameters. More precisely, given $j_{1}$ and $\delta$, we define
\begin{equation}
R_{n}(j_{1},\delta) = \frac{\sum_{k =0}^{2^{j_{0}}-1}  (\hat{c}_{j_{0},k} - c_{j_{0},k})^{2}   + \sum_{j = j_{0}}^{j_{1}}  \sum_{k=0}^{2^{j}-1}   (\hat{\beta}_{j,k} \1_{\{|  \hat{\beta}_{j,k}  | \geqslant \tau_{j,k}\}} - \beta_{j,k})^{2}  + \sum_{j = j_{1} +1}^{J-1}   \sum_{k=0}^{2^{j}-1}\beta_{j,k}^{2}}{ \sum_{k =0}^{2^{j_{0}}-1} \sigma^{2}_{j_{0},k} + \sum_{j = j_{0}}^{j_{1}}  \sum_{k=0}^{2^{j}-1}   \min(\beta_{j,k}^{2} ,\sigma^{2}_{j,k})  + \sum_{j = j_{1} +1}^{J-1}   \sum_{k=0}^{2^{j}-1}\beta_{j,k}^{2}} \label{eq:riskj1}
\end{equation}
To illustrate the usefulness of taking random thresholds $\tau_{j,k}$ depending of the location $(j,k)$, we compare $R_{n}(j_{1},\delta)$ with the risk  $\tilde{R}_{n}(j_{1},\delta)$ of the wavelet estimator  obtained by taking a level-dependent threshold of the form $\delta  \sqrt{ j/n}$, as originally proposed  by \citeA{DJKP96aos}, where 
$$
\tilde{R}_{n}(j_{1},\delta) = \frac{\sum_{k =0}^{2^{j_{0}}-1}  (\hat{c}_{j_{0},k} - c_{j_{0},k})^{2}   + \sum_{j = j_{0}}^{j_{1}}  \sum_{k=0}^{2^{j}-1}   (\hat{\beta}_{j,k} \1_{\{|  \hat{\beta}_{j,k}  | \geqslant \delta \sqrt{j/n}\}} - \beta_{j,k})^{2}  + \sum_{j = j_{1} +1}^{J-1}   \sum_{k=0}^{2^{j}-1}\beta_{j,k}^{2}}{ \sum_{k =0}^{2^{j_{0}}-1} \sigma^{2}_{j_{0},k} + \sum_{j = j_{0}}^{j_{1}}  \sum_{k=0}^{2^{j}-1}   \min(\beta_{j,k}^{2} ,\sigma^{2}_{j,k})  + \sum_{j = j_{1} +1}^{J-1}   \sum_{k=0}^{2^{j}-1}\beta_{j,k}^{2}}
$$

\begin{figure}[h!]
\centering
\subfigure[]
{ \includegraphics[width=3.5cm]{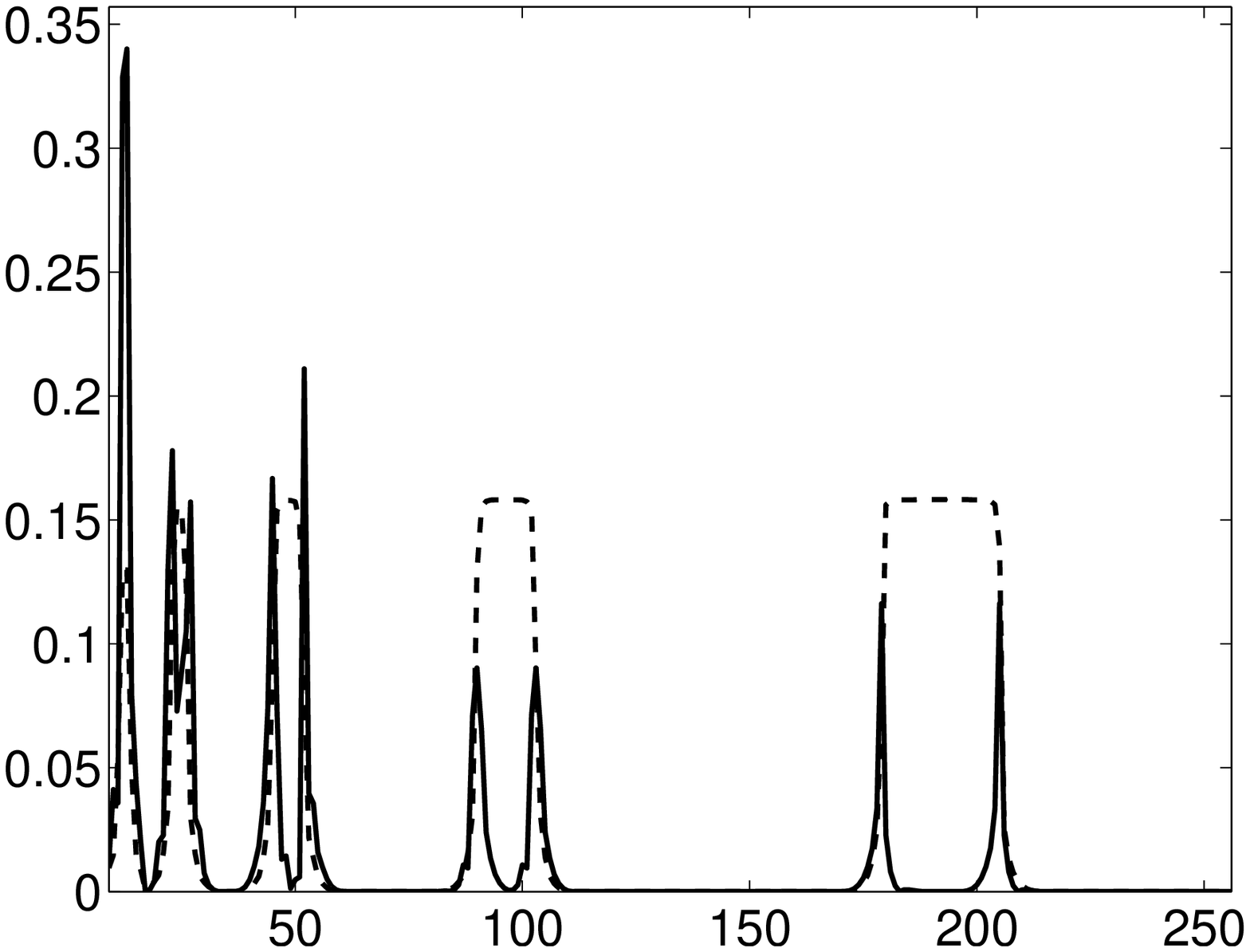} }
\subfigure[]
{ \includegraphics[width= 3.5cm]{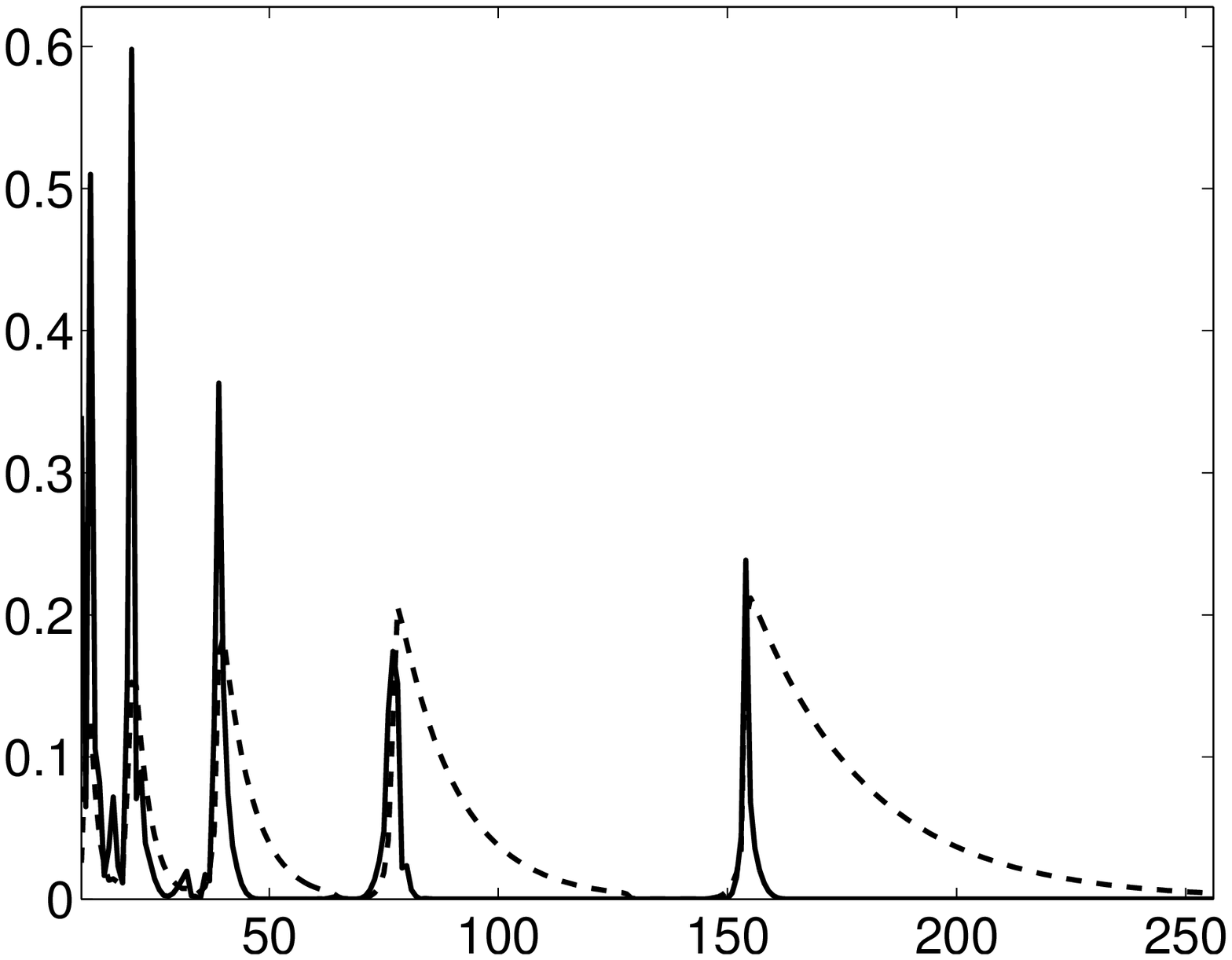} }
\subfigure[]
{ \includegraphics[width= 3.5cm]{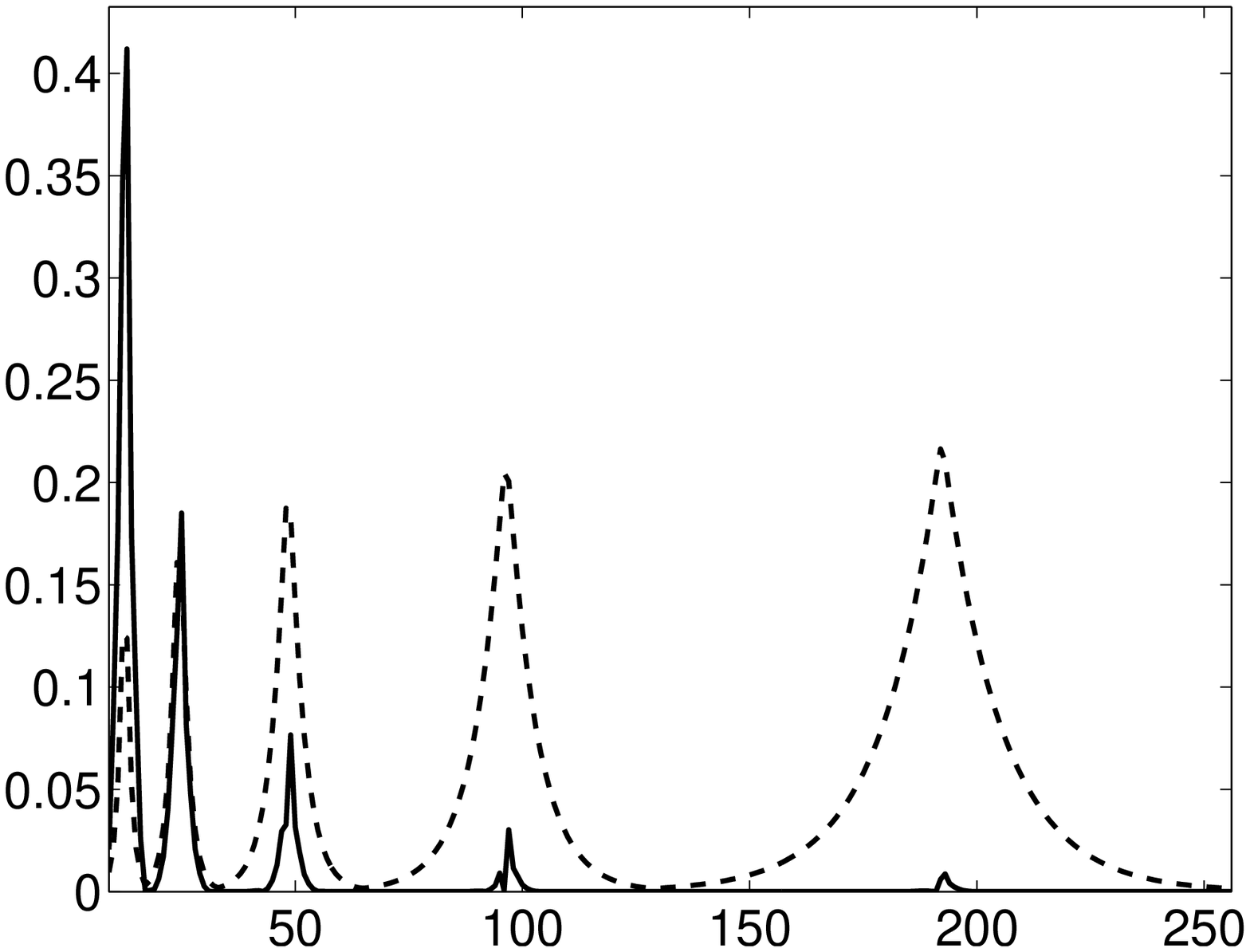} }
\subfigure[]
{ \includegraphics[width= 3.5cm]{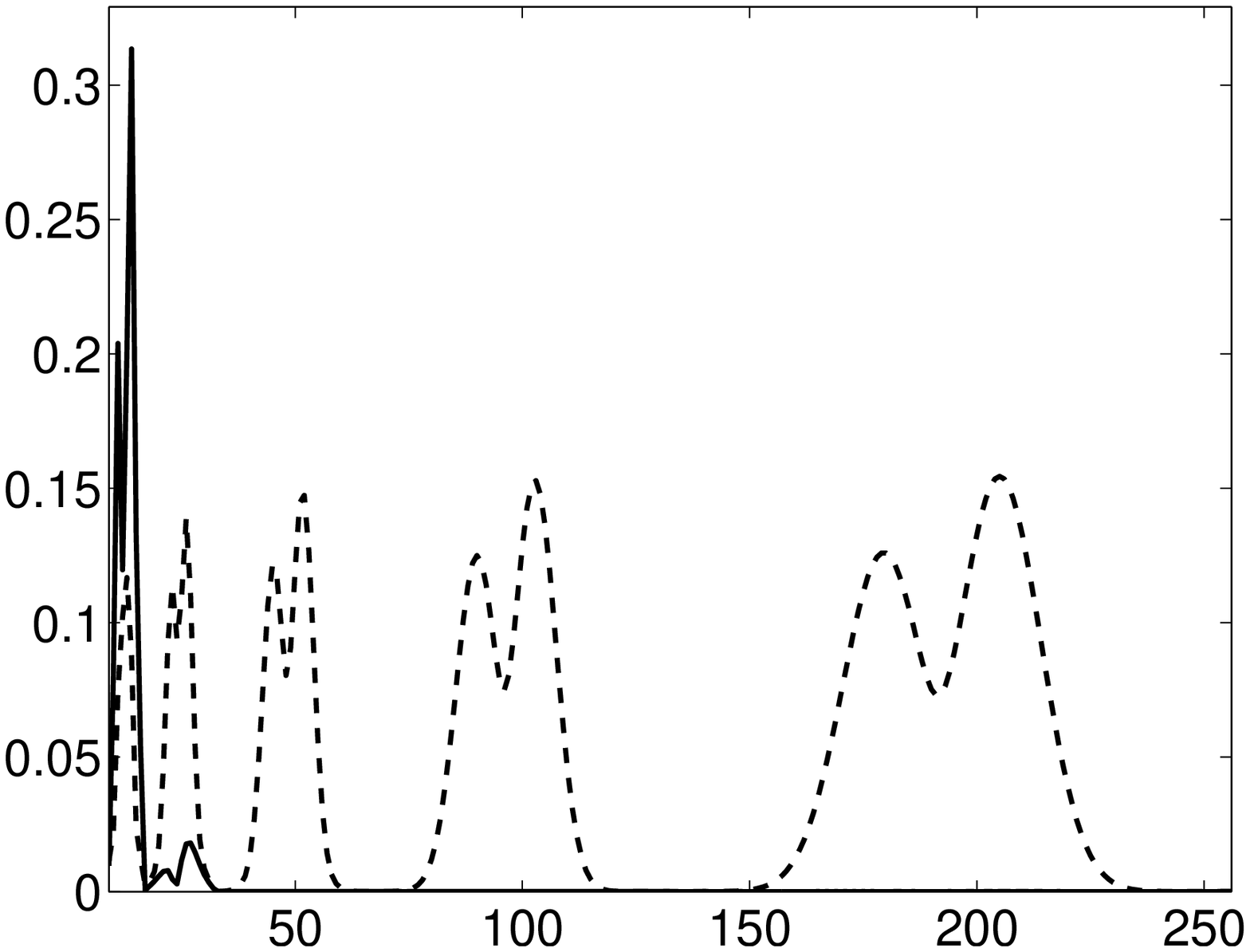} }
\caption{True wavelet coefficients $|\beta_{j,k}|$ (solid curve),  and standard deviation  $\sigma_{j,k}$ (dashed curves) as a function of $2^j+k$ for $j = j_{0},\ldots,J-1, k = 0,\ldots,2^j-1$ (with $j_{0} = 3$, $J=8$),  for each test density: (a) Uniform, (b) Exponential, (c) Laplace, (d) MixtGauss. } \label{fig:oracle}
\end{figure}

\begin{figure}[htbp]
\centering
\subfigure[Uniform - $j_{1} = 4$]
{ \includegraphics[width=3.5cm]{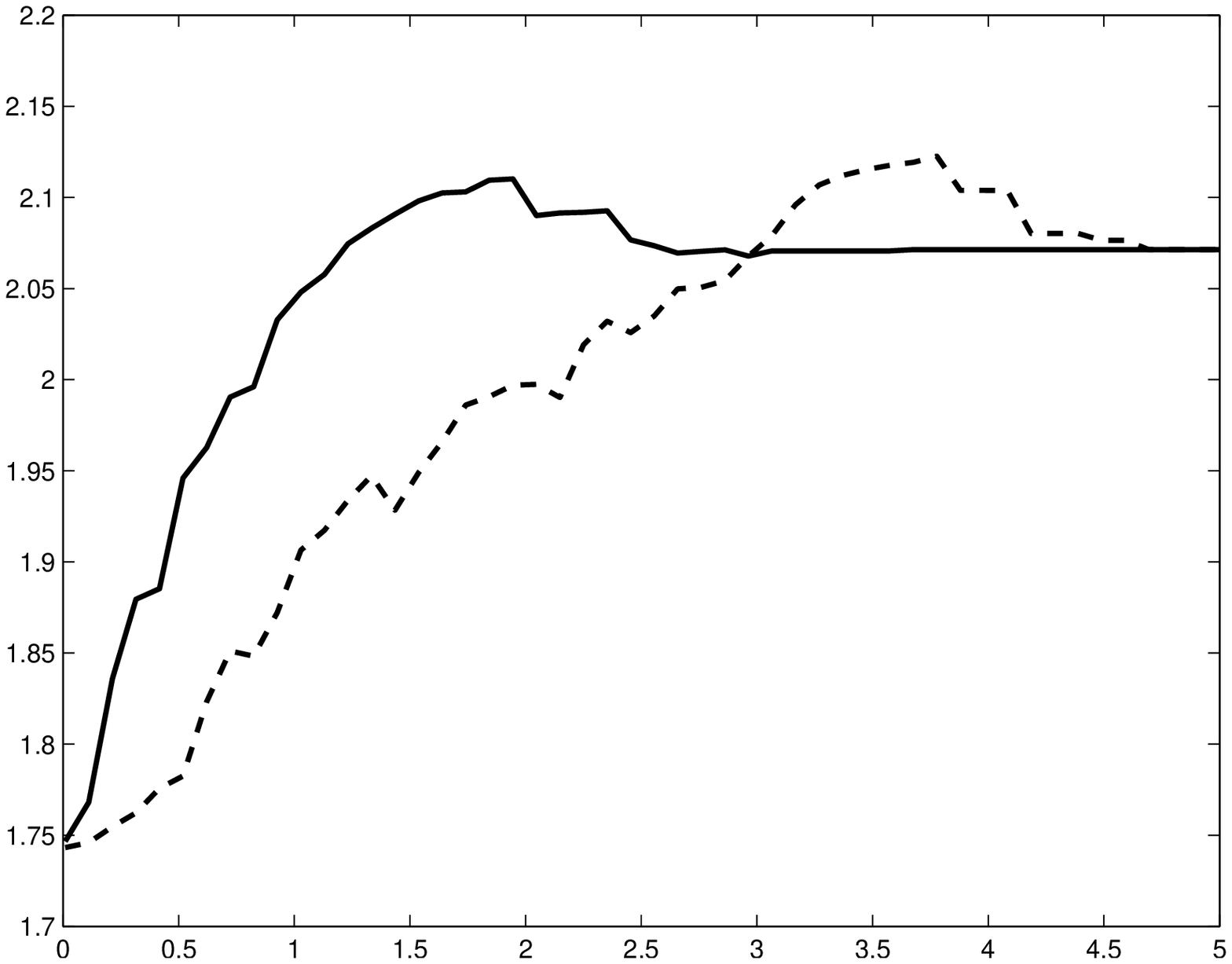} }
\subfigure[Uniform - $j_{1} = 5$]
{ \includegraphics[width=3.5cm]{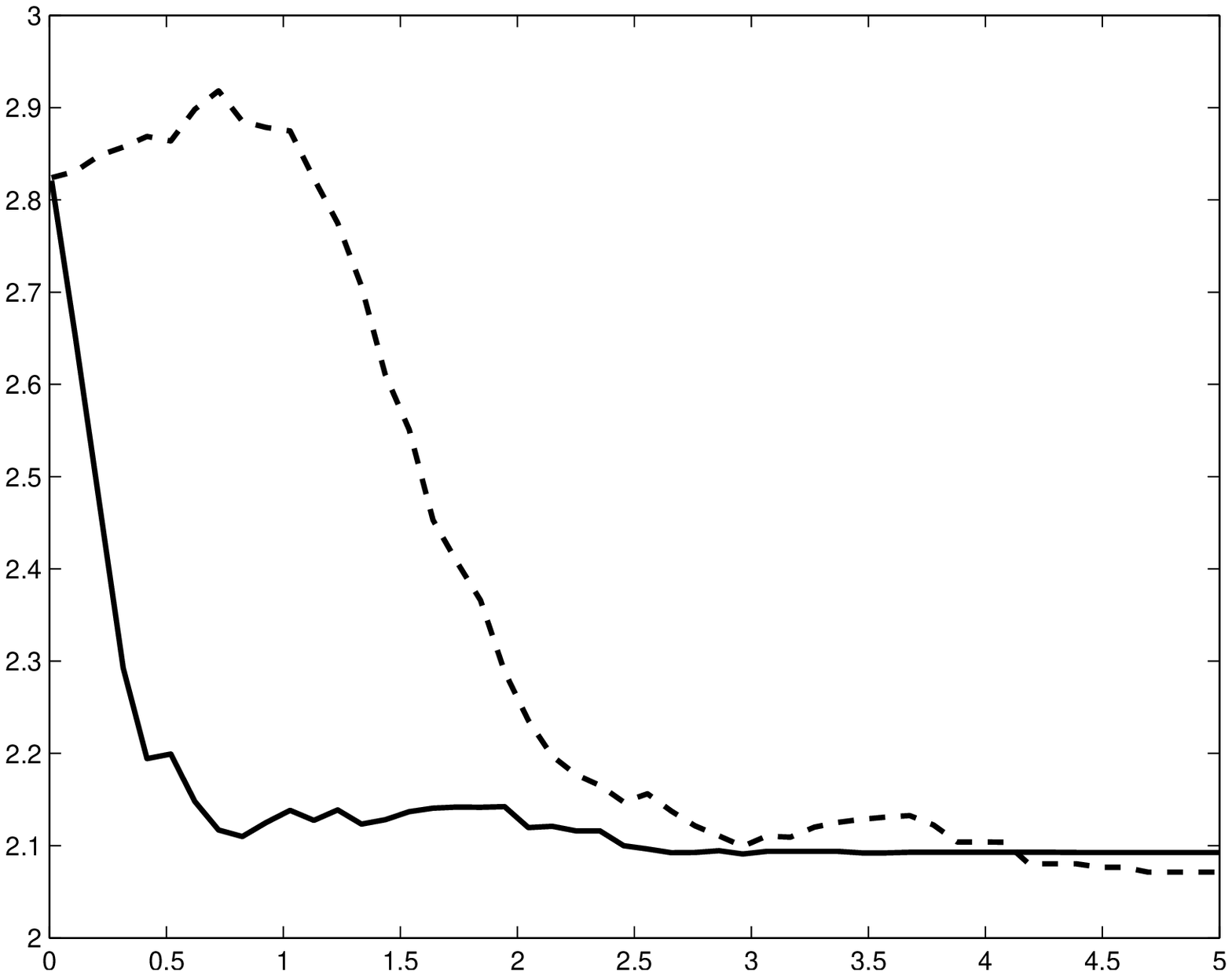} }
\subfigure[Uniform - $j_{1} = 6$]
{ \includegraphics[width=3.5cm]{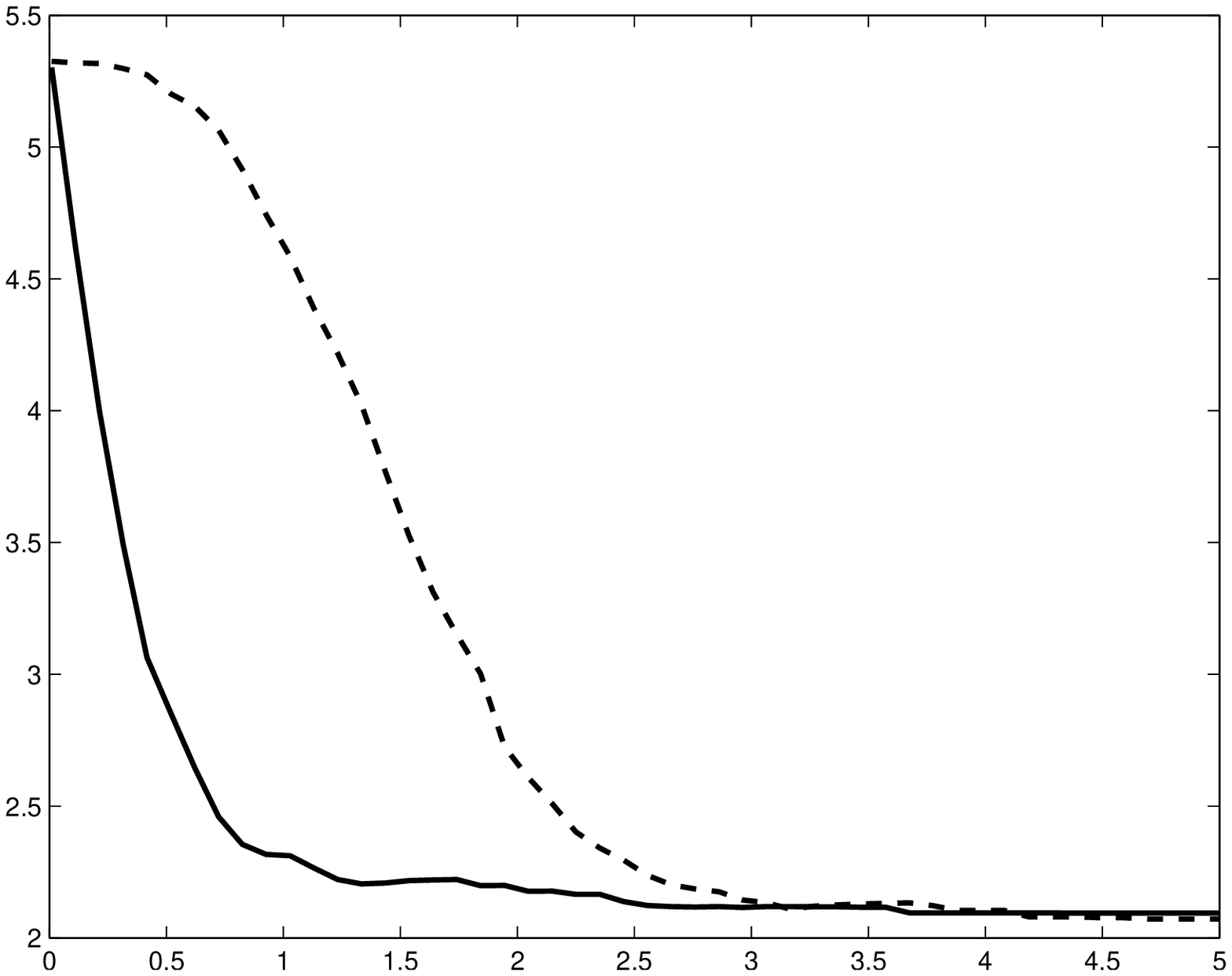} }

\subfigure[Exponential - $j_{1} = 4$]
{ \includegraphics[width=3.5cm]{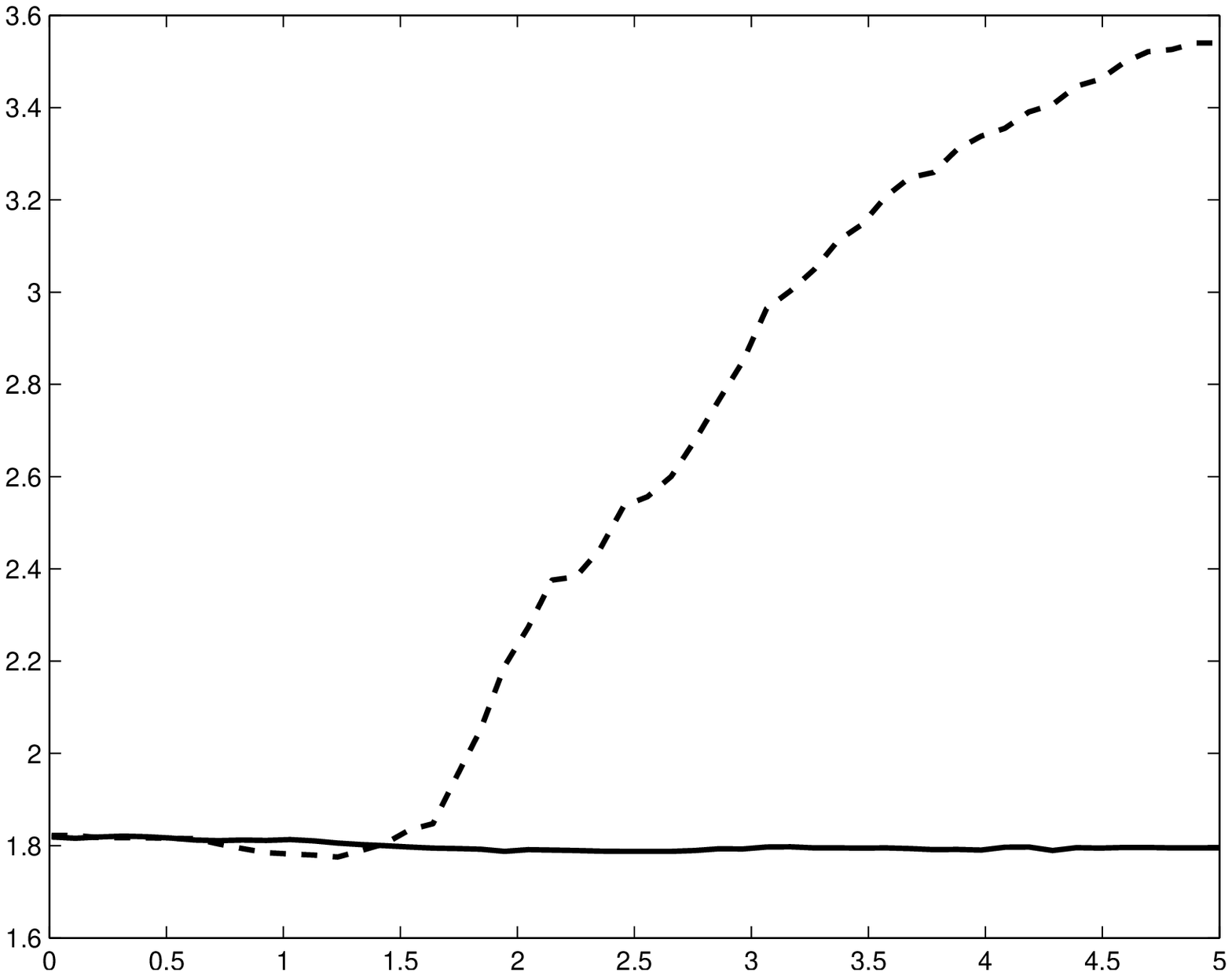} }
\subfigure[Exponential - $j_{1} = 5$]
{ \includegraphics[width=3.5cm]{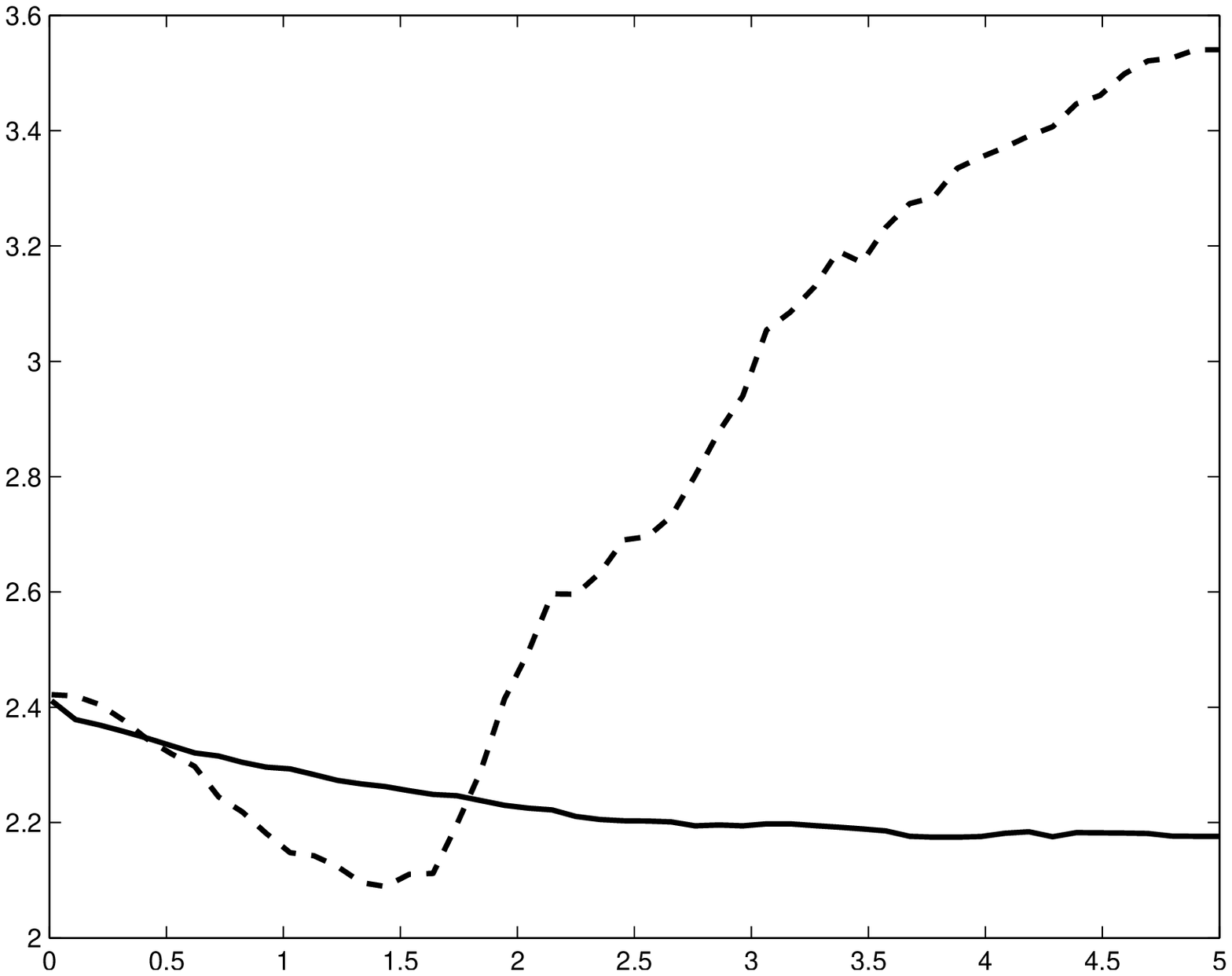} }
\subfigure[Exponential - $j_{1} = 6$]
{ \includegraphics[width=3.5cm]{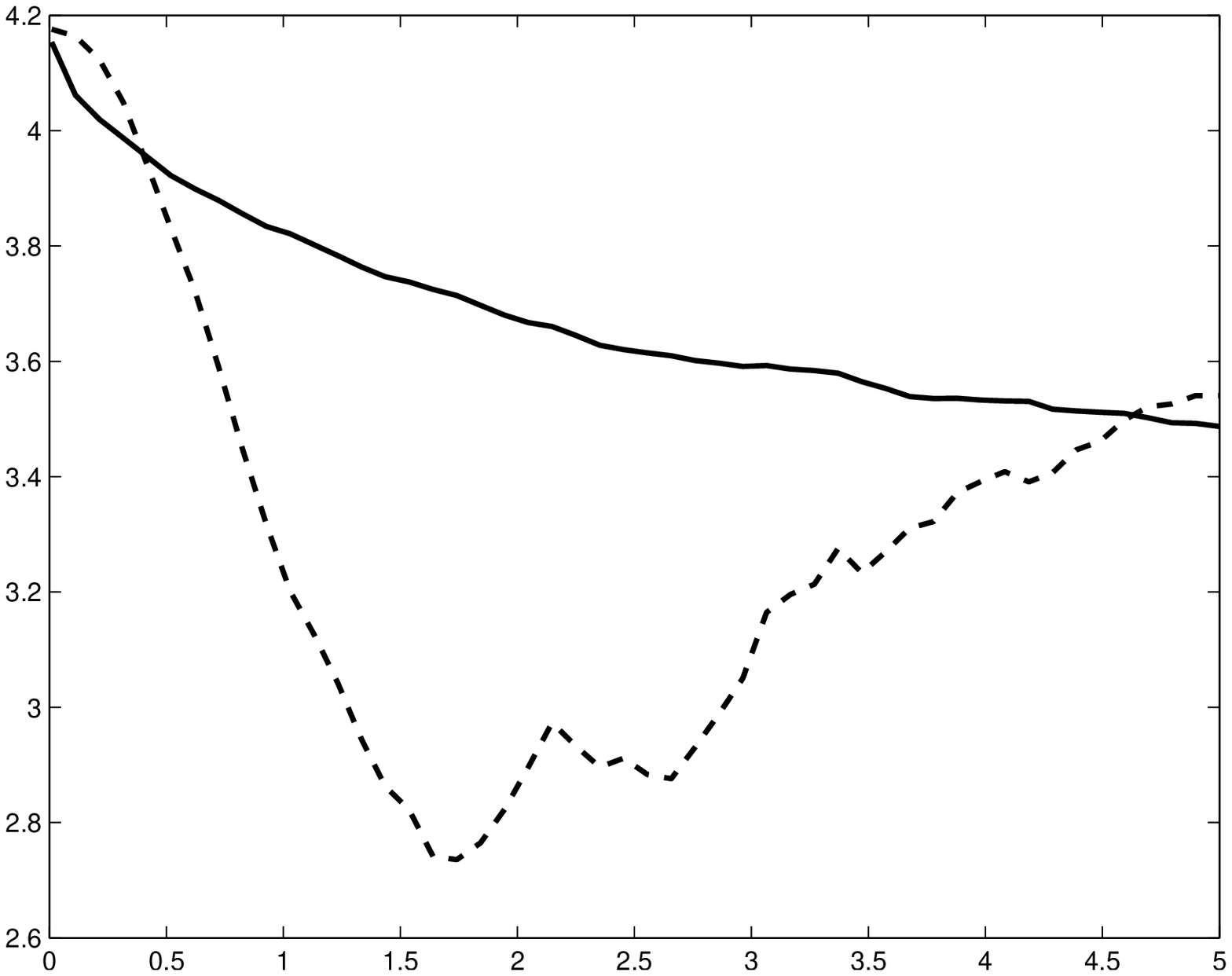} }

\subfigure[Laplace - $j_{1} = 4$]
{ \includegraphics[width=3.5cm]{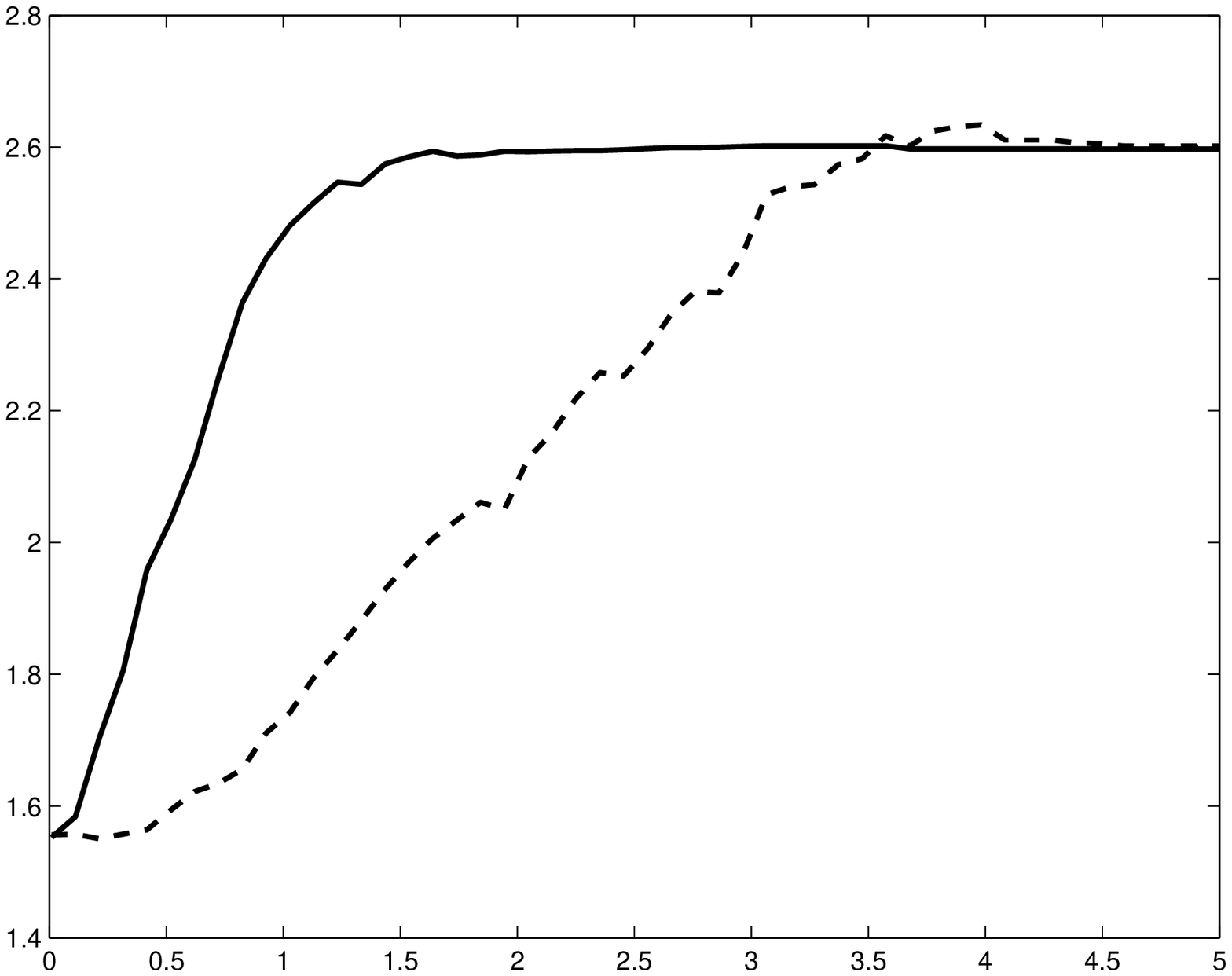} }
\subfigure[Laplace - $j_{1} = 5$]
{ \includegraphics[width=3.5cm]{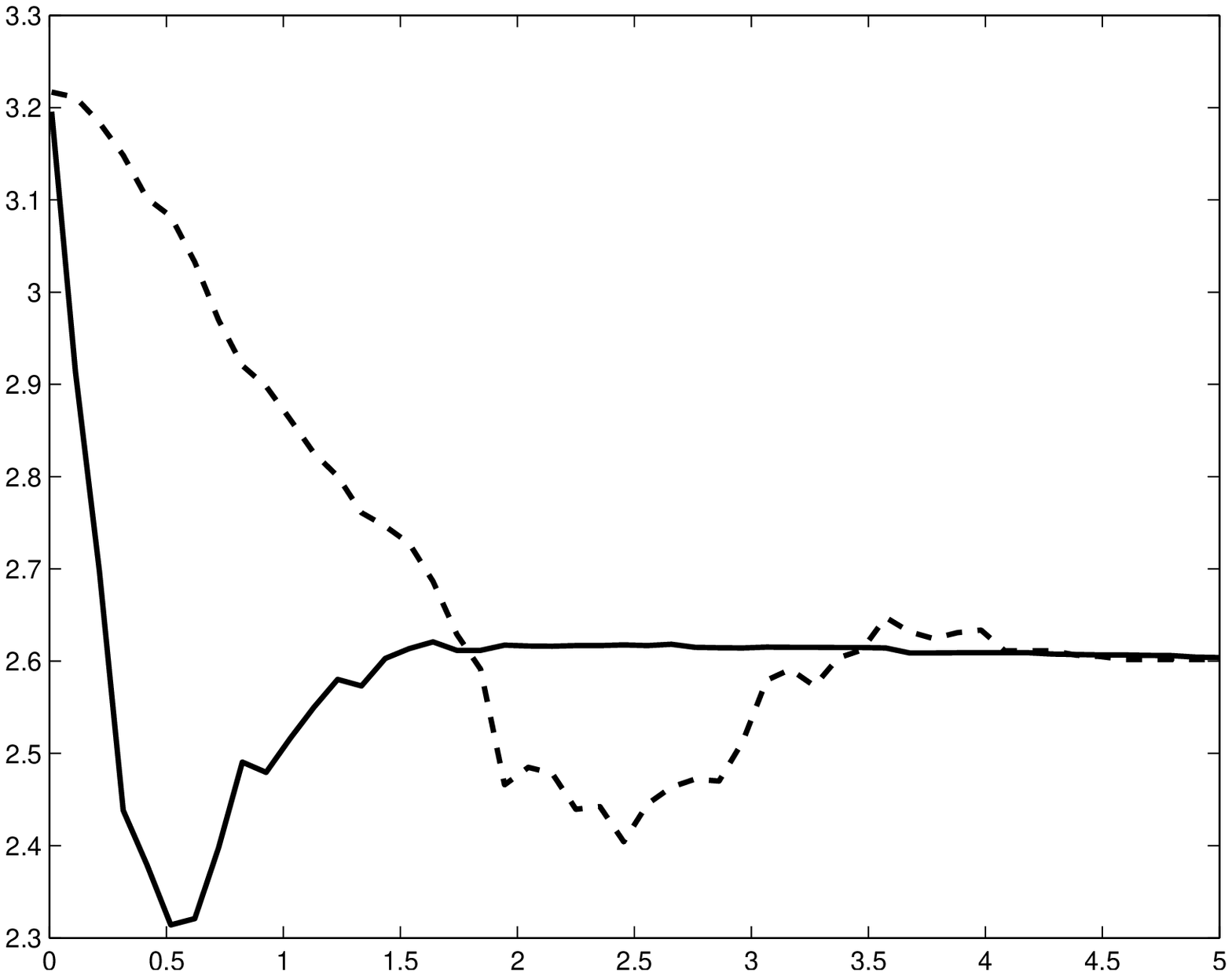} }
\subfigure[Laplace - $j_{1} = 6$]
{ \includegraphics[width=3.5cm]{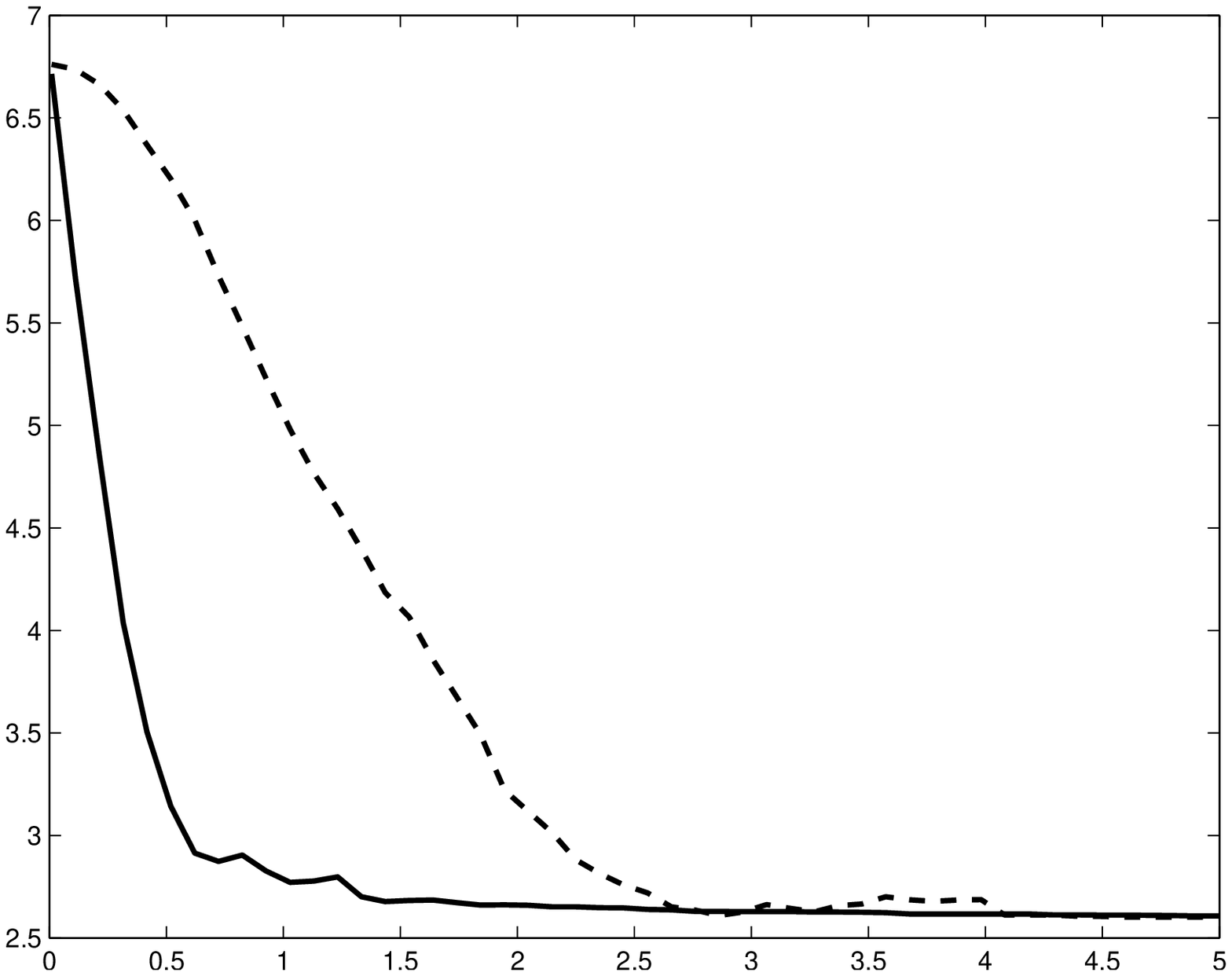} }

\subfigure[MixtGauss - $j_{1} = 4$]
{ \includegraphics[width=3.5cm]{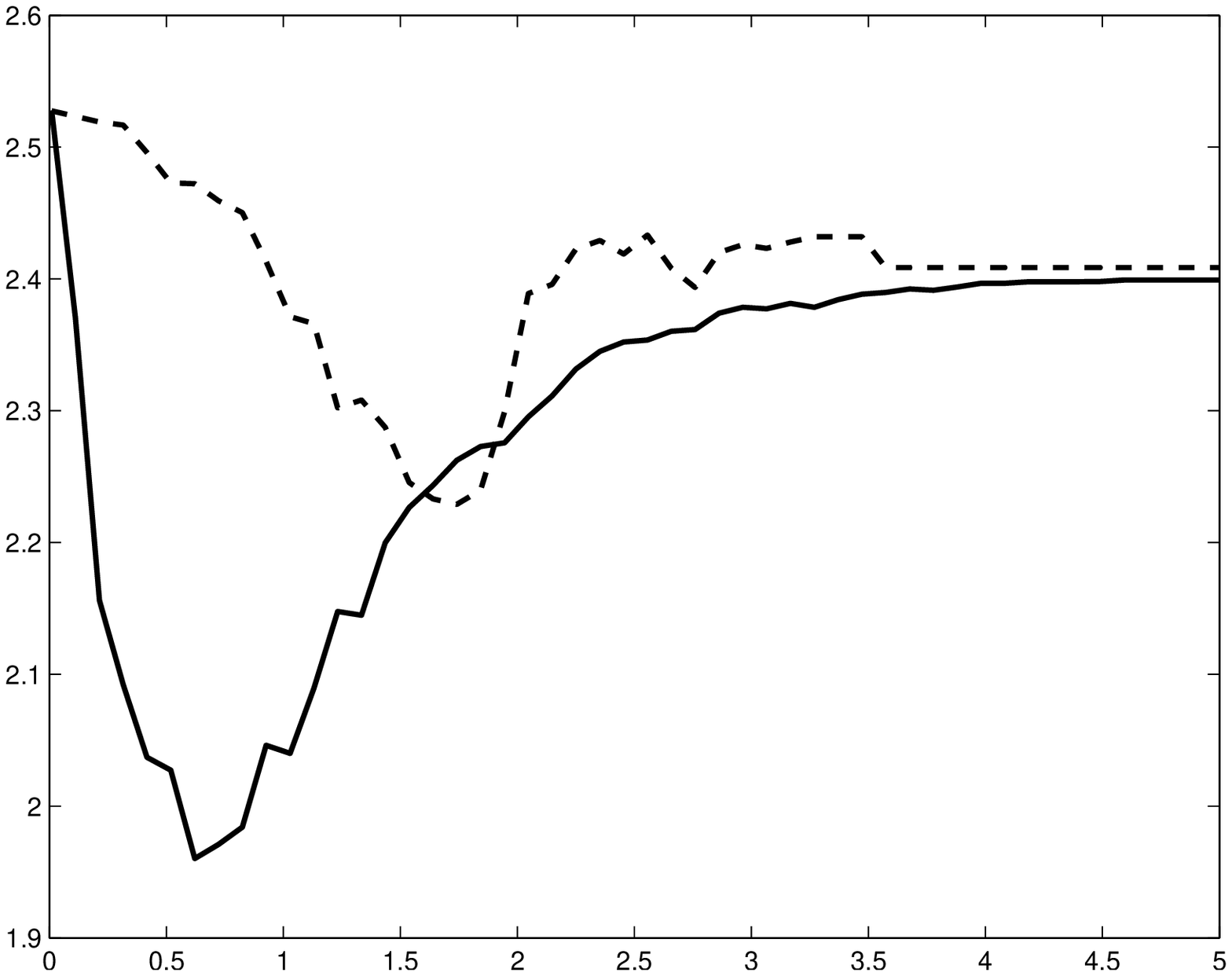} }
\subfigure[MixtGauss - $j_{1} = 5$]
{ \includegraphics[width=3.5cm]{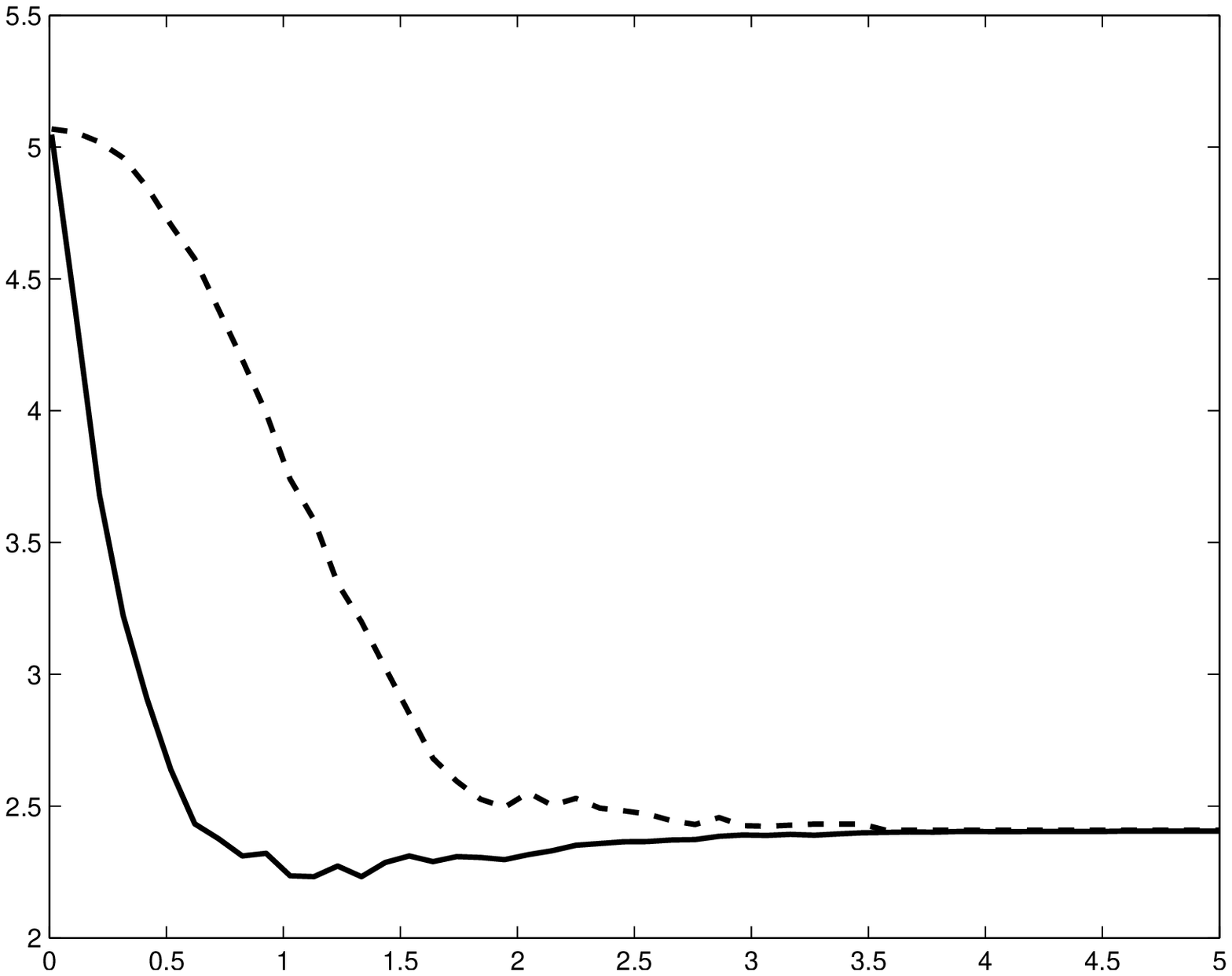} }
\subfigure[MixtGauss - $j_{1} = 6$]
{ \includegraphics[width=3.5cm]{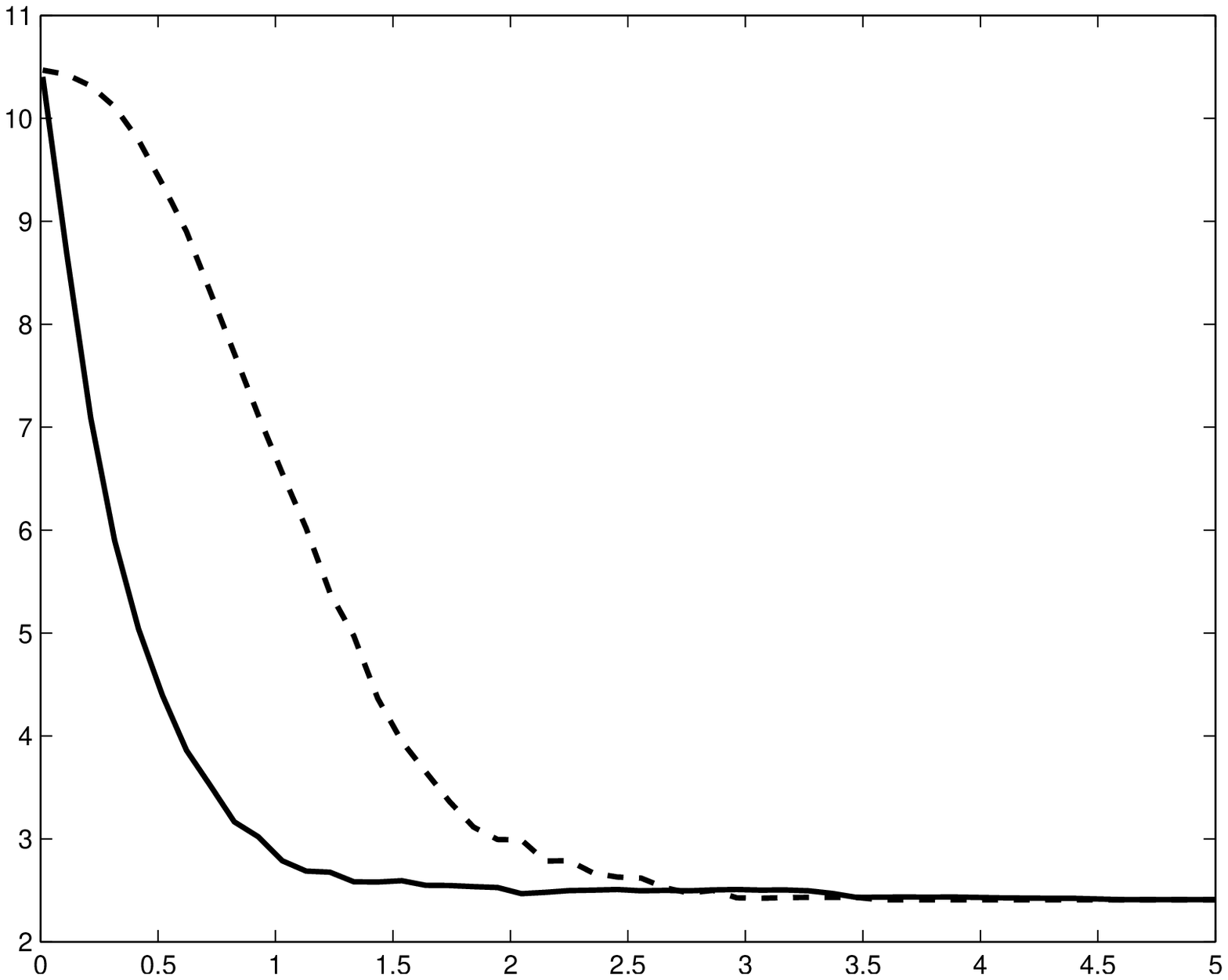} }

\caption{Direct density estimation with $n=100$. Evolution of  $R_{n}(j_{1},\delta)$ (solid curves) and $\tilde{R}_{n}(j_{1},\delta)$ (dashed curves) as a function of $\delta \in [0,5]$  for different values of $j_{1} \geq  j_{1}^{\ast} = \lfloor \frac{1}{2} \log_{2}(n) \rfloor + 1 = 4$} \label{fig:riskn100}
\end{figure}

\begin{figure}[htbp]
\centering
\subfigure[Uniform - $j_{1} = 4$]
{ \includegraphics[width=3.5cm]{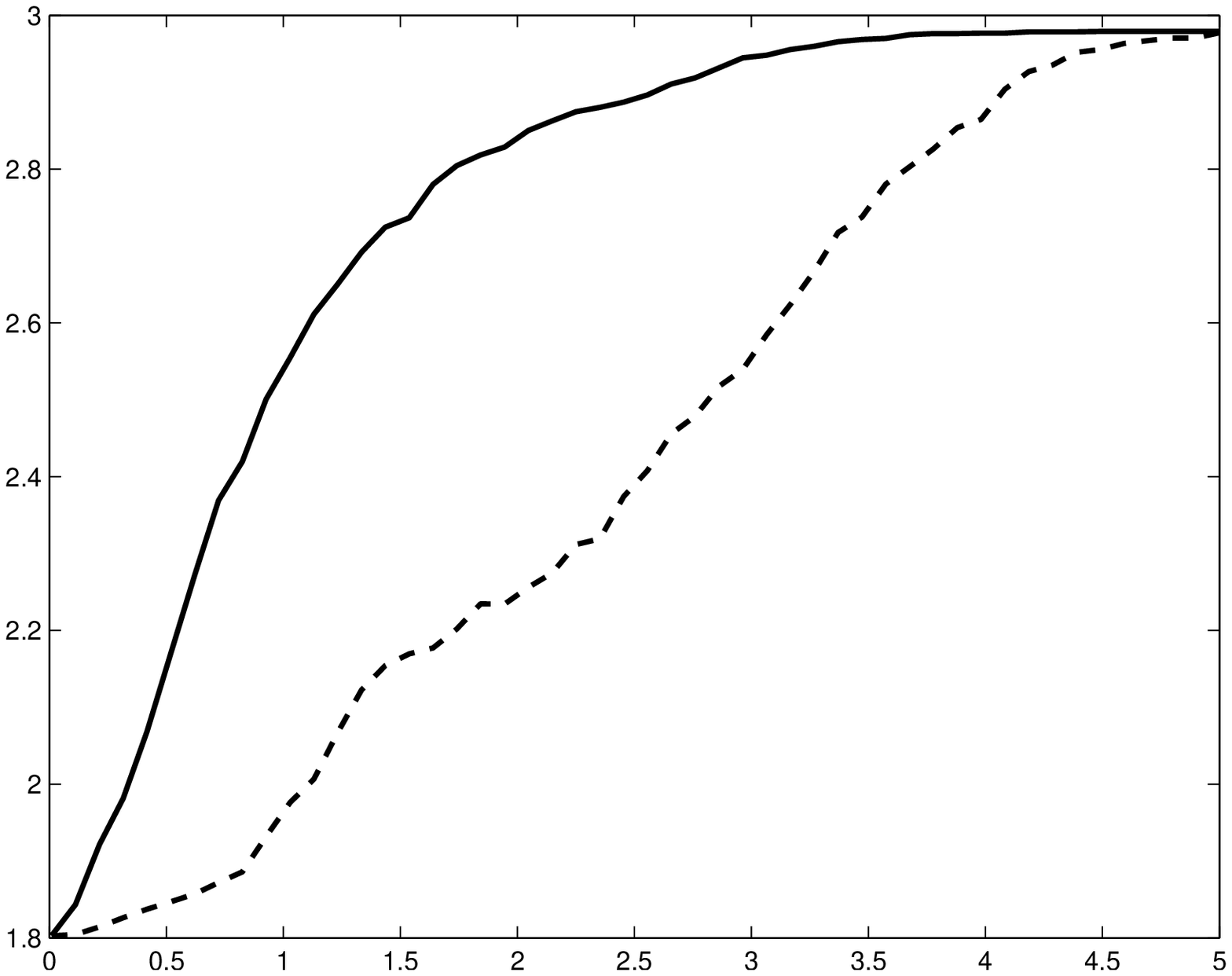} }
\subfigure[Uniform - $j_{1} = 5$]
{ \includegraphics[width=3.5cm]{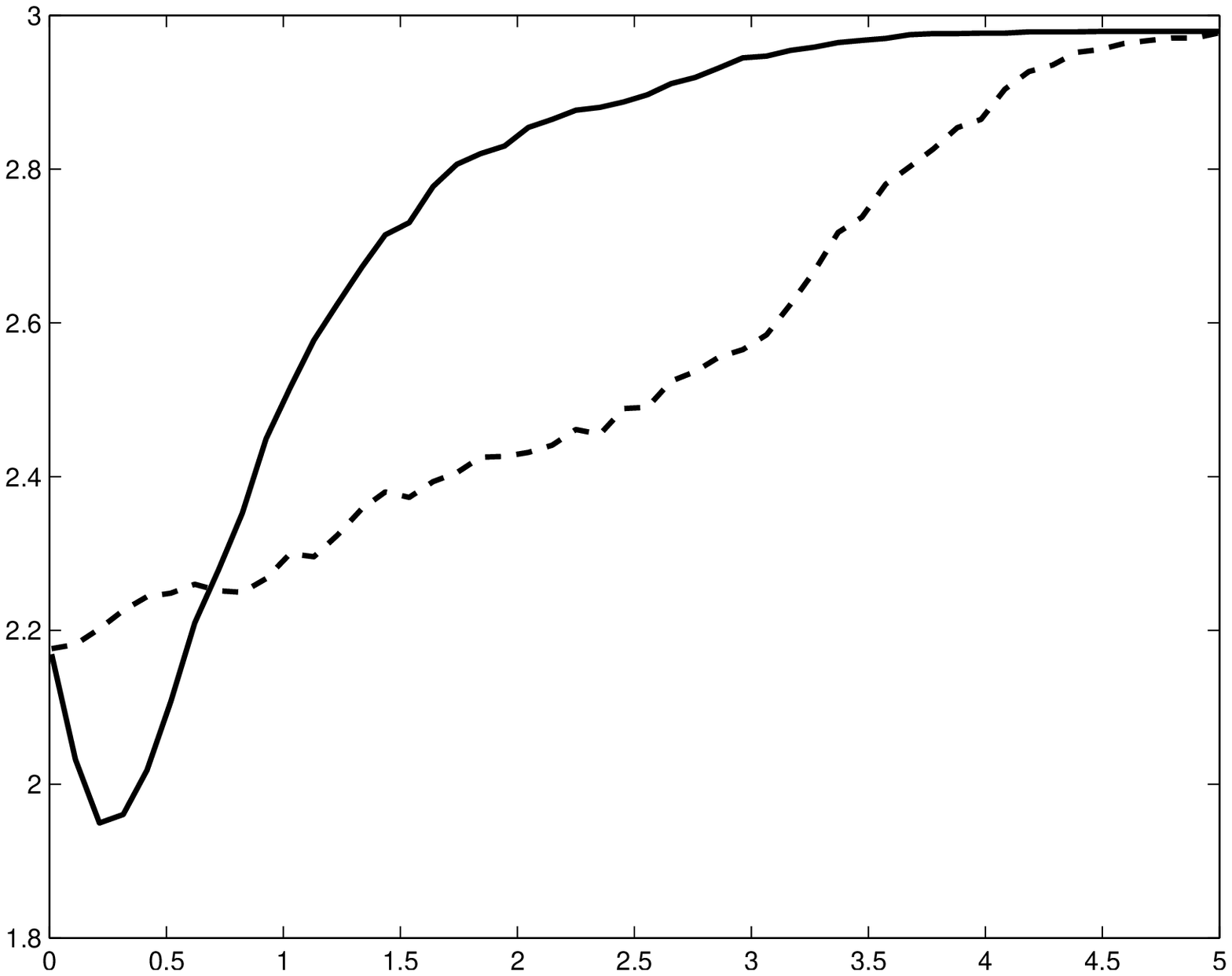} }
\subfigure[Uniform - $j_{1} = 6$]
{ \includegraphics[width=3.5cm]{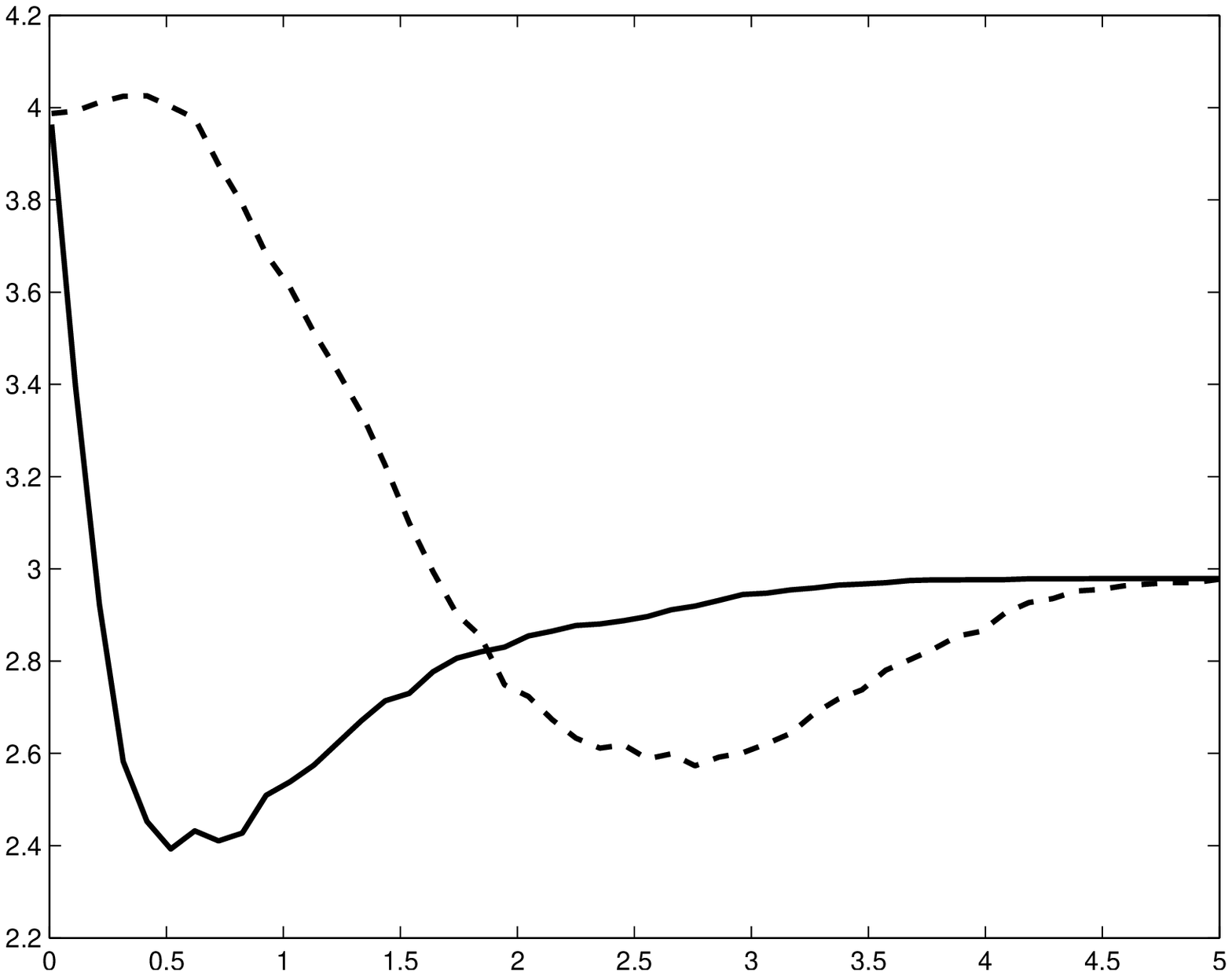} }
%\subfigure[Uniform - $j_{1} = 7$]
%{ \includegraphics[width=3.5cm]{Unif_n200_j1_7.eps} }

\subfigure[Exponential - $j_{1} = 4$]
{ \includegraphics[width=3.5cm]{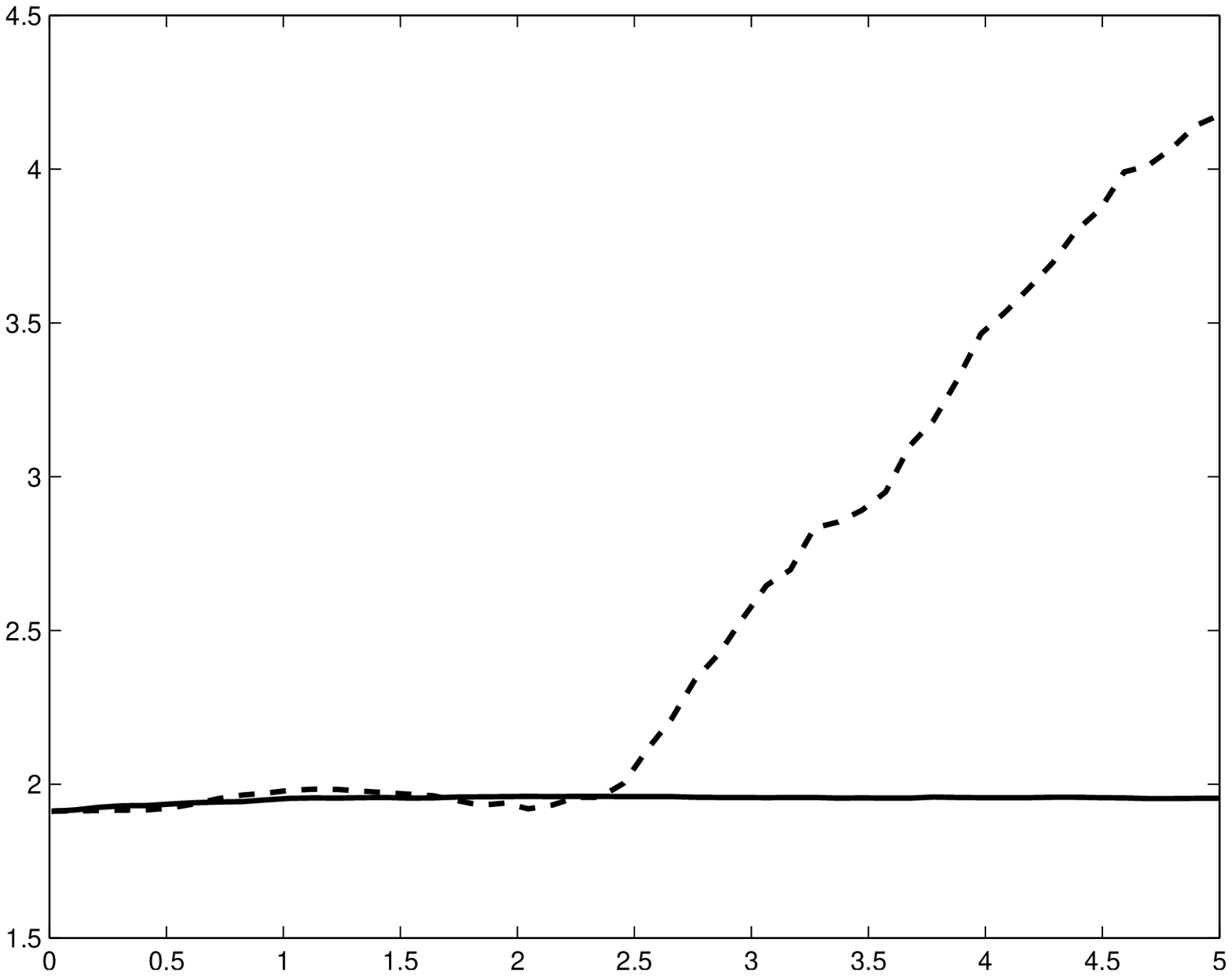} }
\subfigure[Exponential - $j_{1} = 5$]
{ \includegraphics[width=3.5cm]{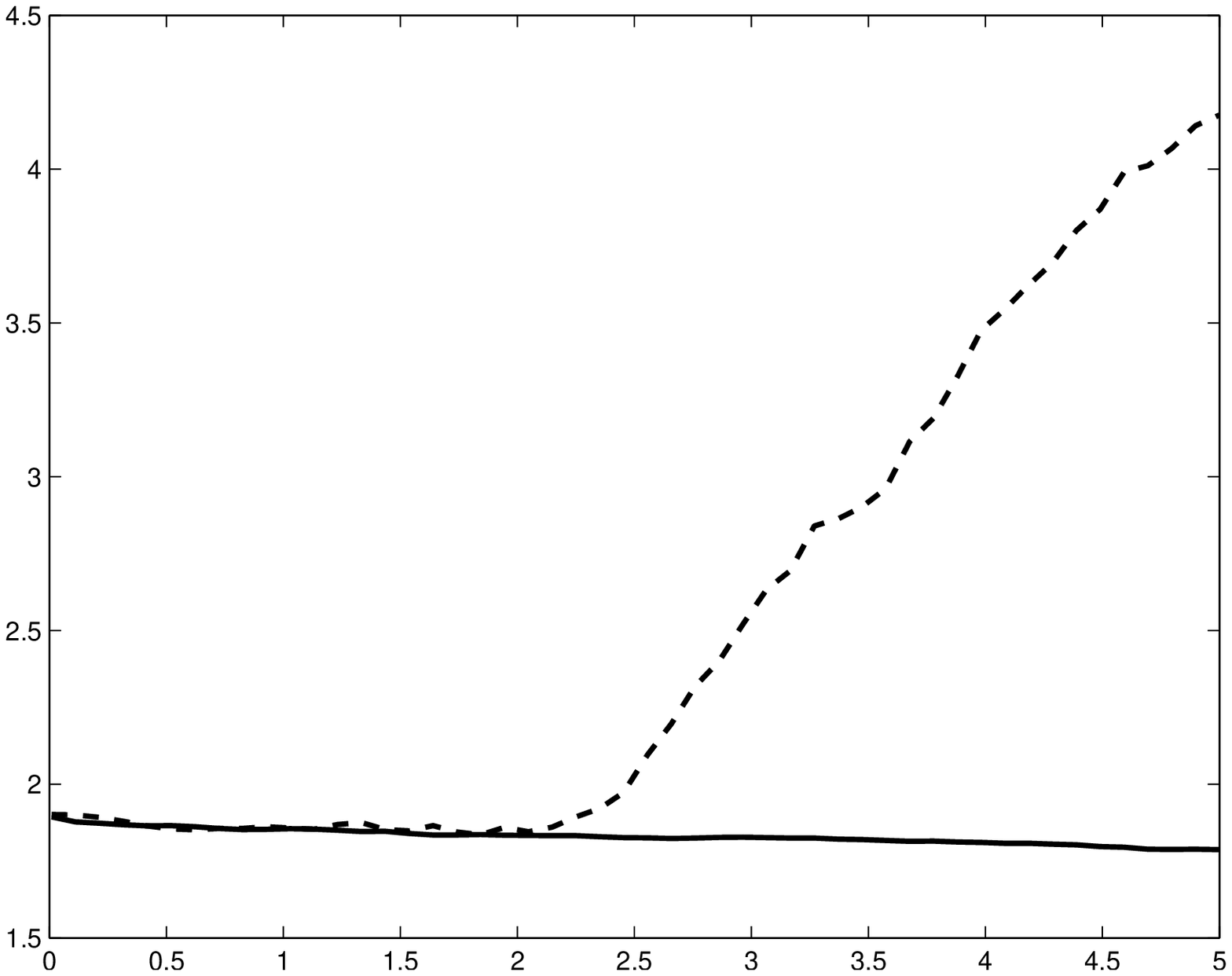} }
\subfigure[Exponential - $j_{1} = 6$]
{ \includegraphics[width=3.5cm]{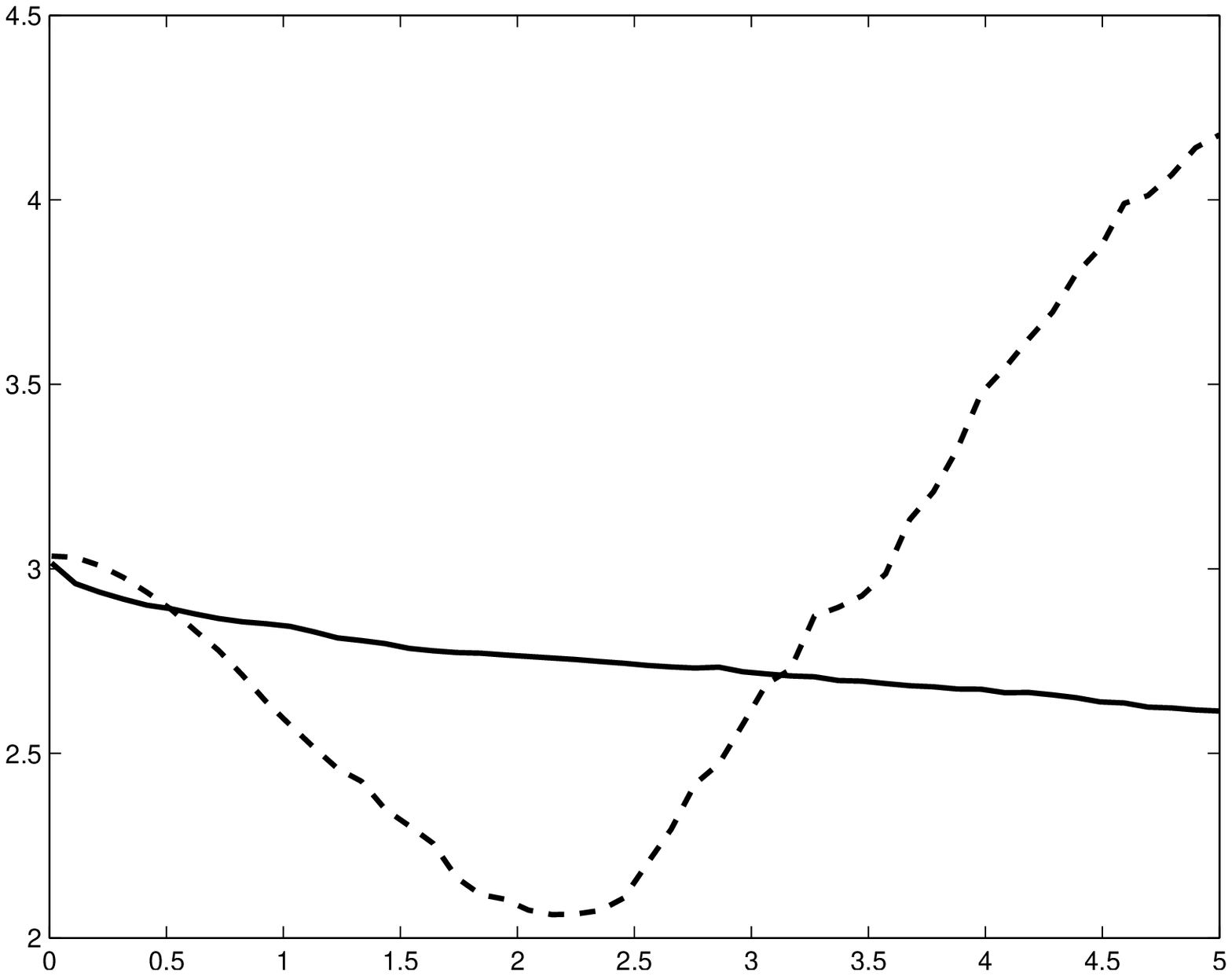} }
%\subfigure[Exponential - $j_{1} = 7$]
%{ \includegraphics[width=3.5cm]{Exp_n200_j1_7.eps} }

\subfigure[Laplace - $j_{1} = 4$]
{ \includegraphics[width=3.5cm]{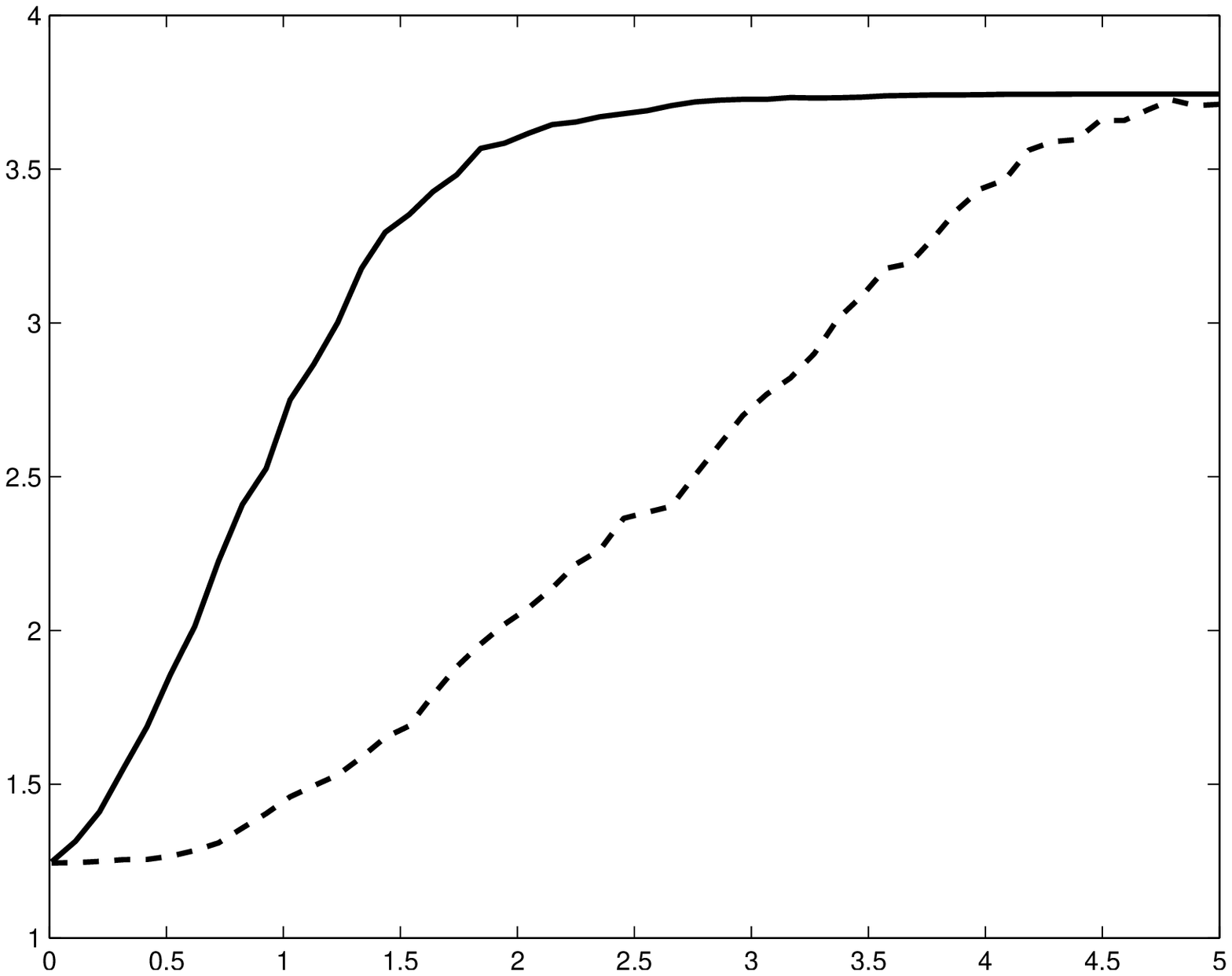} }
\subfigure[Laplace - $j_{1} = 5$]
{ \includegraphics[width=3.5cm]{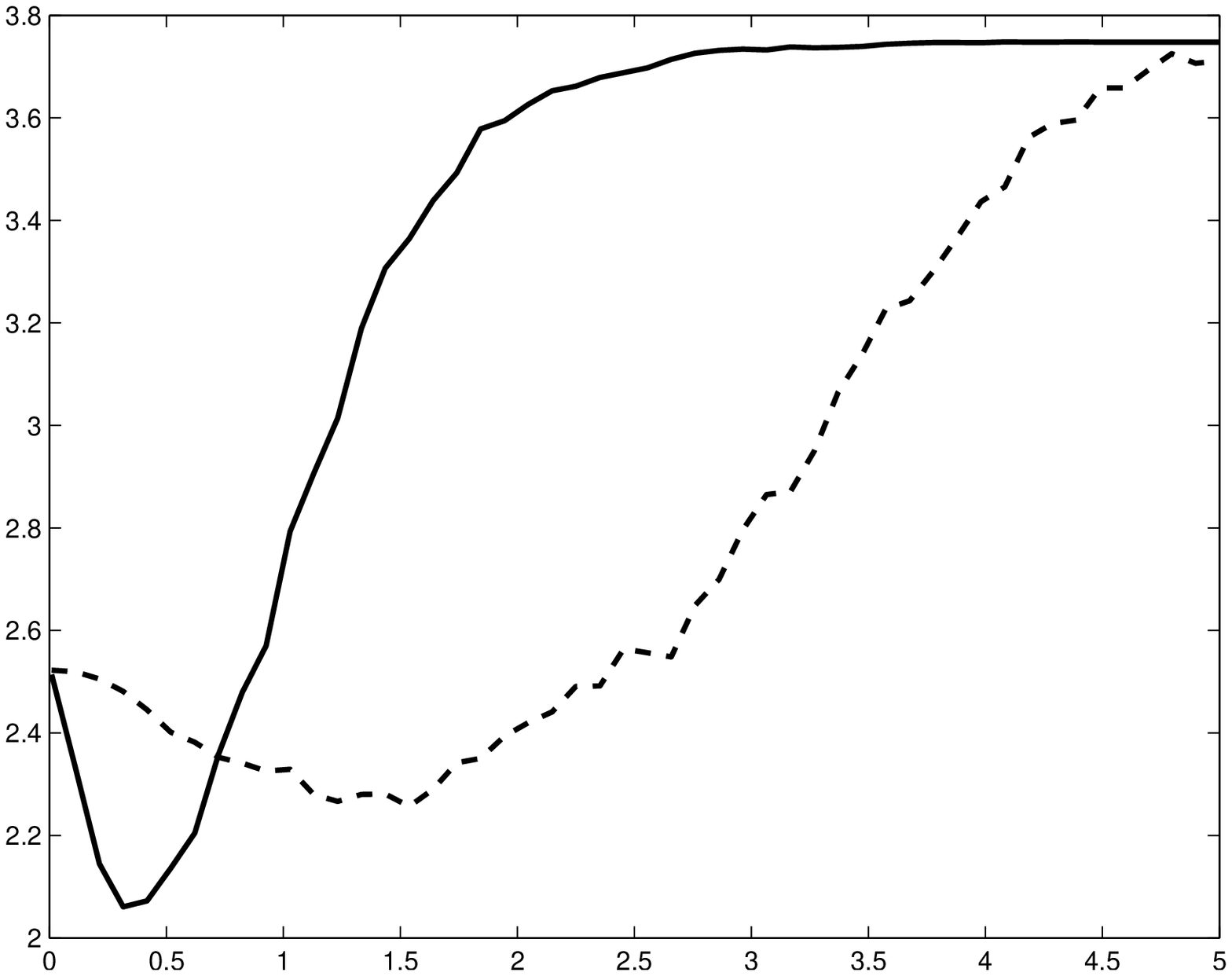} }
\subfigure[Laplace - $j_{1} = 6$]
{ \includegraphics[width=3.5cm]{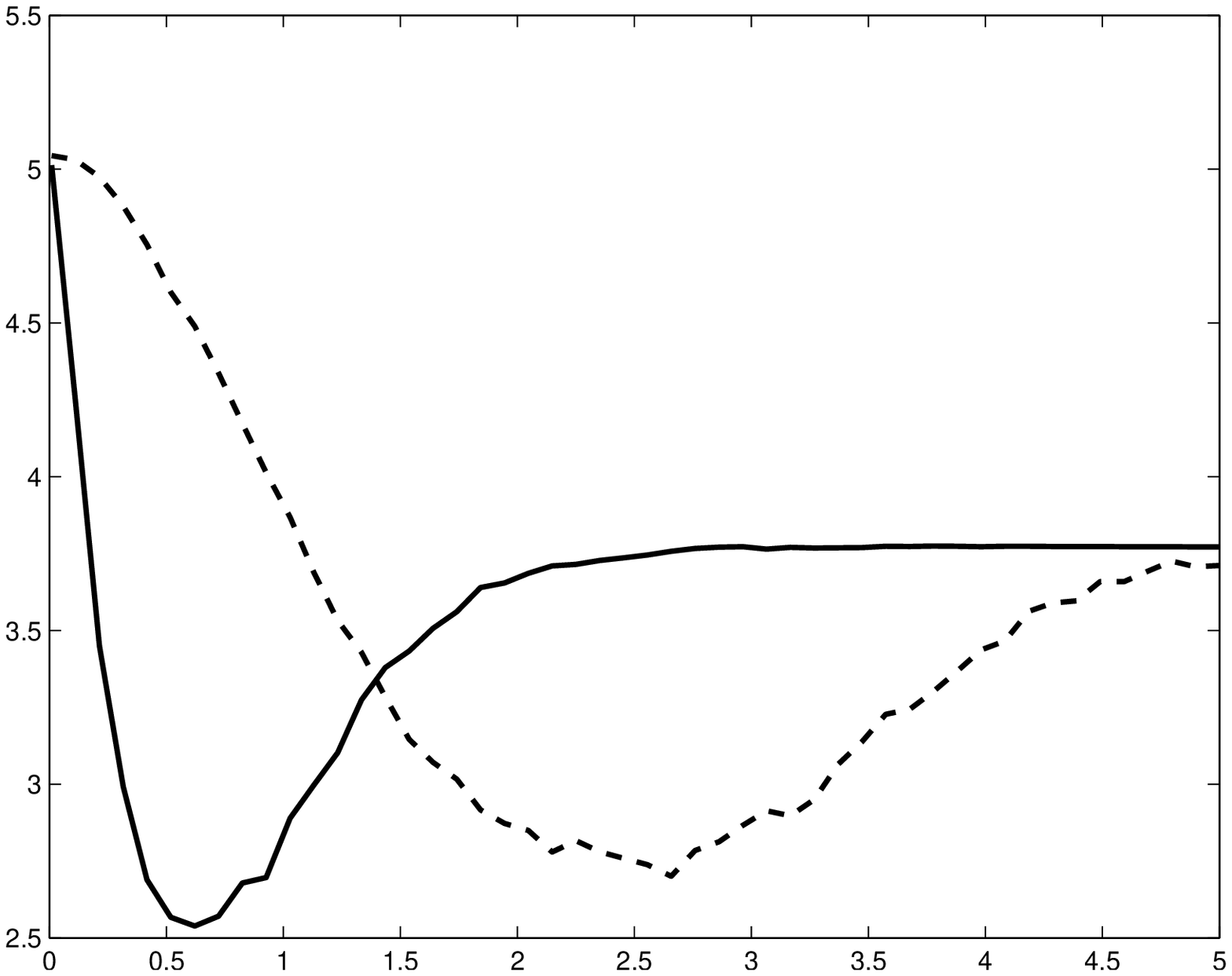} }
%\subfigure[Laplace - $j_{1} = 7$]
%{ \includegraphics[width=3.5cm]{Laplace_n200_j1_7.eps} }

\subfigure[MixtGauss - $j_{1} = 4$]
{ \includegraphics[width=3.5cm]{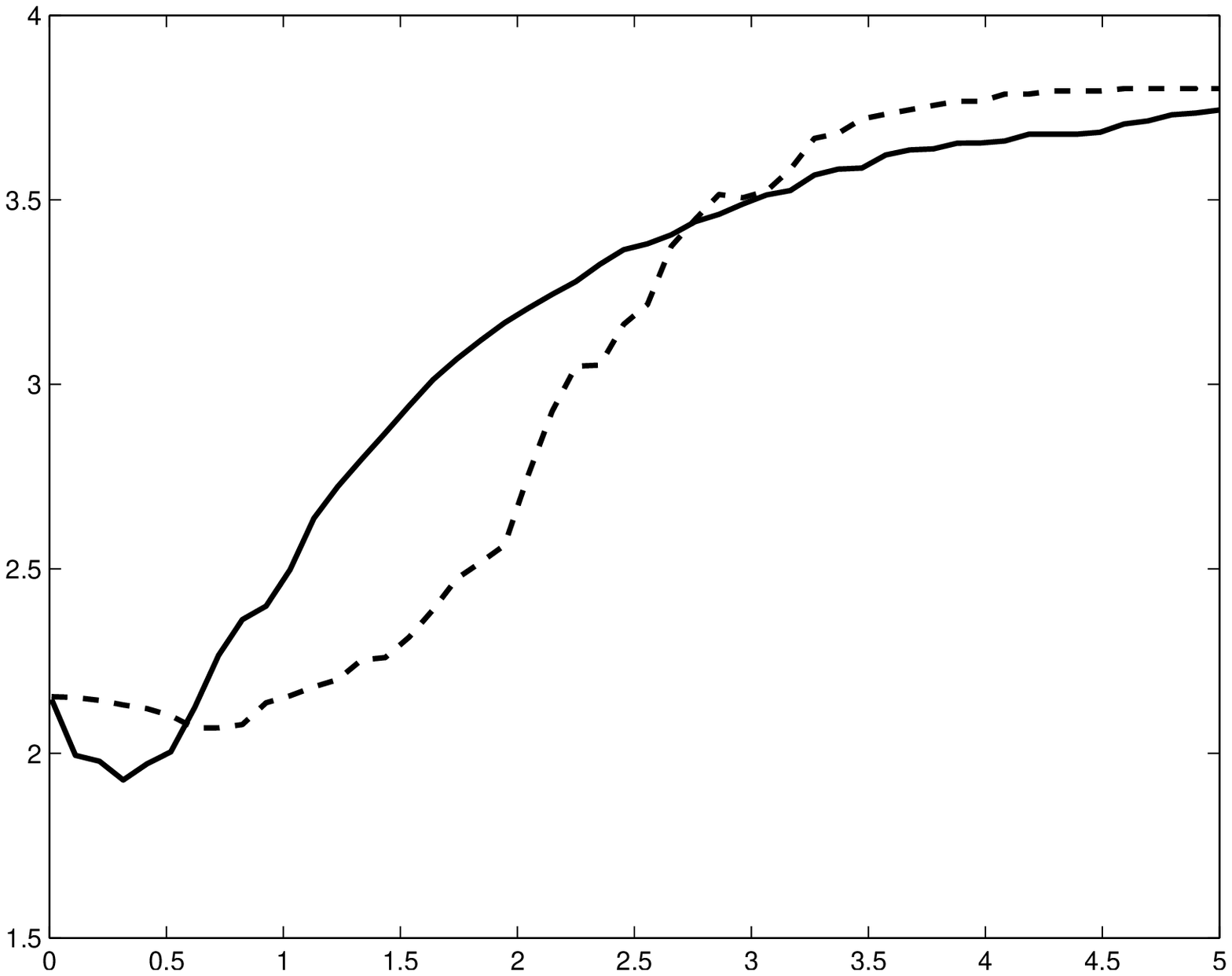} }
\subfigure[MixtGauss - $j_{1} = 5$]
{ \includegraphics[width=3.5cm]{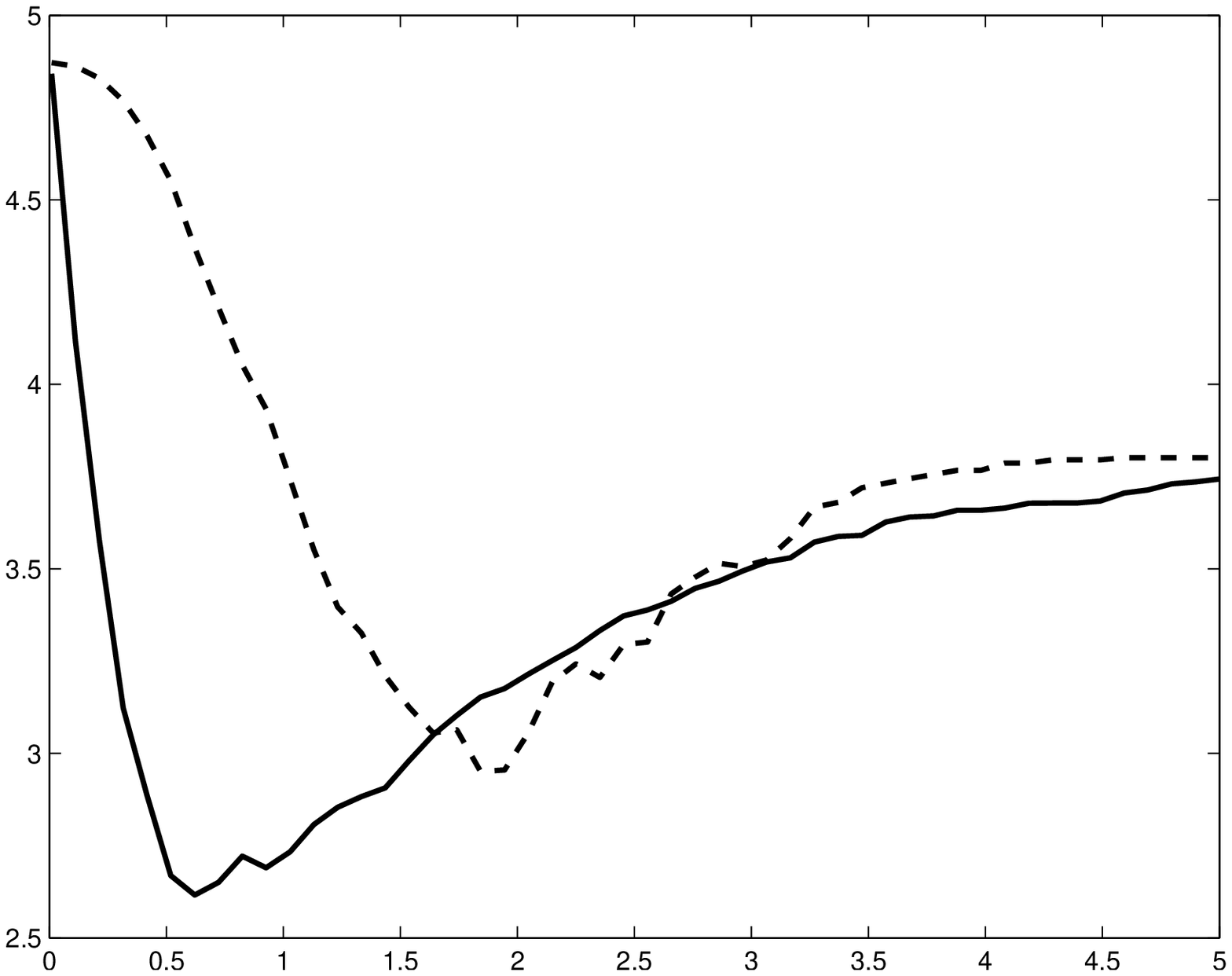} }
\subfigure[MixtGauss - $j_{1} = 6$]
{ \includegraphics[width=3.5cm]{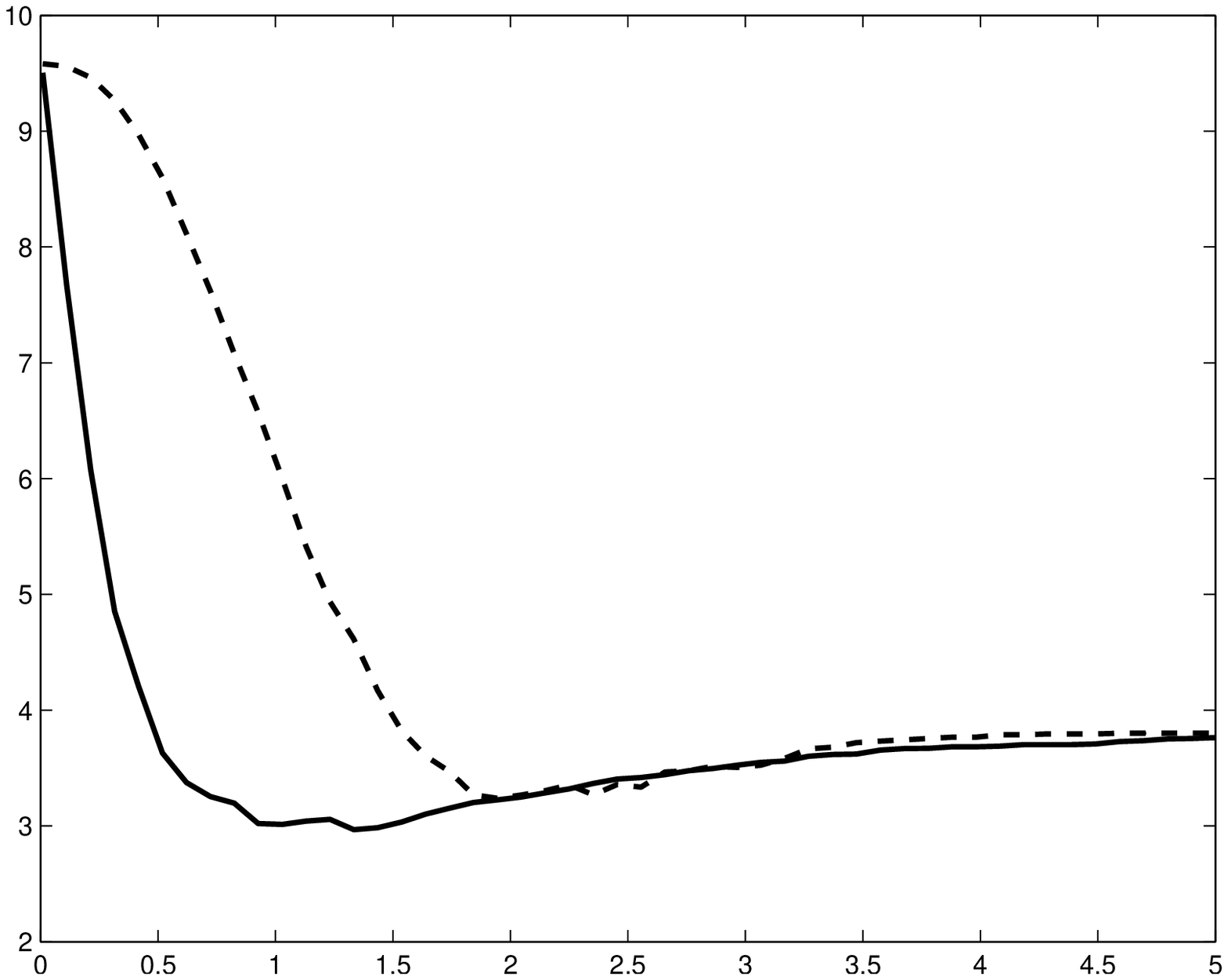} }
%\subfigure[MixtGauss - $j_{1} = 7$]
%{ \includegraphics[width=3.5cm]{MixtGauss_n200_j1_7.eps} }

\caption{Direct density estimation with $n=200$. Evolution of  $R_{n}(j_{1},\delta)$ (solid curves) and $\tilde{R}_{n}(j_{1},\delta)$ (dashed curves)  as a function of $\delta \in [0,5]$  for different values of $j_{1} \geq  j_{1}^{\ast} = \lfloor \frac{1}{2} \log_{2}(n) \rfloor + 1 = 4$.} \label{fig:riskn200}
\end{figure}

Then, for each test density, $M=100$ independent samples of size $n = 100, 200$ are drawn. Empirical average of $R_{n}(j_{1},\delta)$ and $\tilde{R}_{n}(j_{1},\delta)$ over these $M$ repetitions are plotted in Figure \ref{fig:riskn100}  ($n=100$) and Figure \ref{fig:riskn200} ($n=200$) with $\delta \in [0,5]$ and $ j_{1}^{\ast} \leq j_{1} \leq  j_{1}^{\ast} + 2$. For $n = 200$ and $N=256$, we also display in Figure \ref{fig:oracle} the true wavelet coefficients  $\beta_{j,k}$ and the standard deviation  $\sigma_{j,k}$ to show the ideal thresholding performed by the oracle estimator $\tilde{f}_{n}$. For the Uniform, Exponential and Laplace distributions and $j_{1} = 4$, it can be seen that $|\beta_{j,k}| \geq \sigma_{j,k}$ for all $j_{0} \leq j \leq j_{1}$ and $0 \leq k \leq 2^{j}-1$ , which means that the oracle estimator $\tilde{f}_{n}$ performs no thresholding and keeps all empirical wavelet coefficients $\beta_{j,k}$  for $j \leq j_{1}$. For the MixtGauss distribution and $j_{1} = 4$ the behavior of the  oracle estimator is different as for some $0 \leq k \leq 2^{j_{1}}-1$, one can observe that $|\beta_{j,k}| < \sigma_{j,k}$. One retrieves this behavior in the first column of  Figure \ref{fig:riskn200} (case $j_{1} = 4$) for the curves $R_{n}(j_{1},\delta) $ and  $\tilde{R}_{n}(j_{1},\delta)  $ which both have a minimum at $\delta = 0$ (no thresholding is done)  for the Uniform, Exponential and Laplace distributions. For the MixtGauss distribution,   $R_{n}(j_{1},\delta) $ has a minimum at $\delta \approx 0.4$ while  $\tilde{R}_{n}(j_{1},\delta) $ has a minimum at   $\delta \approx 0.8$. For $5 \leq j_{1} \leq 7$ and  the  Uniform, Laplace and MixtGauss distributions,   $R_{n}(j_{1},\delta) $ has a minimum at some $\delta \in [0.4,1]$ and the value of  $R_{n}(j_{1},\delta) $ at this point is smaller that the minimum of $\tilde{R}_{n}(j_{1},\delta)$. This indicates that taking thresholds depending on the location $(j,k)$ can improve the quality of the estimation. For the Exponential distribution, the estimator with a level-dependent threshold of the form $\delta \sqrt{j/n}$ performs generally better than $\hat{f}_{n}$.  Similar comments can be made for the behavior of the curves displayed in Figure \ref{fig:riskn100} (case $n = 100$).

Moreover, it can been that the smallest value of the risk of $\hat{f}_{n}$ relative to the risk of the oracle  $\tilde{f}_{n}$ are obtained for  $j_{1} = j_{1}^{\ast} = \lfloor \frac{1}{2} \log_{2}(n) \rfloor + 1$. This indicates that introducing wavelet coefficients at resolution level larger than $j_{1}^{\ast}$ generally deteriorates the quality of the estimation.  Finally, note that the curves in Figure \ref{fig:riskn100} and Figure \ref{fig:riskn200} tends to confirm that the choice (\ref{eq:choice}) for $j_{1}$ and $\delta$ is reasonable, and leads to very satisfactory estimators.

\subsection{Density deconvolution}

In the case of density deconvolution, we propose to compare the performances of our wavelet approach with those of the adaptive density deconvolution estimator of \citeA{CRT06jscs,CRT06cjs} that is based on penalized contrast minimization over a collection of models containing square integrable functions with Fourier transform having a compact support. Such an estimator is therefore a band-limited function, see \citeA{CRT06cjs} for further details. Moreover, this method can be viewed as a kind of adaptive linear wavelet estimator with a Shannon wavelet basis which is also a band-limited function like the Meyer wavelet but less localized in the time domain. \citeA{CRT06cjs} have shown that the model selection procedure performs very well for finite samples, compared with other standard estimators. In particular, this estimator outperforms the kernel estimator, even when the bandwidth parameter is selected in a data-driven way which makes this procedure as the most challenging competitor in our simulations.

Observations $Y_i, i=1,\ldots,n$ are generated from the additive model
	$
	Y_{i} = X_{i} + \epsilon_{i}
	$, 
where $X_{i}$ are independent realizations from one of the test functions $f$ displayed in Figure \ref{figdens}, and the $\epsilon_{i}$'s  are i.d.d. additive errors with density $h$. Results are presented with a Laplace measurement error, that is 
	$
	h(x) = (\sqrt{2}\sigma_\epsilon)^{-1} \exp ( -\sqrt{2} |x|/\sigma_\epsilon )$, for $x \in \RR$ and where $\sigma_\epsilon$ is the standard deviation of the errors.  The Fourier coefficients of $h$ are given by $h_{\ell} = (1+2 \sigma^{2}_\epsilon\pi^{2} \ell^{2})^{-1}$, $\ell \in \ZZ$ which thus corresponds to the case of ordinary smooth deconvolution with $\nu = 2$. The main quantities in the simulations are the sample size $n$ and the root signal-to-noise ratio defined by $s2n = \sigma_{X}/\sigma_\epsilon$ with $\sigma_{X} = \sqrt{\var(X_{1})}$.
	
	According to Theorem \ref{theo:oracledec}, one can take $j_{0} = \lfloor \log_{2}(\log(n)) \rfloor + 1$ and important quantities to control the quality of estimation by wavelet thresholding are the highest resolution level $j_{1}$ and the tuning parameter $\delta$. If $j_{1}$ is such that $2^{j_{1}} > n^{\eta/(\nu +1)} (\log n)^{\alpha} \geq 2^{j_{1}-1}$ for some $\eta > 0$, $\alpha \geq 0$, then Theorem \ref{theo:oracledec} suggests to take $\delta = \eta \left(1 + \frac{\nu}{\nu +1}\right)$. As already remarked, for ill-posed inverse problems, a smaller $j_{1}$ than in the direct case should be used. However, for $n=200$, following the asymptotic considerations in Section \ref{sec:minimax}, the choices $\eta = 1/2$, $\alpha = 0$ and $\nu =2$ yield to $j_{1} = \lfloor \frac{1}{6} \log_{2}(n) \rfloor + 1 = 2$ which is smaller than $j_{0} =  \lfloor \log_{2}(\log(n)) \rfloor + 1 = 3$. Hence, setting in advance values for $\eta$ and $\alpha$ may yield a theoretical choice for $j_{1}$ that cannot be used in practice. Note that this issue has been noticed in several papers on deconvolution by wavelet thresholding, see  \citeA{JKPR04jrssb}, \citeA{BVB07}.
	
	 Alternatively, let us argue as in the direct case, by assuming that $j_{1}$ is given.  If one sets $\alpha = 0$, then for choosing $\delta$ one can take $\eta^{\ast} = (\nu +1)(j_{1}-1)/\log_{2}(n)$ which is the smallest constant $\eta$ that satisfies $2^{j_{1}} > n^{\eta/(\nu +1)}  \geq 2^{j_{1}-1}$, and then take $\delta  = \left(1 + \frac{\nu}{\nu +1}\right)\eta^{\ast}  =  \left( 2\nu + 1 \right) (j_{1}-1)/\log_{2}(n)  $. A smaller value for $\delta$ can also be made, by choosing $\alpha \neq 0$ and by taking $\delta  = \left(1 + \frac{\nu}{\nu +1}\right)\eta^{\ast}$ with $\eta^{\ast} = (\nu +1)(j_{1}-1-\alpha \log_{2}(\log(n)) )/\log_{2}(n)$.
	
\begin{figure}[h!]
\centering
\subfigure
{ \includegraphics[width=3.5cm]{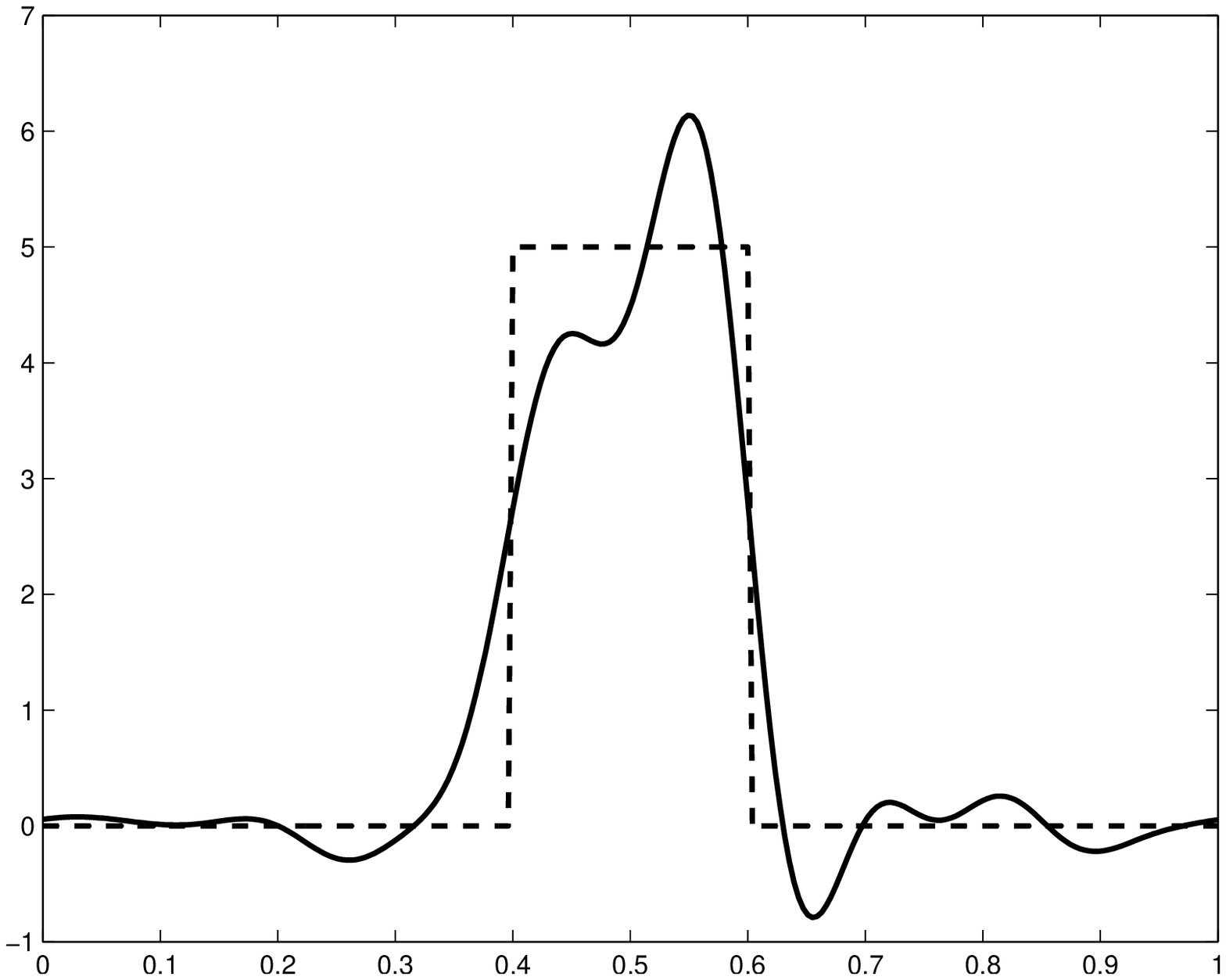} }
\subfigure
{ \includegraphics[width= 3.5cm]{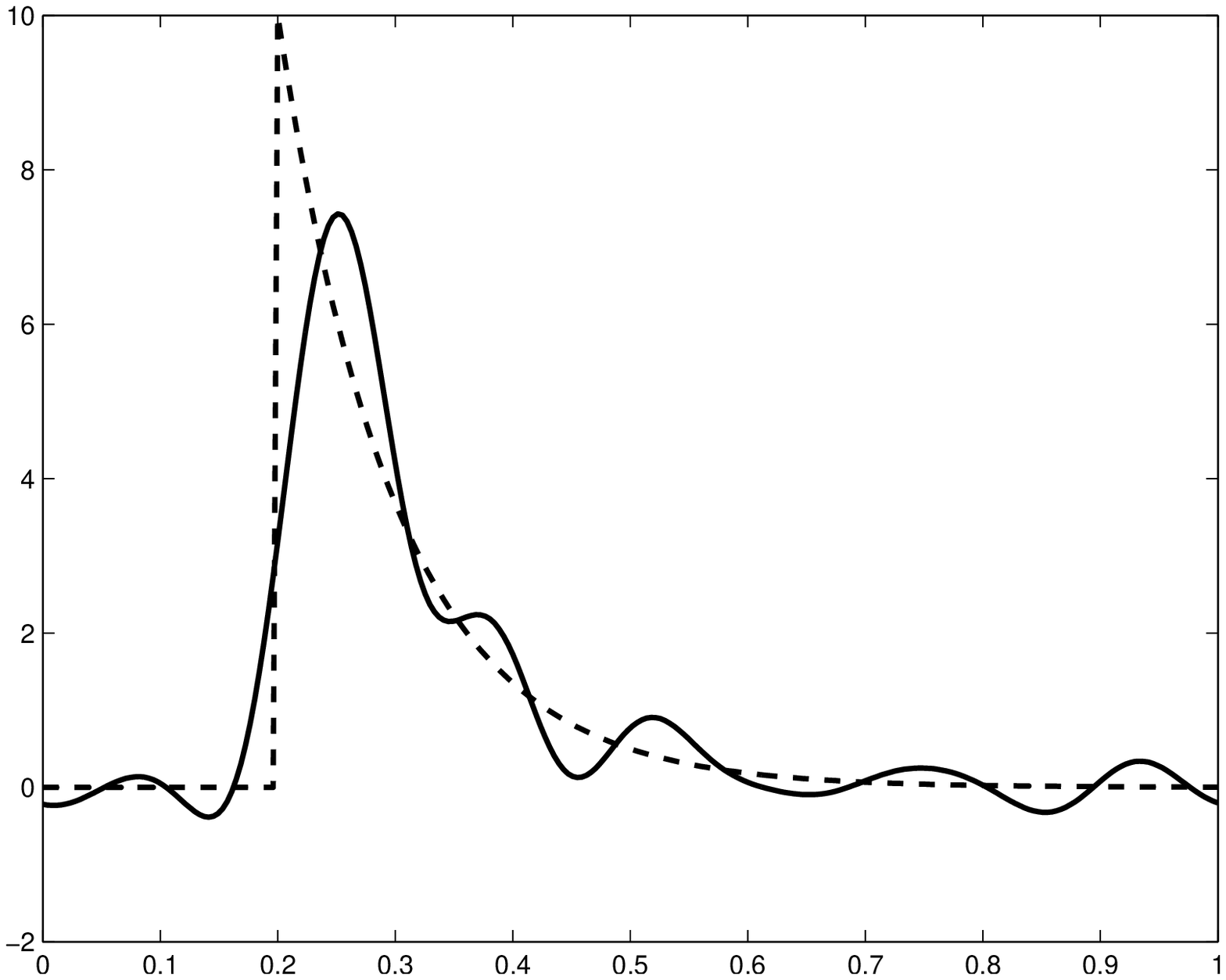} }
\subfigure
{ \includegraphics[width= 3.5cm]{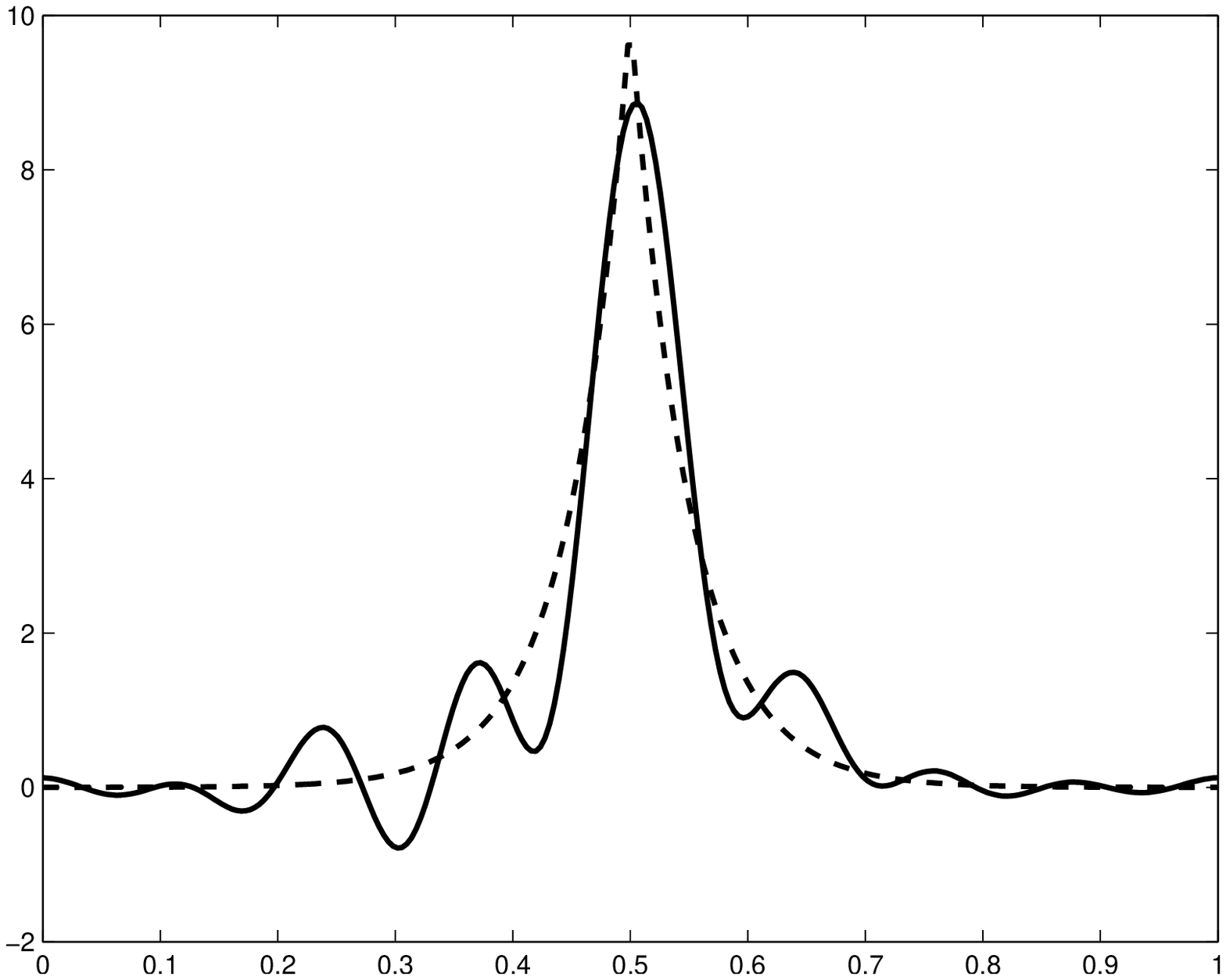} }
\subfigure
{ \includegraphics[width= 3.5cm]{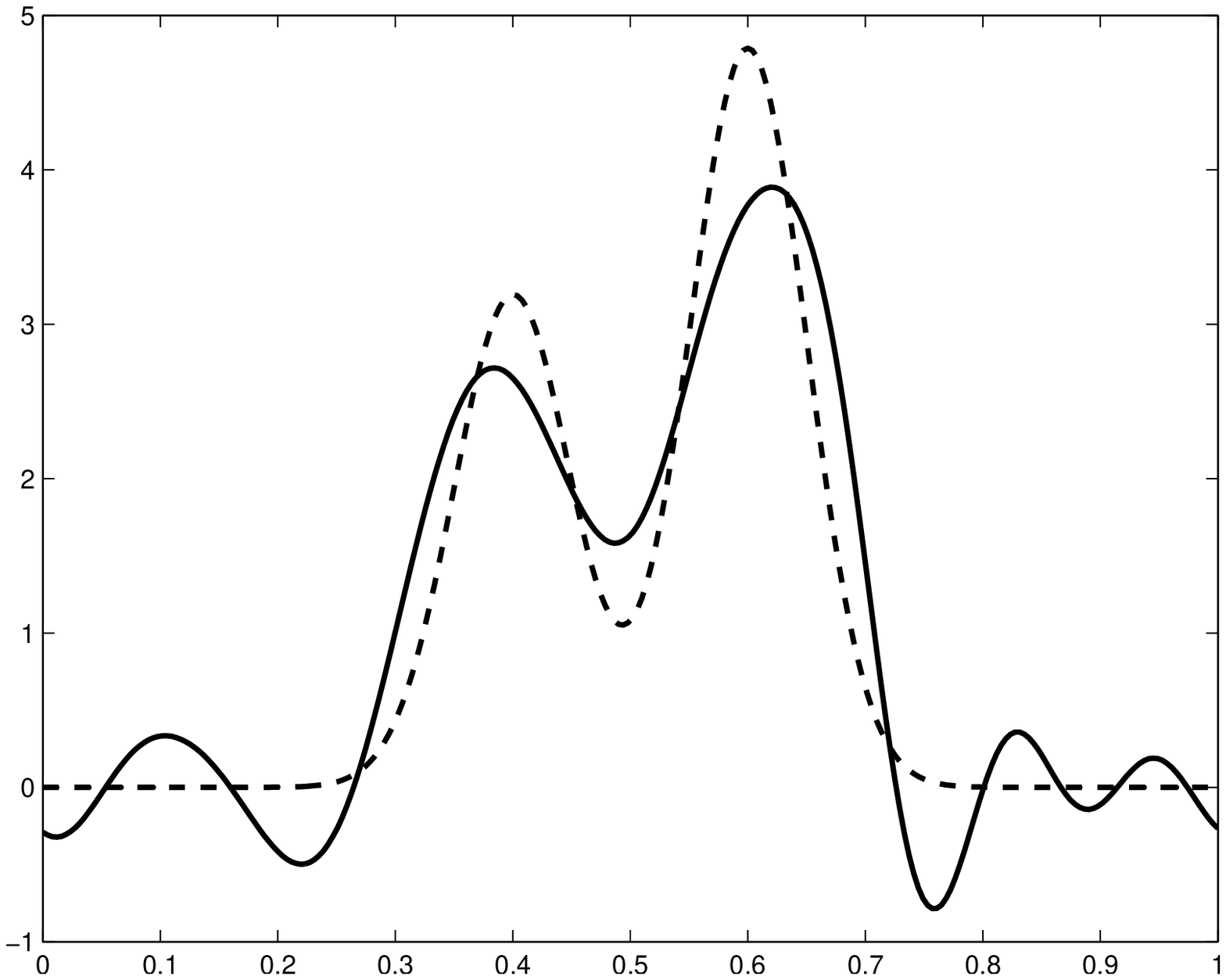} }

\subfigure
{ \includegraphics[width= 3.5cm]{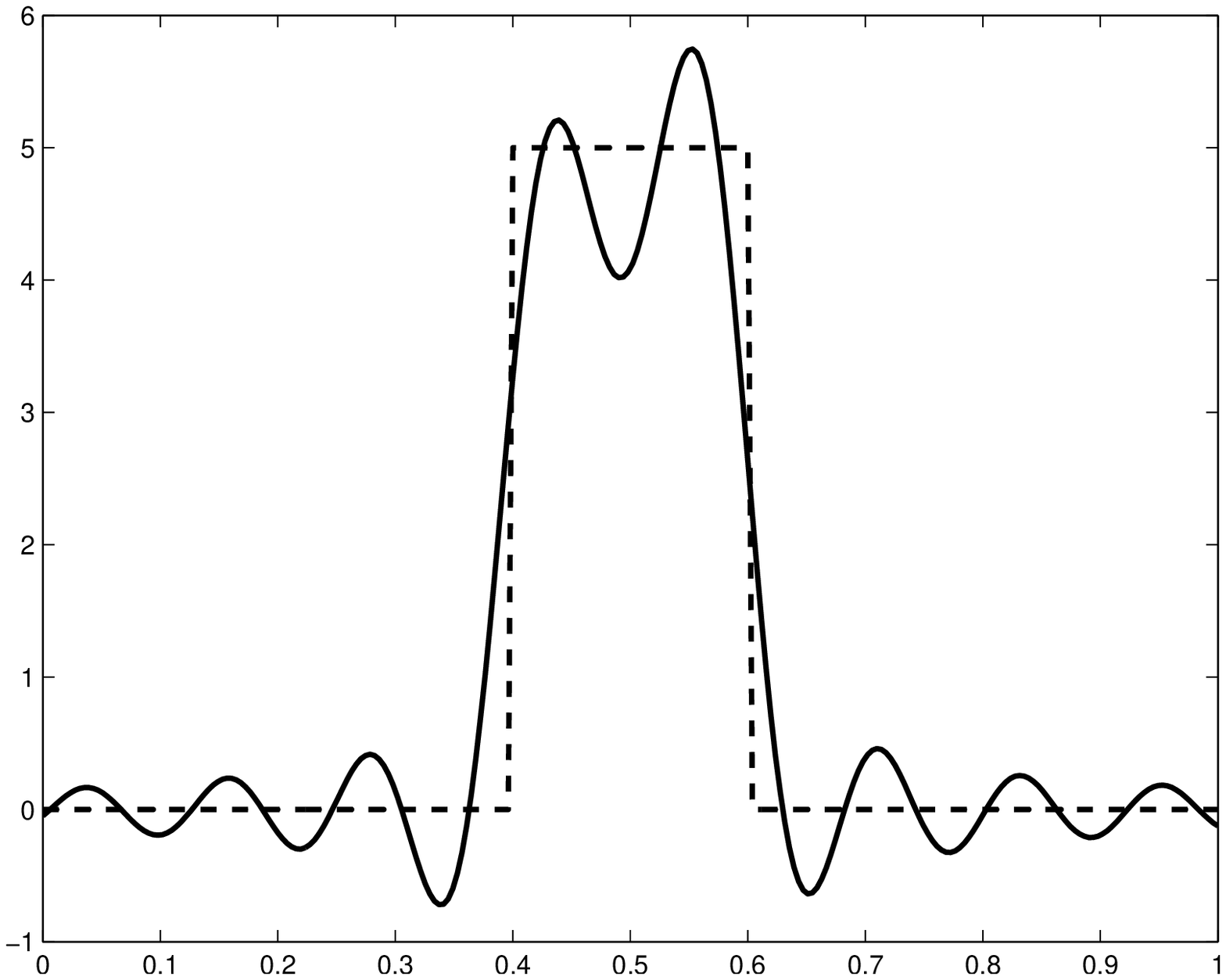} }
\subfigure
{ \includegraphics[width= 3.5cm]{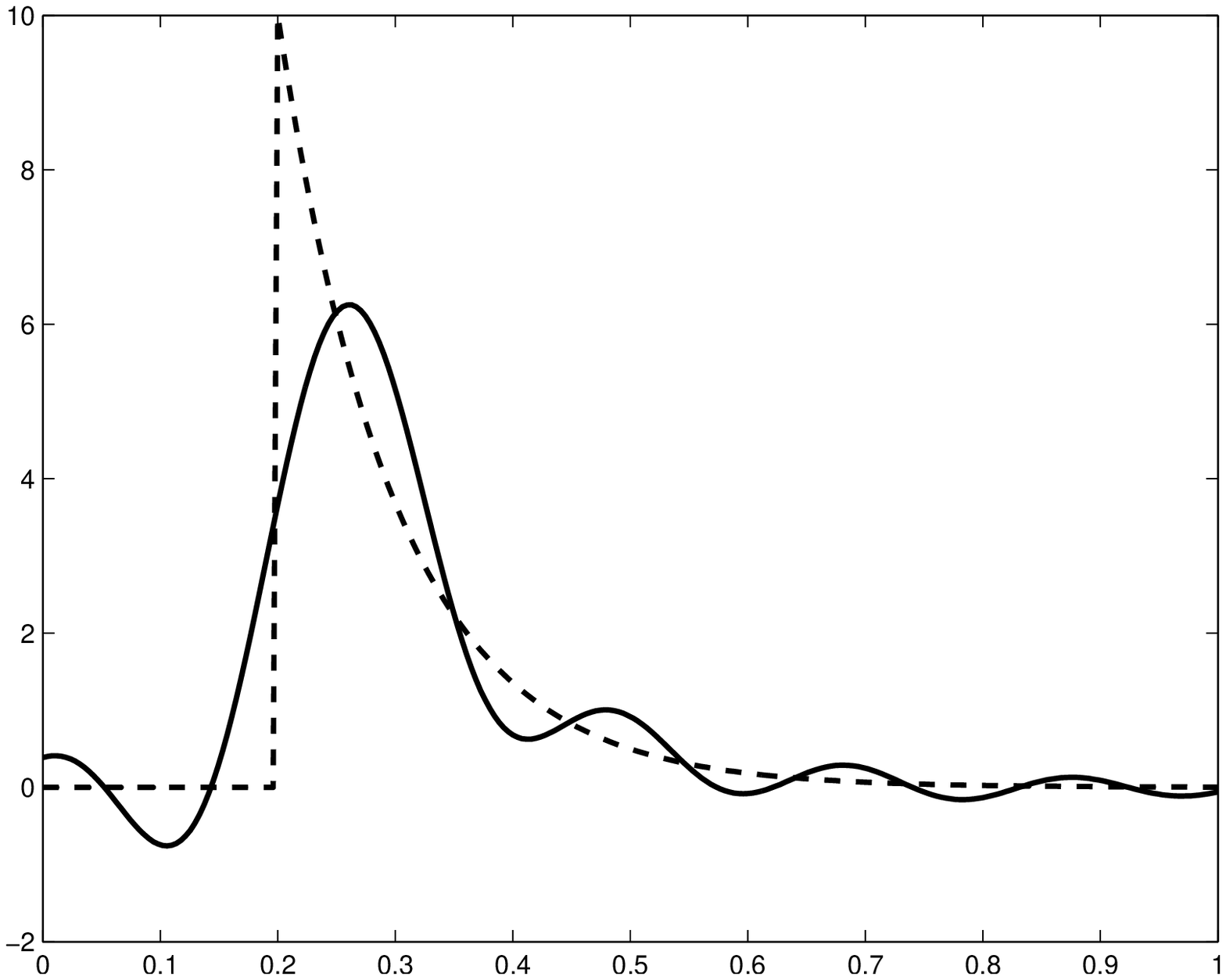} }
\subfigure
{ \includegraphics[width= 3.5cm]{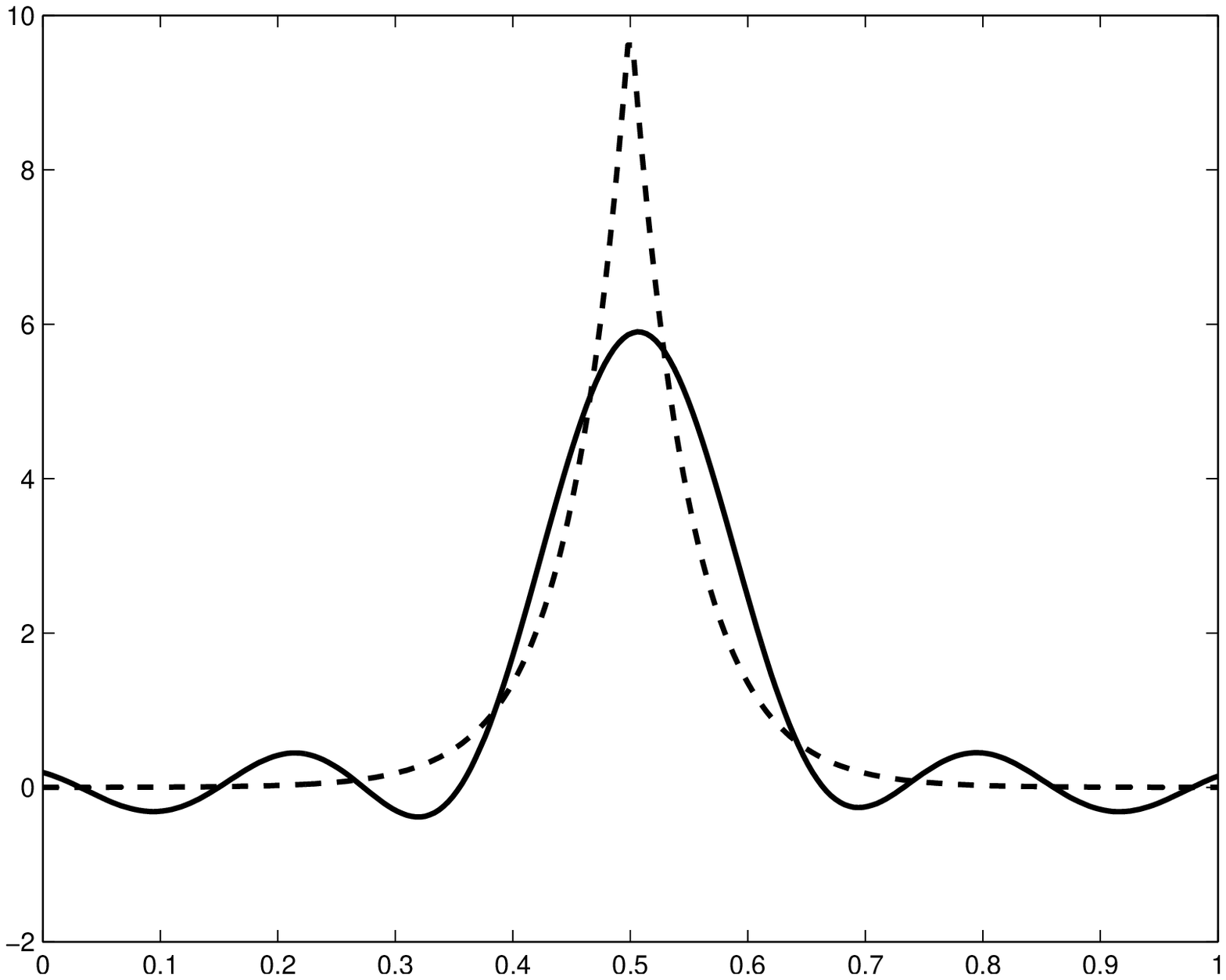} }
\subfigure
{ \includegraphics[width= 3.5cm]{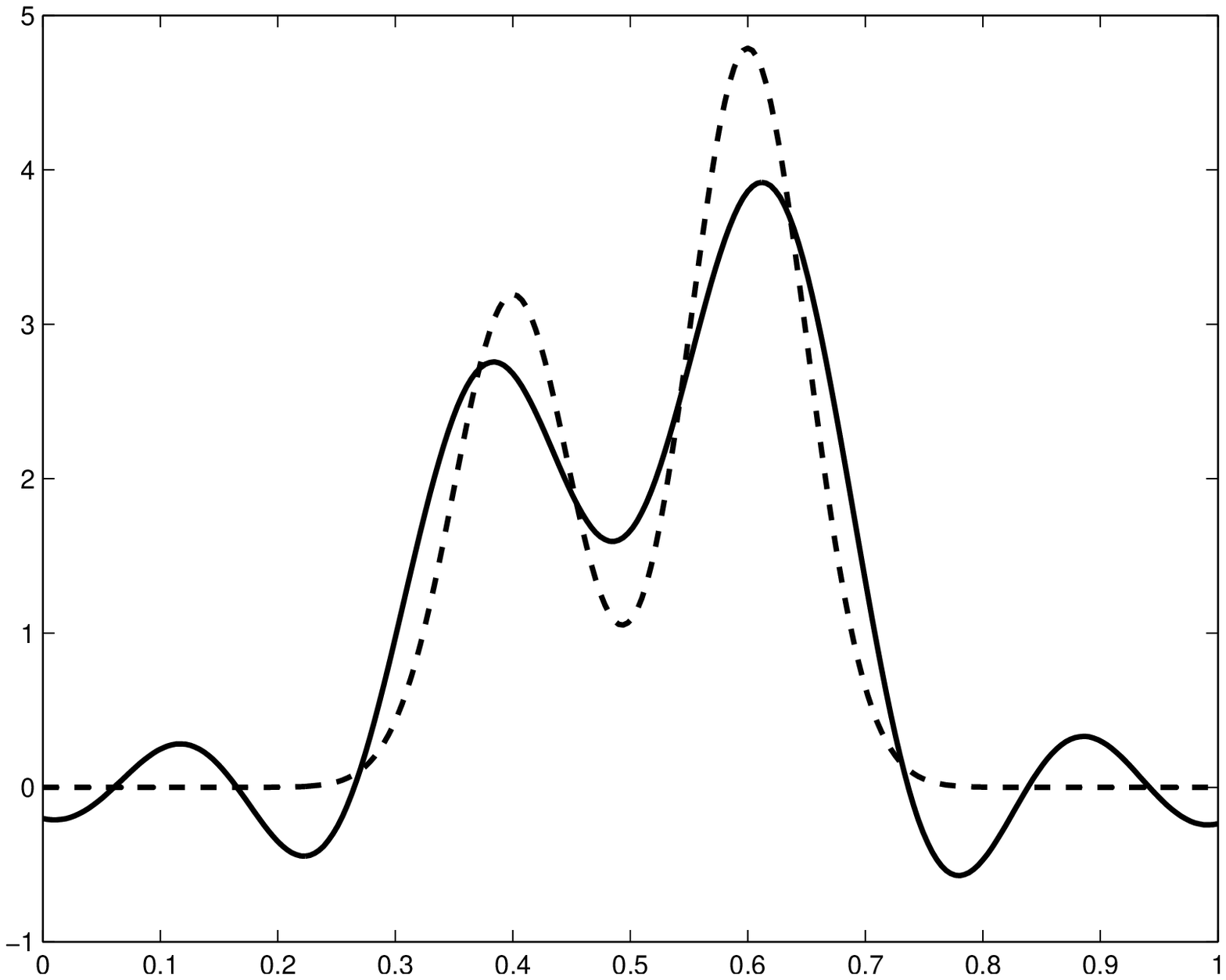} }

\caption{Typical reconstructions from a single simulation with  $n=200$ for the four test densities Uniform (1st column), Exponential  (2nd column),  Laplace  (3rd column), MixtGauss  (4th column) contaminated with Laplace noise, by wavelet thresholding (1st row) with $j_{1} = j_{0}  = 3$, $\delta \approx 0.5215$, and model selection (2nd row). } \label{fig:exampledec}
\end{figure}

	We report results for $n=100, 200$ and $s2n = 3$, which a relatively large signal-to-noise ratio.  To give an idea of the quality of $\hat{f}_{n}$ and to compare it with the model selection estimator, a typical example of estimation with $n=200$, $j_{1} = j_{0}  = 3$, $\alpha = 0.5$ and $\delta =  \left(2 \nu + 1 \right) (j_{1}-1-\alpha \log_{2}(\log(n)))/\log_{2}(n)   \approx 0.5215$  is given in  Figure \ref{fig:exampledec}. Both methods perform similarly for the estimation of the smooth density MixtGauss. For the three non-smooth densities, wavelet thresholding performs much better than the model selection estimator. This comparison on a single simulation tends to confirm the superiority of wavelet-based methods over those based on Fourier decompositions for the reconstruction of signals with local singularities.

	However, it is not clear how to choose $j_{1}$  in practice. Indeed, an optimal theoretical level can be too small, but taking a too high  level of resolution may introduce some instability  in the estimation. Moreover, it is of interested to study if $\delta$ can be smaller that $ \eta \left(1 + \frac{\nu}{\nu +1} \right)$. A goal of the following simulation study is thus to identify a reasonable empirical range of values for $j_1$ and $\delta$.
		
\begin{figure}[htbp]
\centering
\subfigure[Uniform - $j_{1} = 3$]
{ \includegraphics[width=3.5cm]{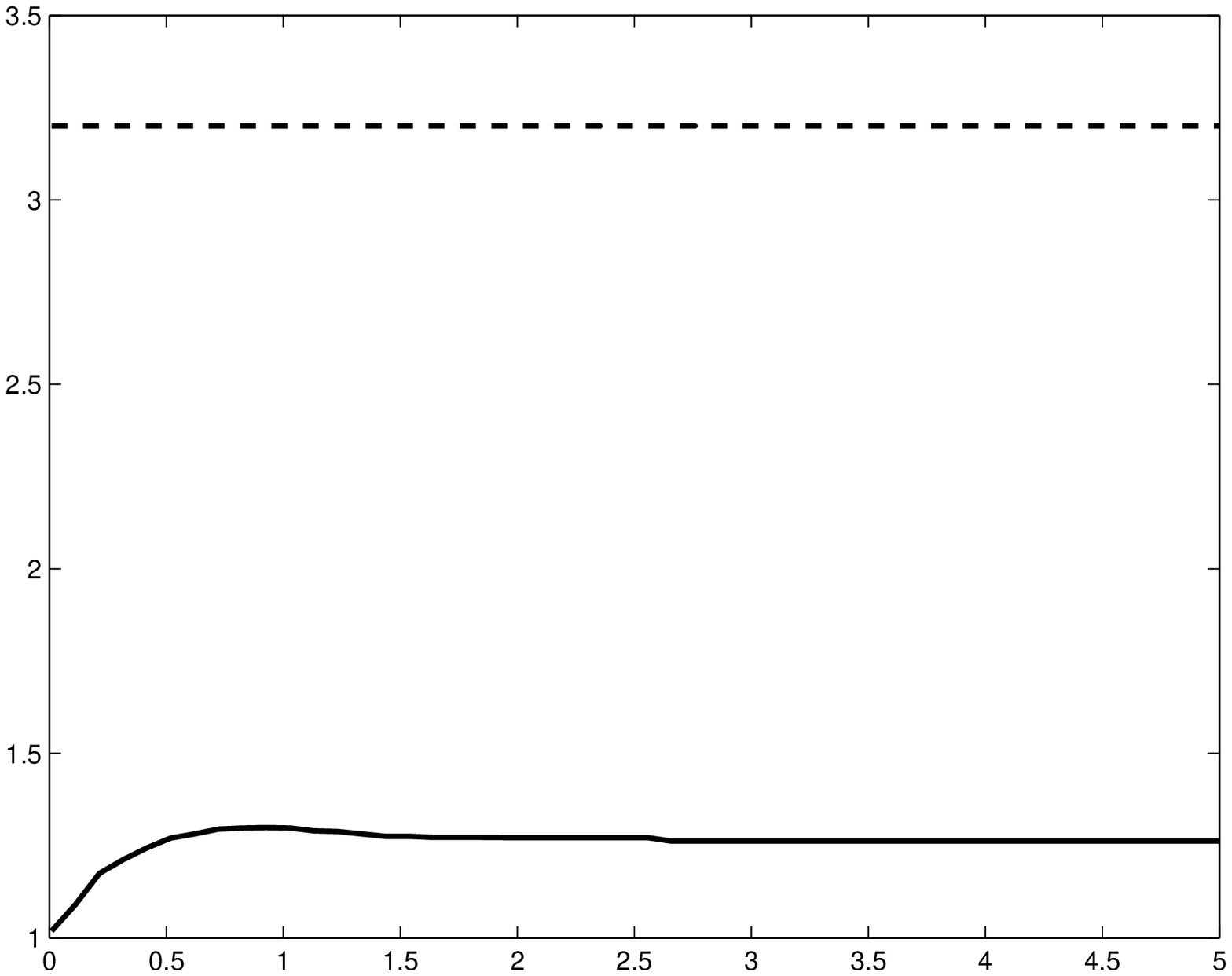} }
\subfigure[Uniform - $j_{1} = 4$]
{ \includegraphics[width=3.5cm]{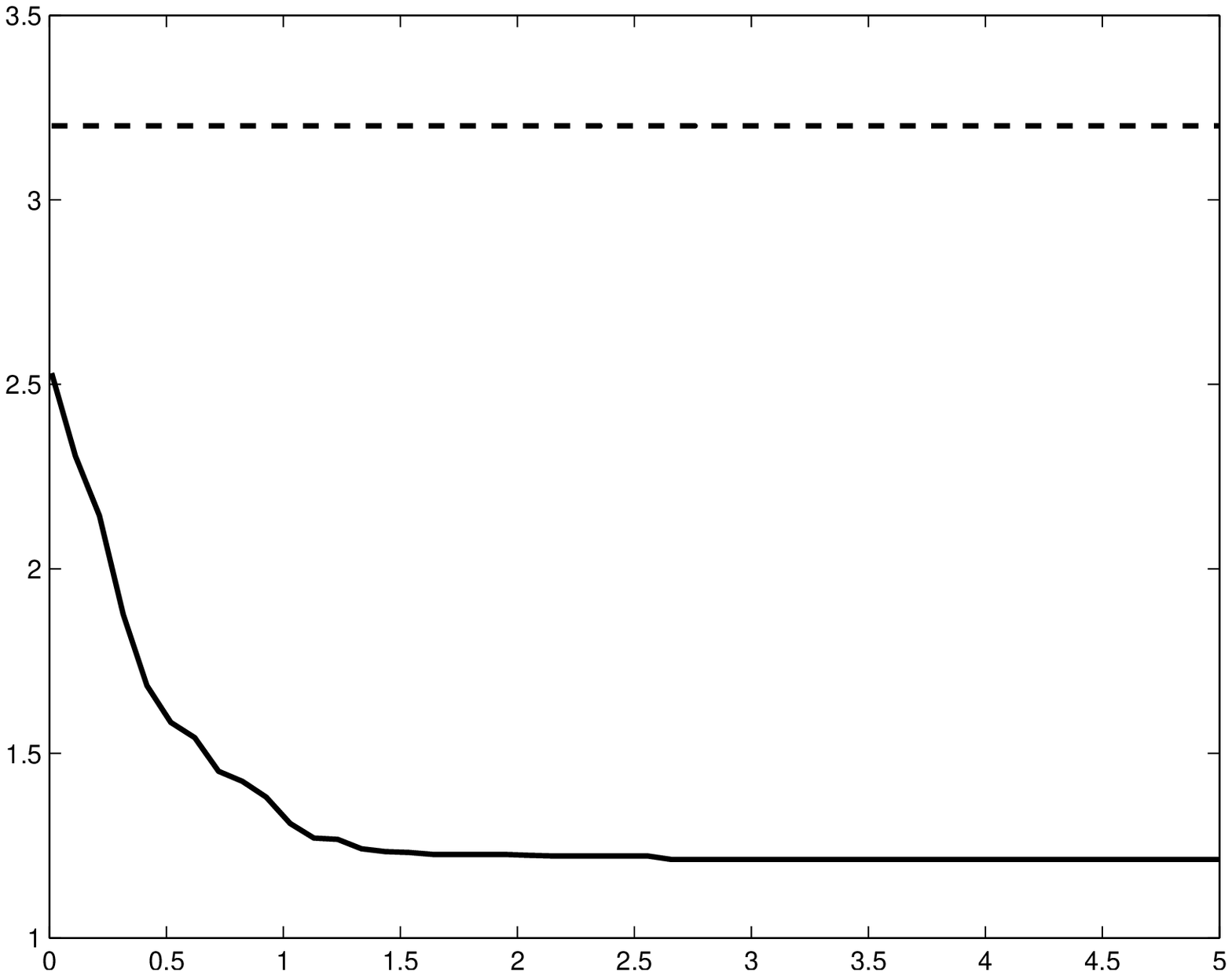} }
\subfigure[Uniform - $j_{1} = 5$]
{ \includegraphics[width=3.5cm]{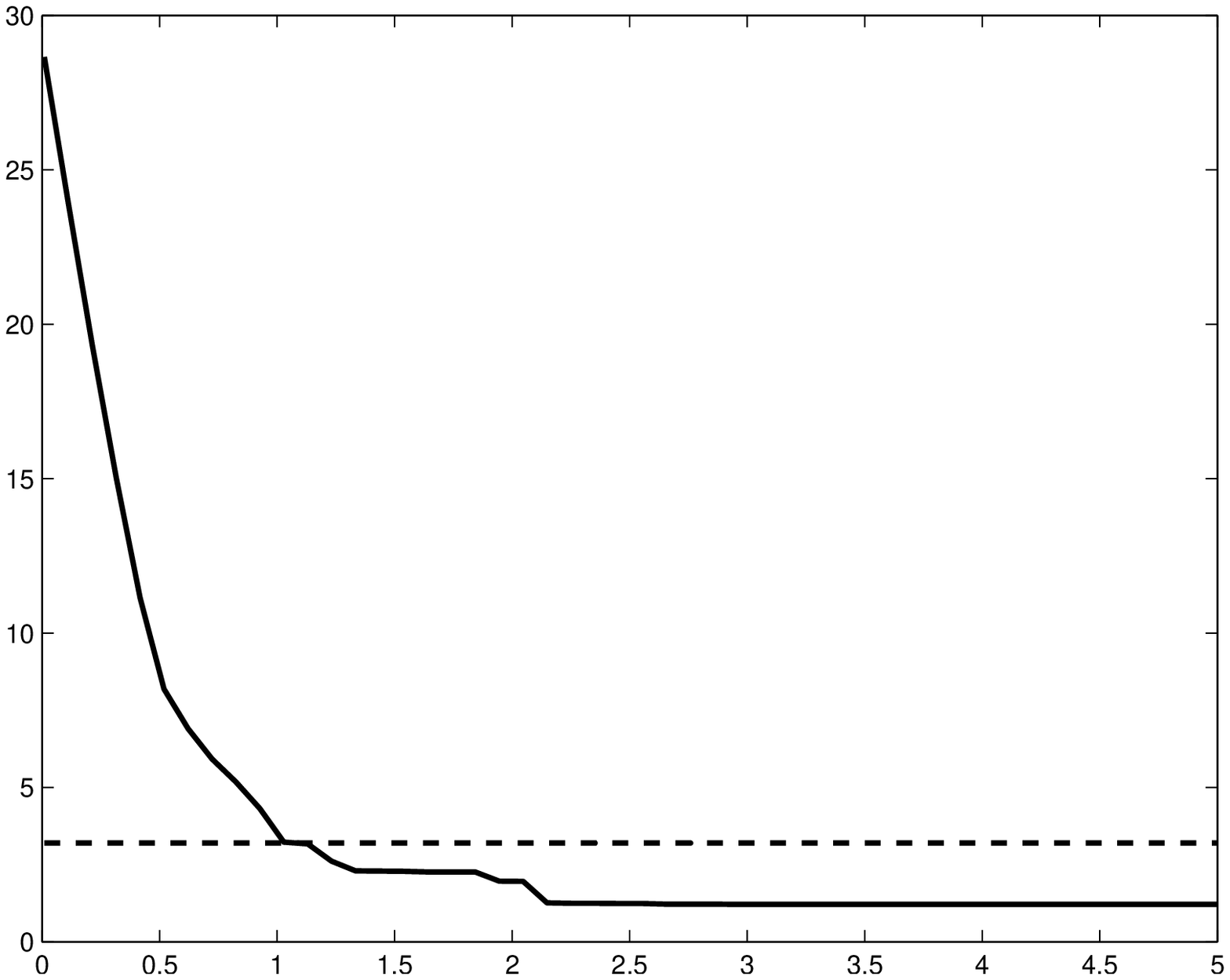} }

\subfigure[Exponential - $j_{1} = 3$]
{ \includegraphics[width=3.5cm]{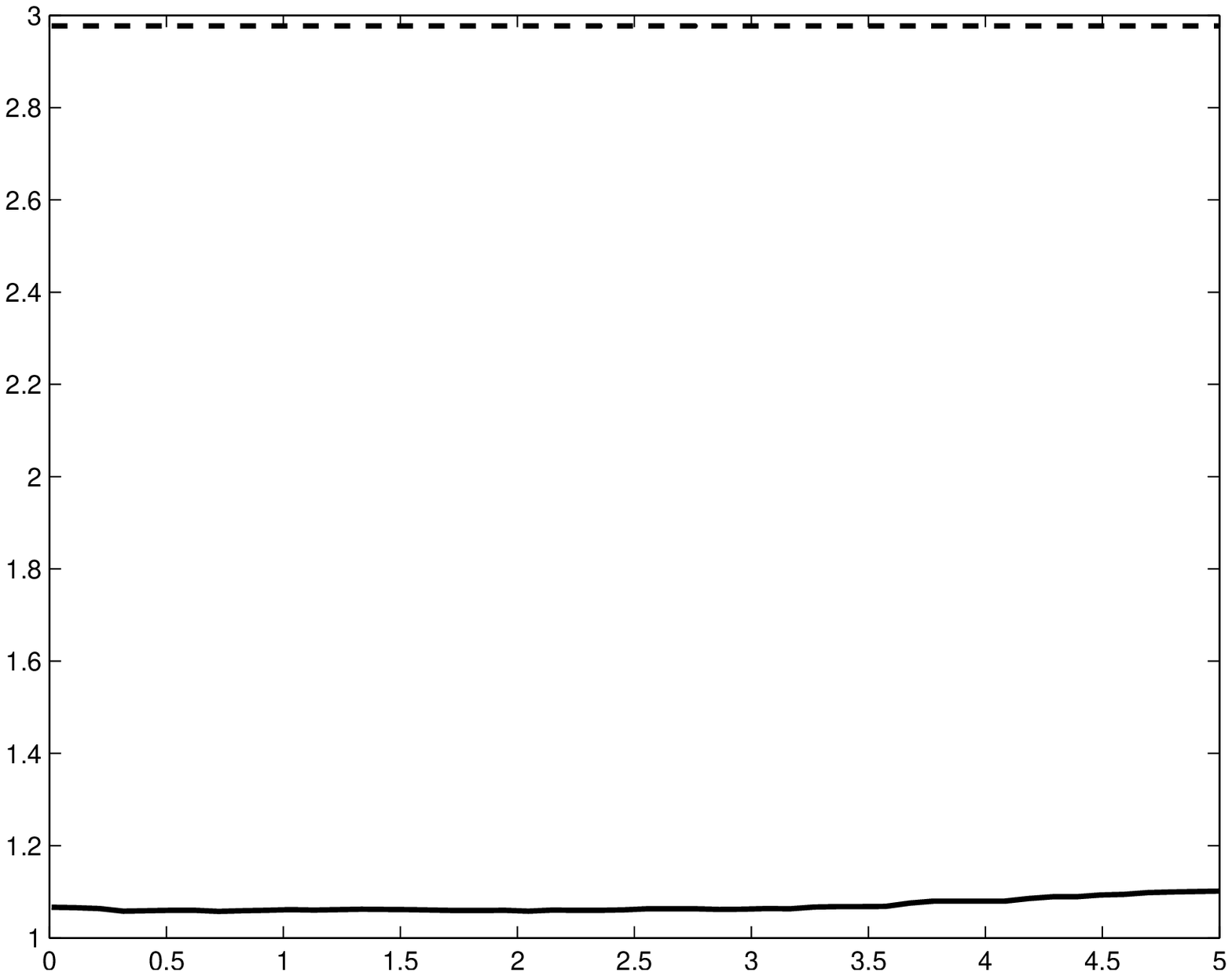} }
\subfigure[Exponential - $j_{1} = 4$]
{ \includegraphics[width=3.5cm]{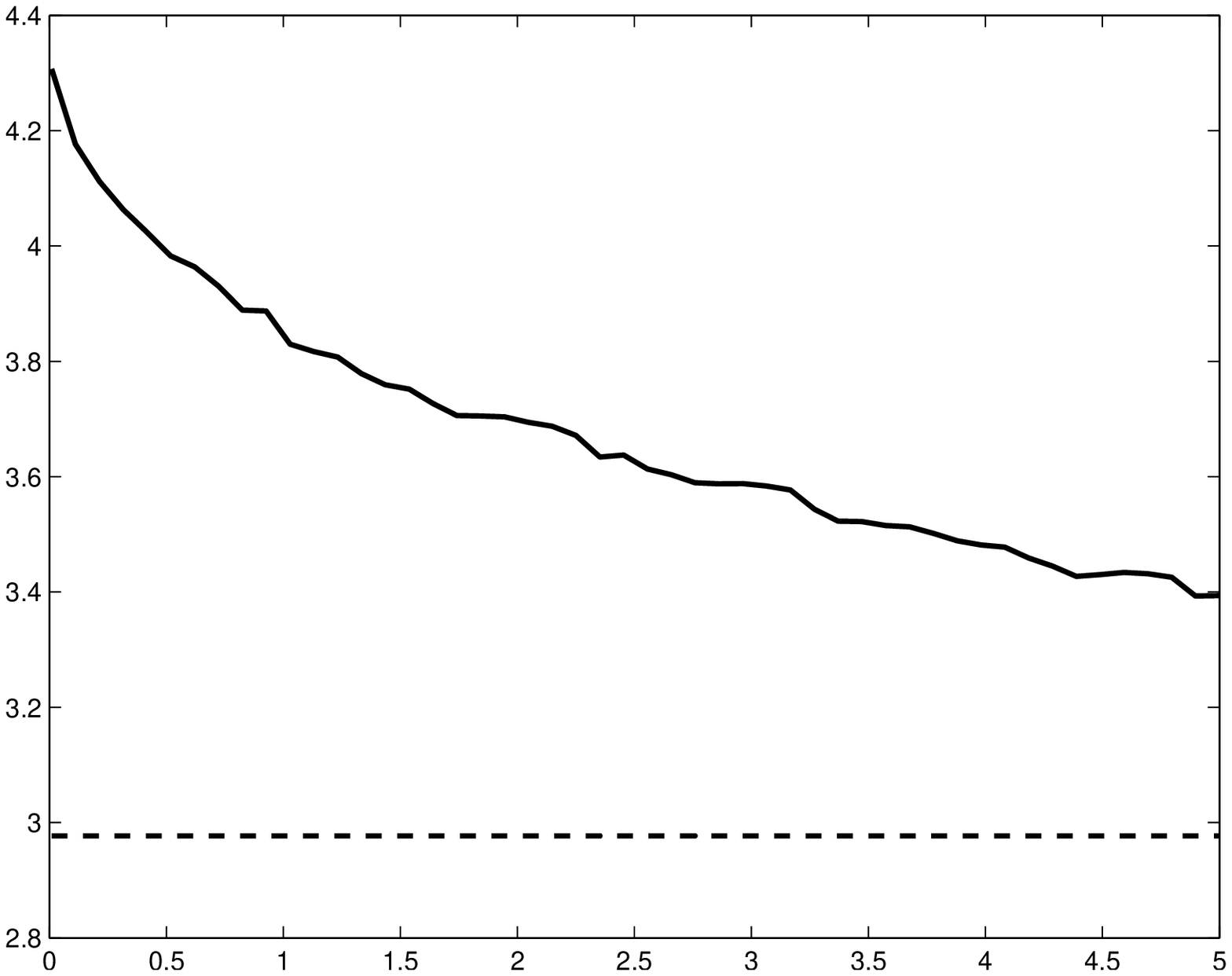} }
\subfigure[Exponential - $j_{1} = 5$]
{ \includegraphics[width=3.5cm]{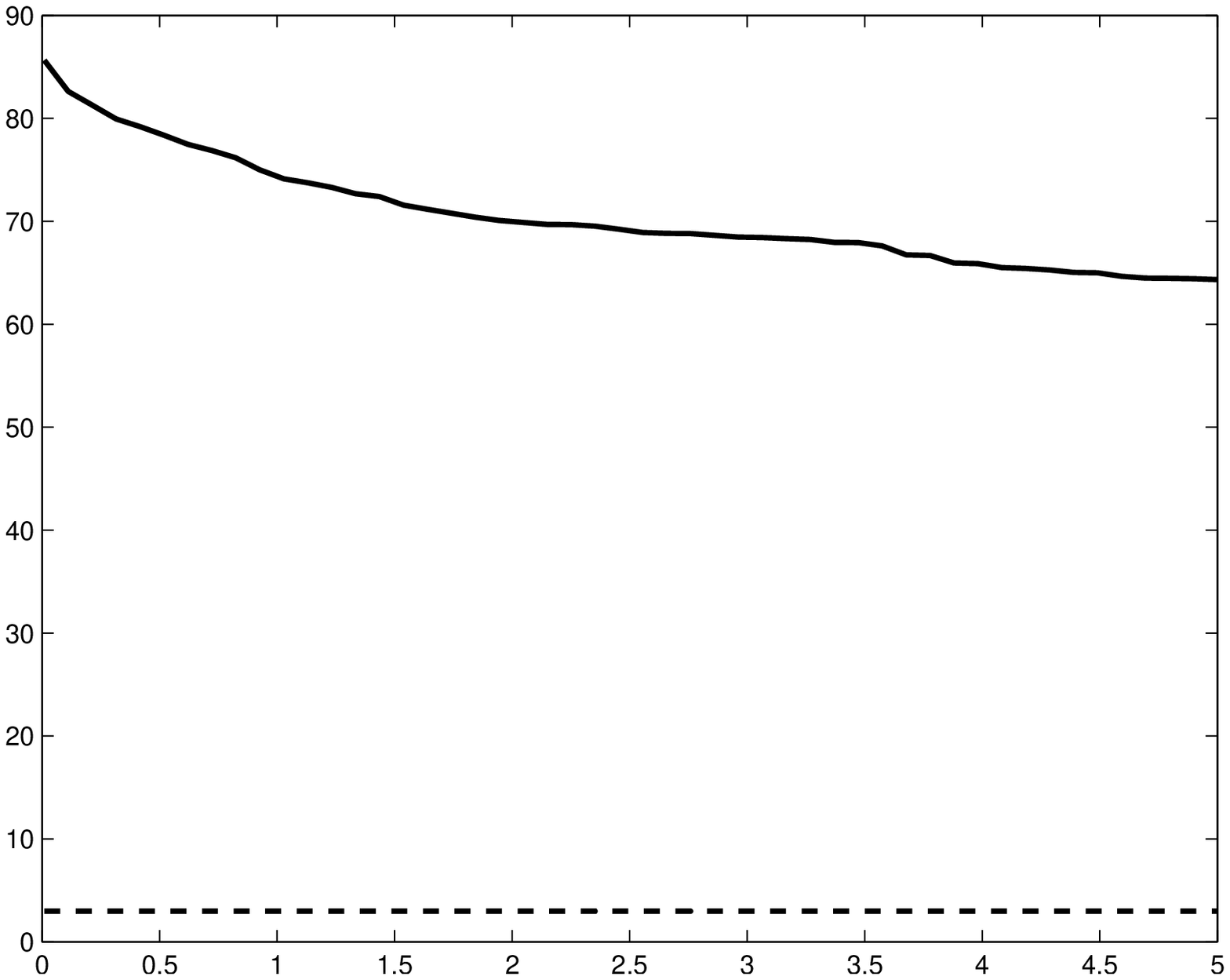} }

\subfigure[Laplace - $j_{1} = 3$]
{ \includegraphics[width=3.5cm]{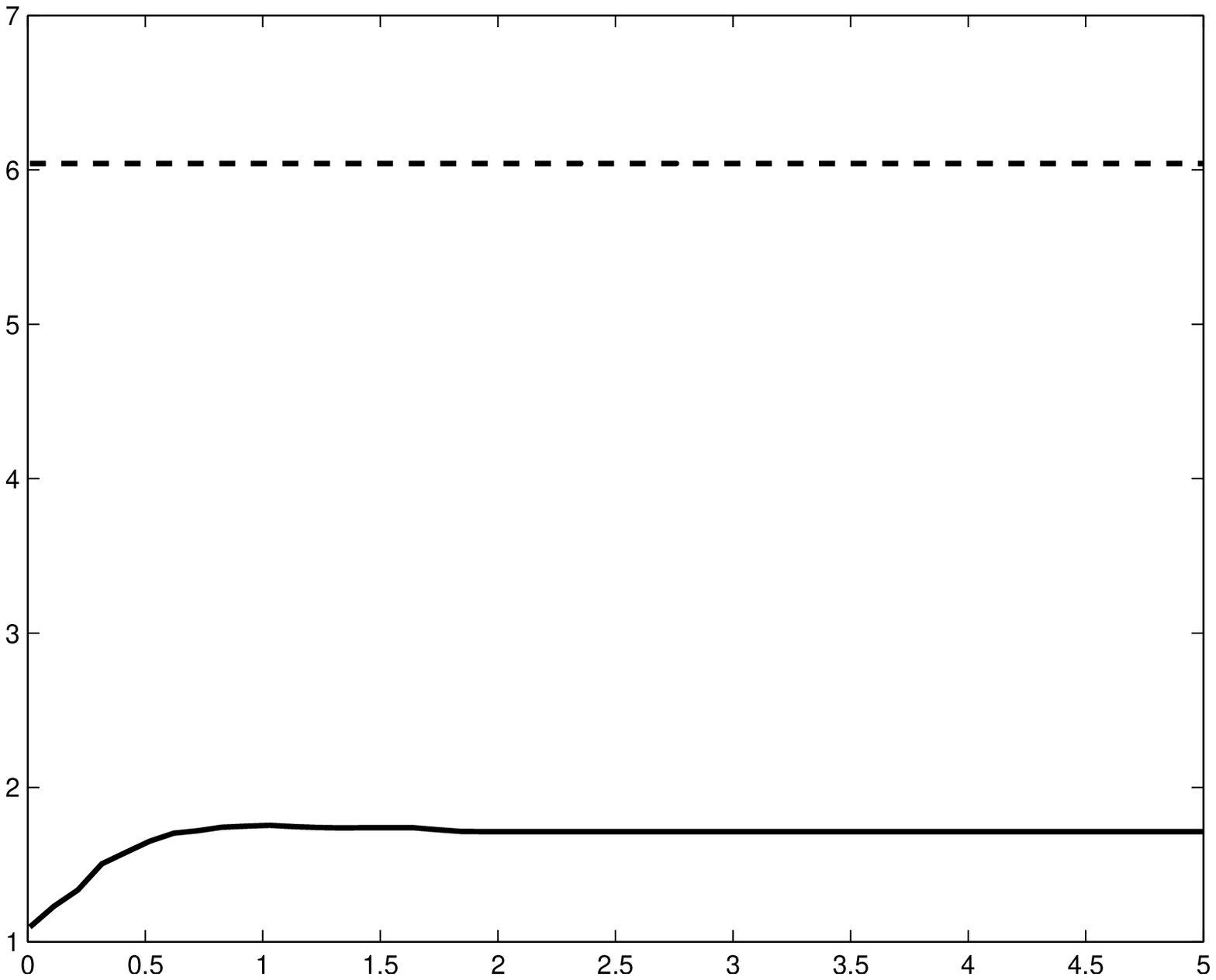} }
\subfigure[Laplace - $j_{1} = 4$]
{ \includegraphics[width=3.5cm]{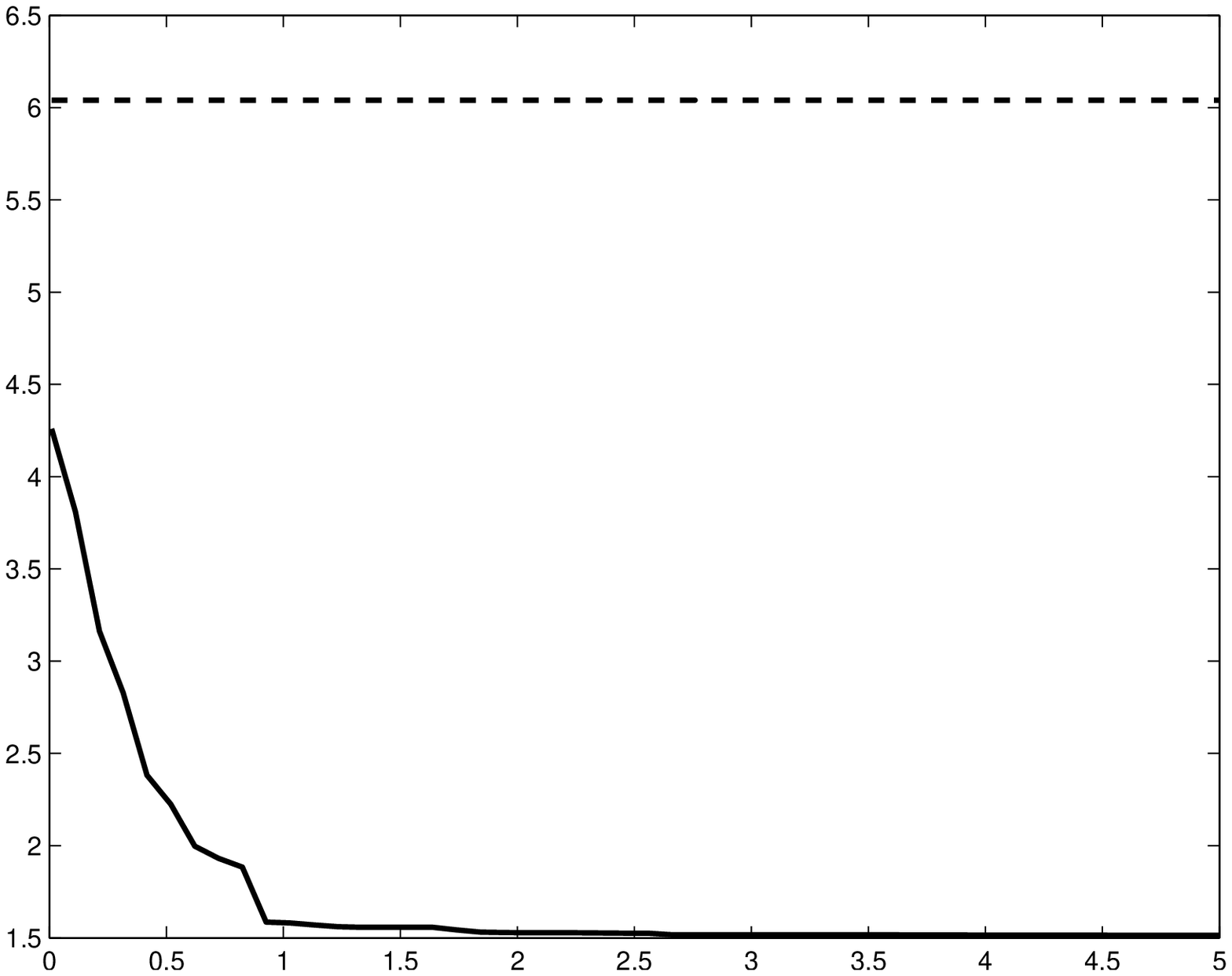} }
\subfigure[Laplace - $j_{1} = 5$]
{ \includegraphics[width=3.5cm]{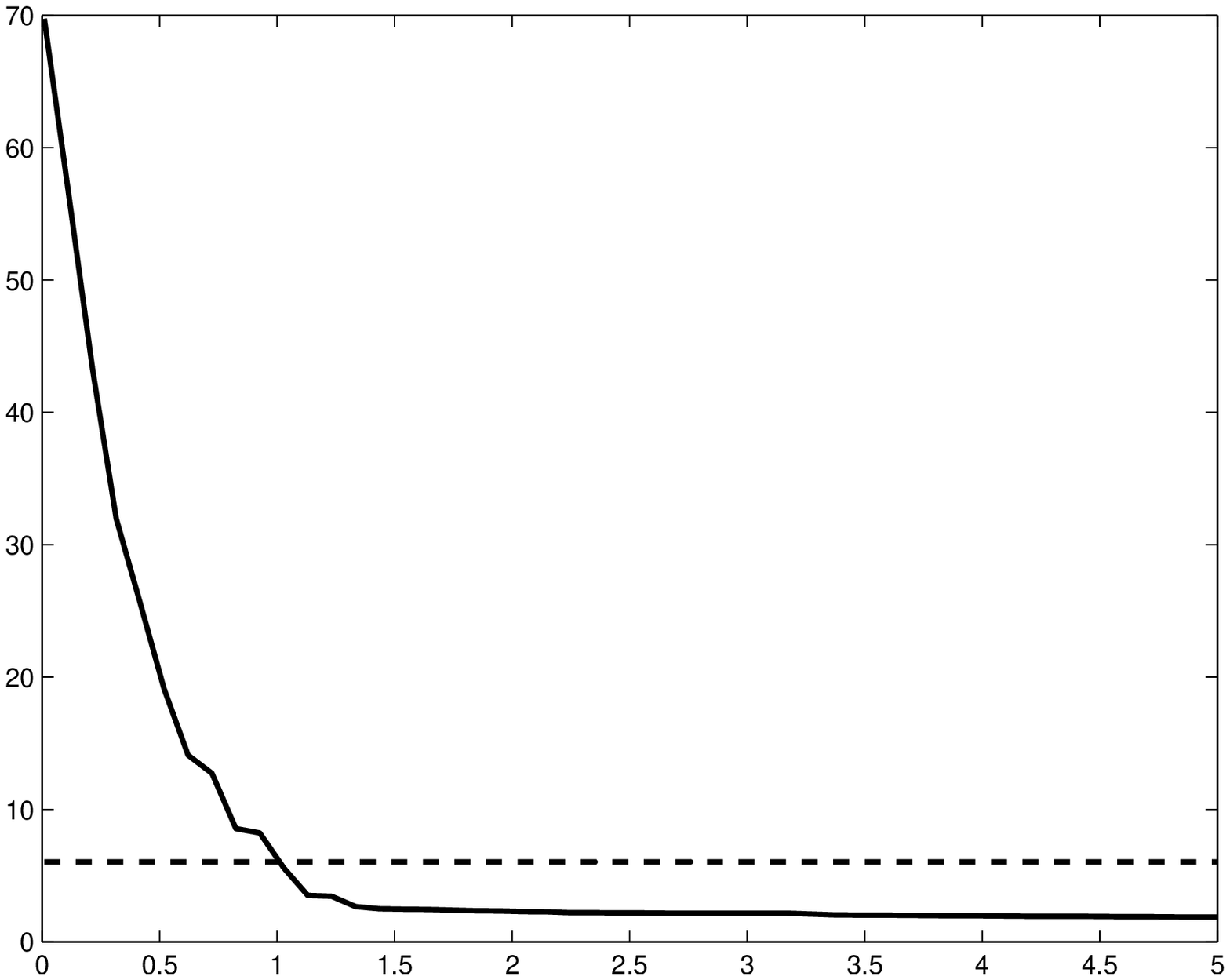} }

\subfigure[MixtGauss - $j_{1} = 3$]
{ \includegraphics[width=3.5cm]{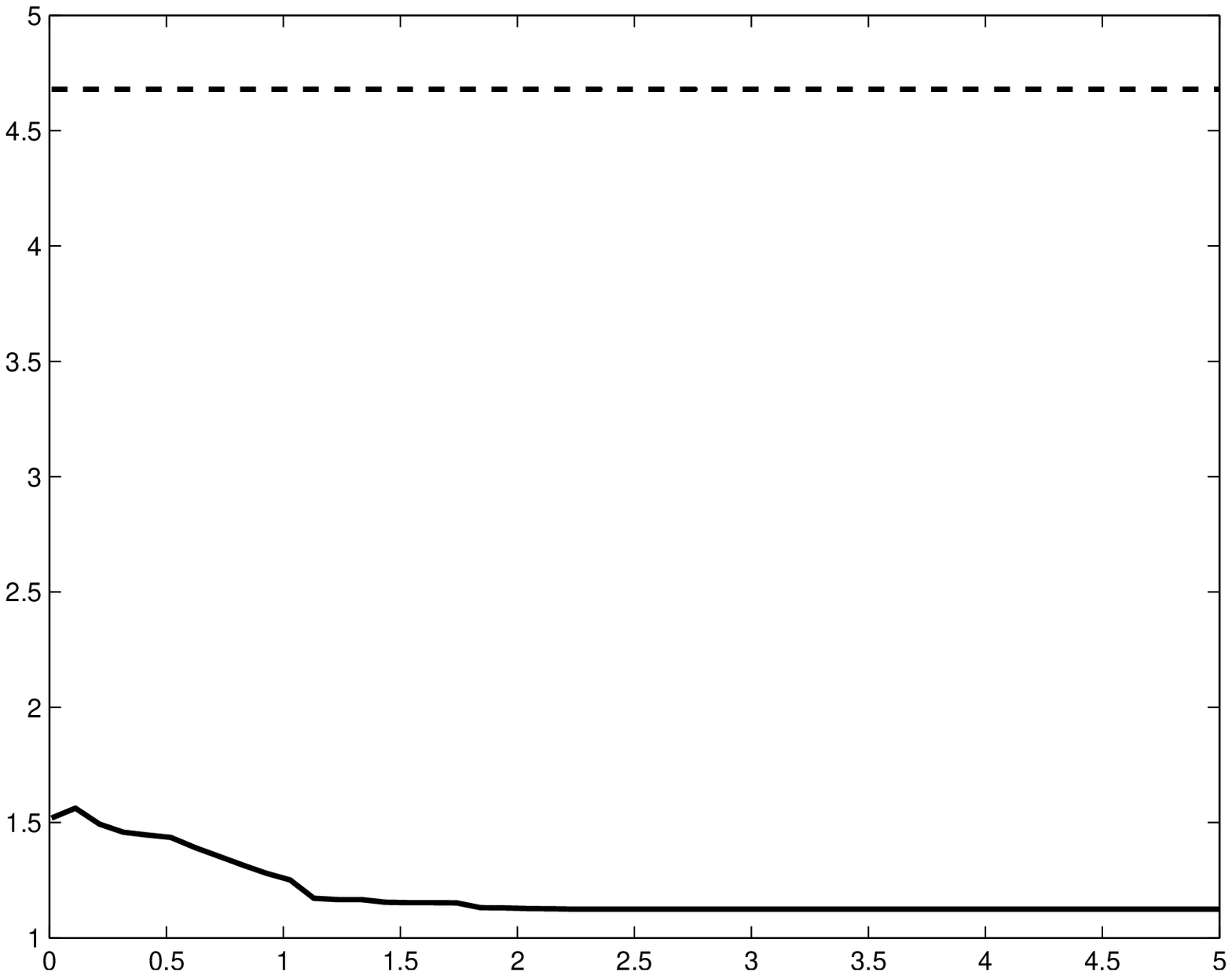} }
\subfigure[MixtGauss - $j_{1} = 4$]
{ \includegraphics[width=3.5cm]{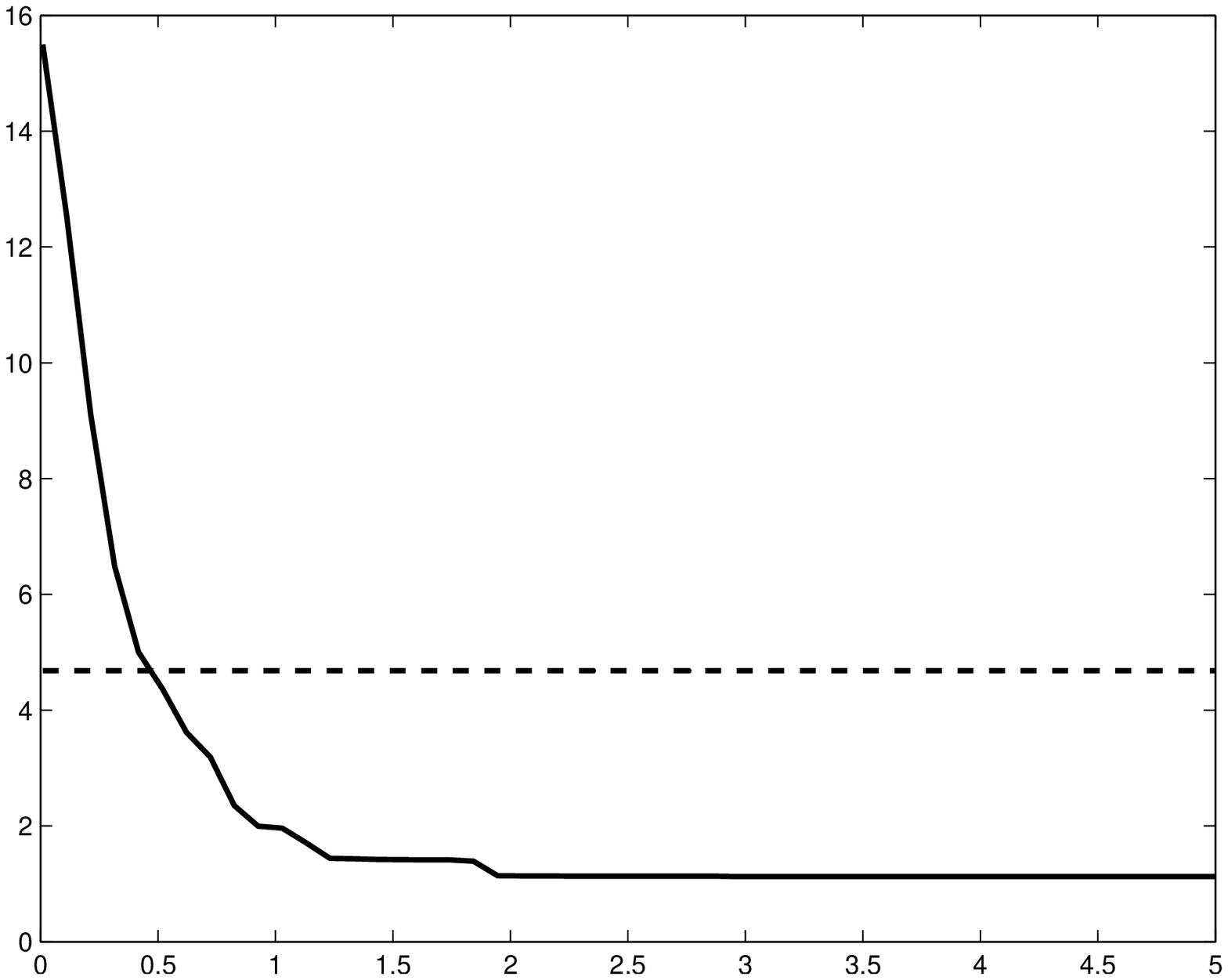} }
\subfigure[MixtGauss - $j_{1} = 5$]
{ \includegraphics[width=3.5cm]{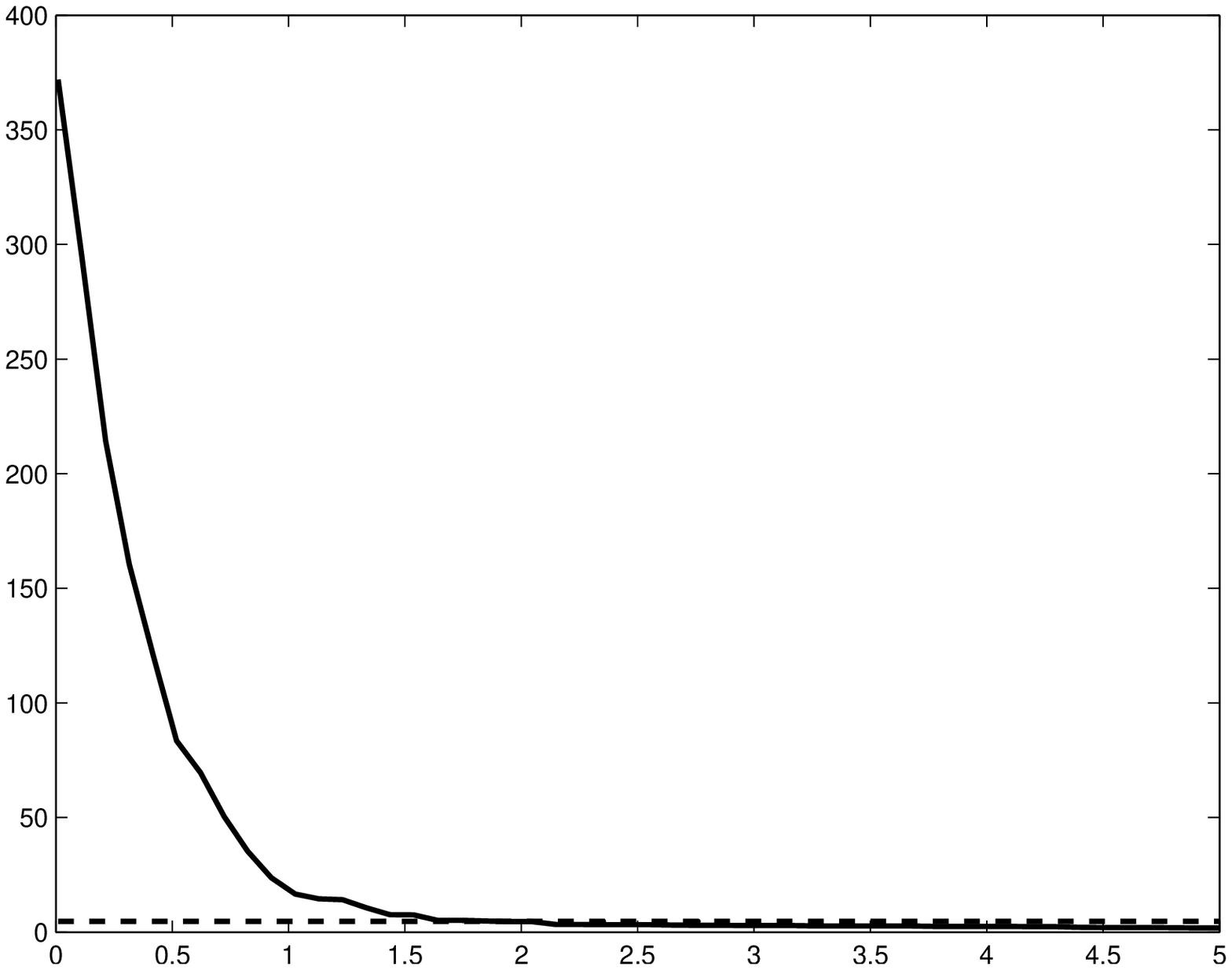} }

\caption{Density deconvolution with $s2n = 3$ and $n = 100$. Evolution of  $R_{n}(j_{1},\delta)$  as a function of $\delta \in [0,5]$  for different values of $j_{1} \geq  j_{0} = 3$ (solid curves). The dotted lines represent the risk of the model selection estimator divided by the risk of the oracle.} \label{fig:riskdecn100s2n3}
\end{figure}

\begin{figure}[htbp]
\centering
\subfigure[Uniform - $j_{1} = 3$]
{ \includegraphics[width=3.5cm]{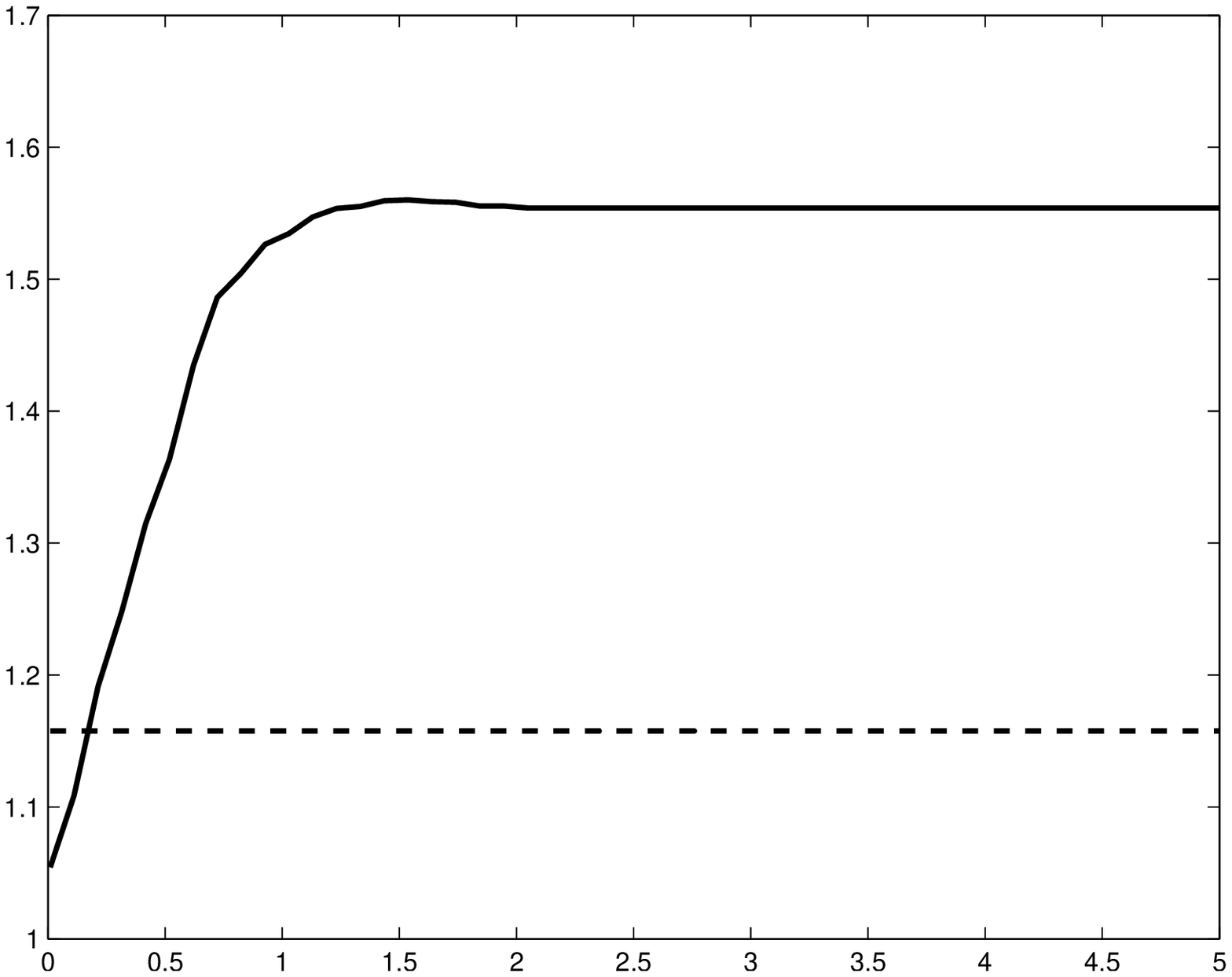} }
\subfigure[Uniform - $j_{1} = 4$]
{ \includegraphics[width=3.5cm]{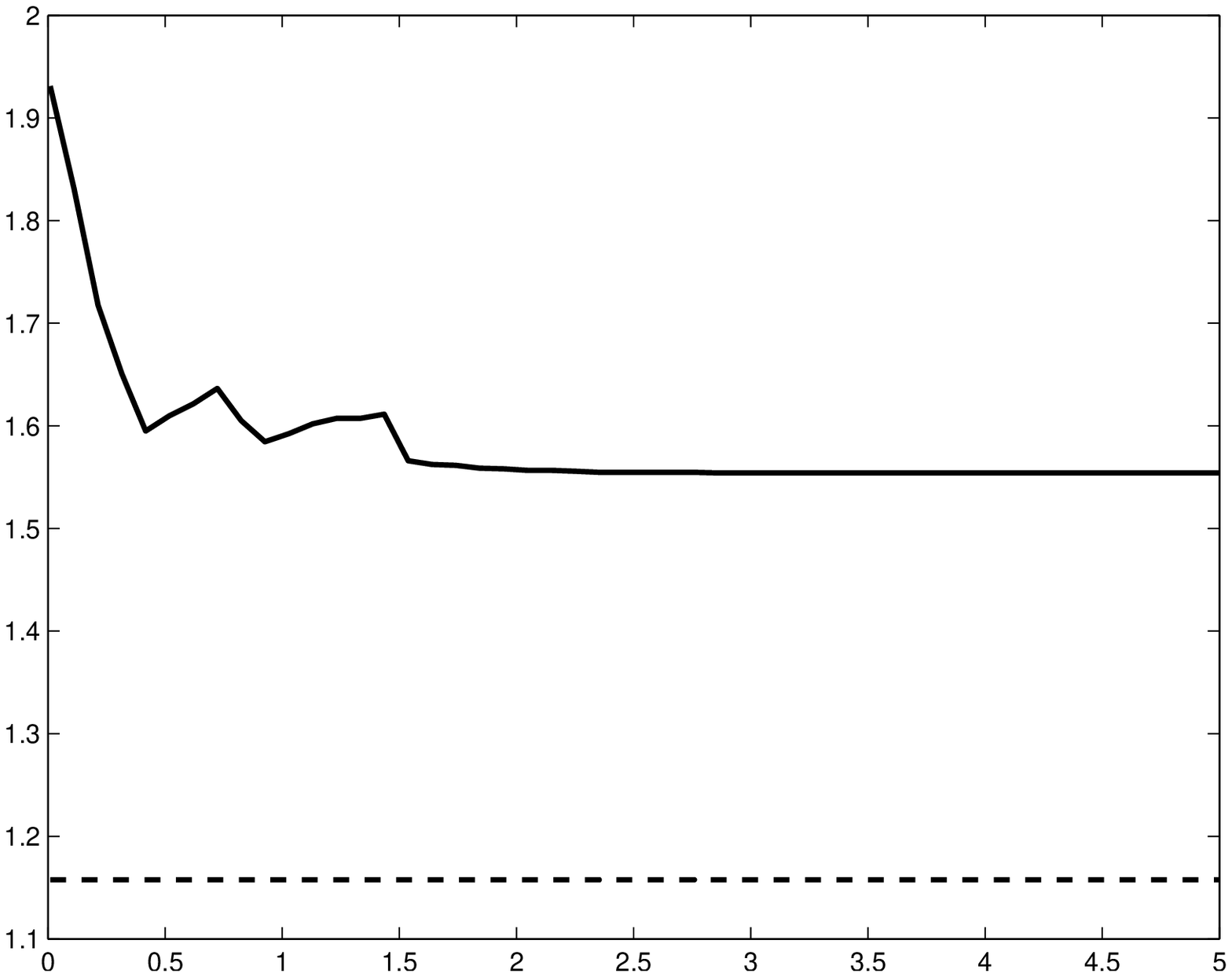} }
\subfigure[Uniform - $j_{1} = 5$]
{ \includegraphics[width=3.5cm]{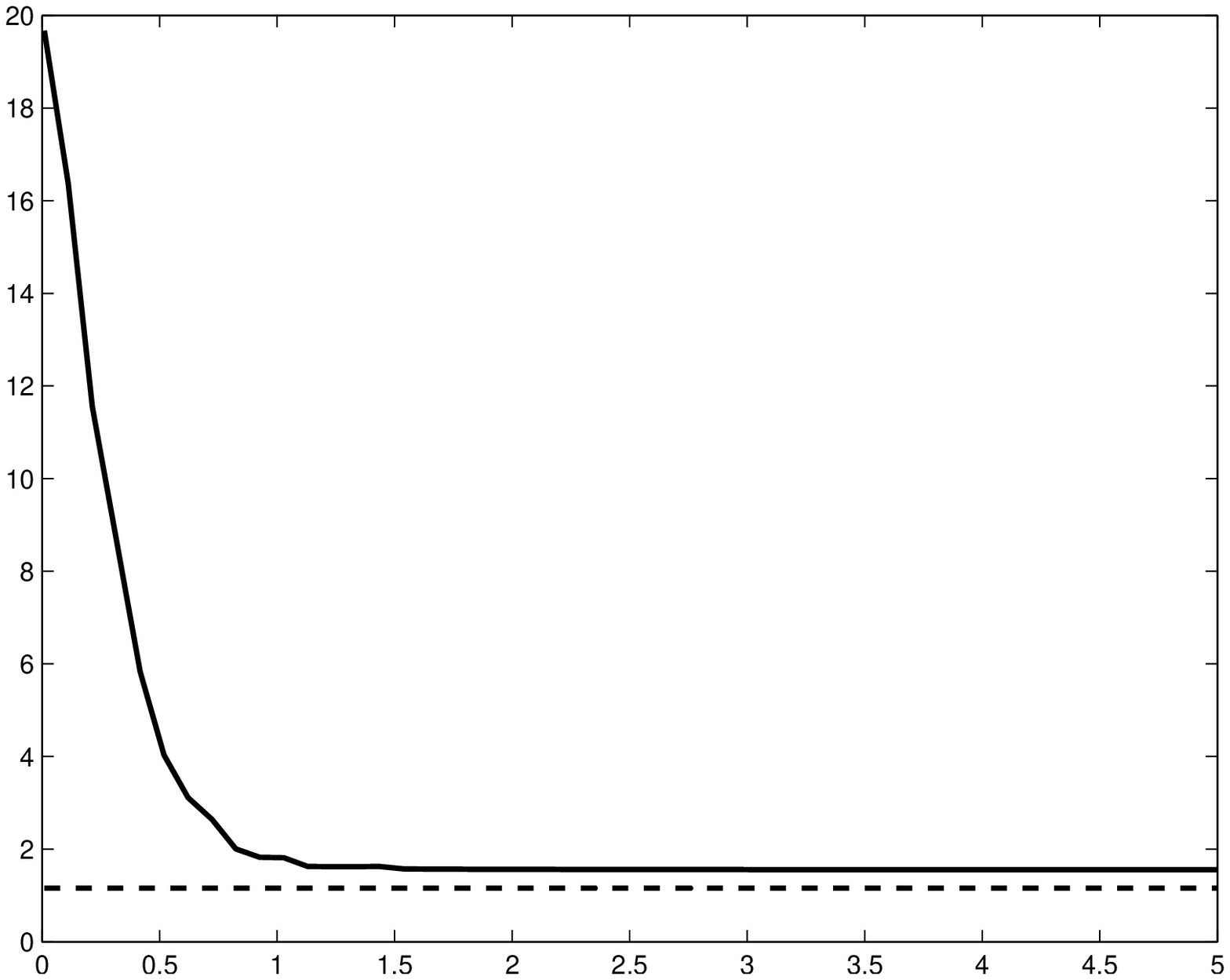} }

\subfigure[Exponential - $j_{1} = 3$]
{ \includegraphics[width=3.5cm]{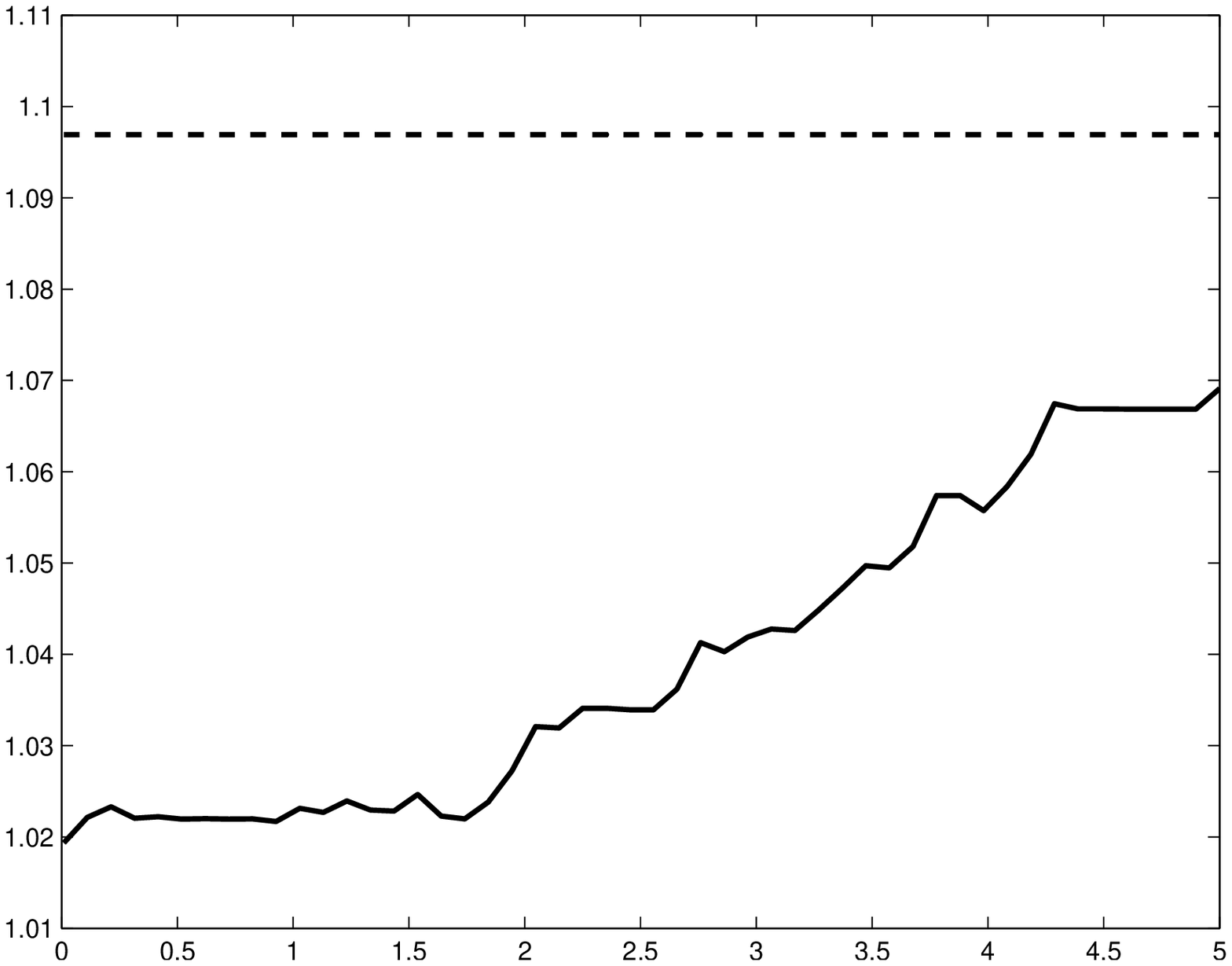} }
\subfigure[Exponential - $j_{1} = 4$]
{ \includegraphics[width=3.5cm]{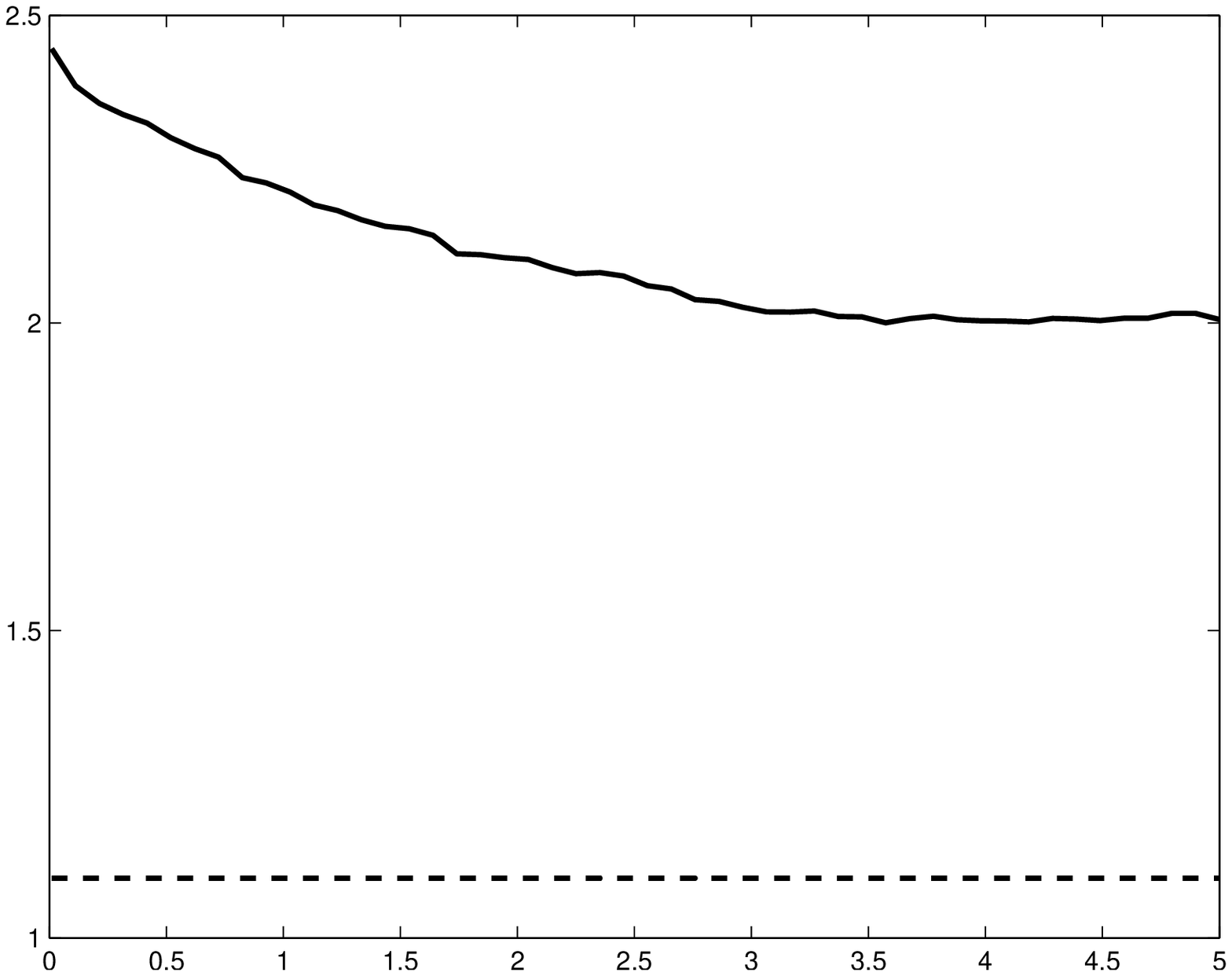} }
\subfigure[Exponential - $j_{1} = 5$]
{ \includegraphics[width=3.5cm]{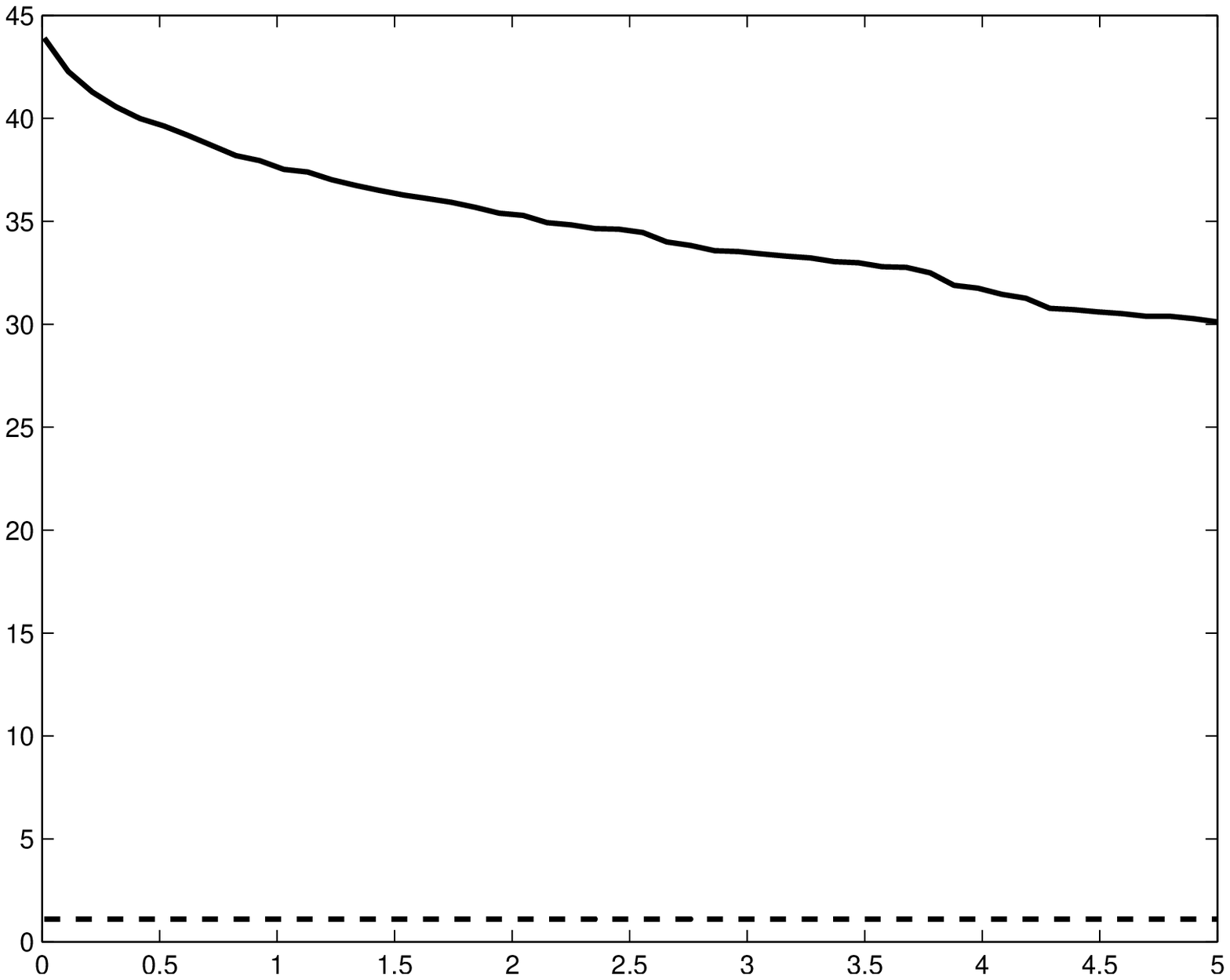} }

\subfigure[Laplace - $j_{1} = 3$]
{ \includegraphics[width=3.5cm]{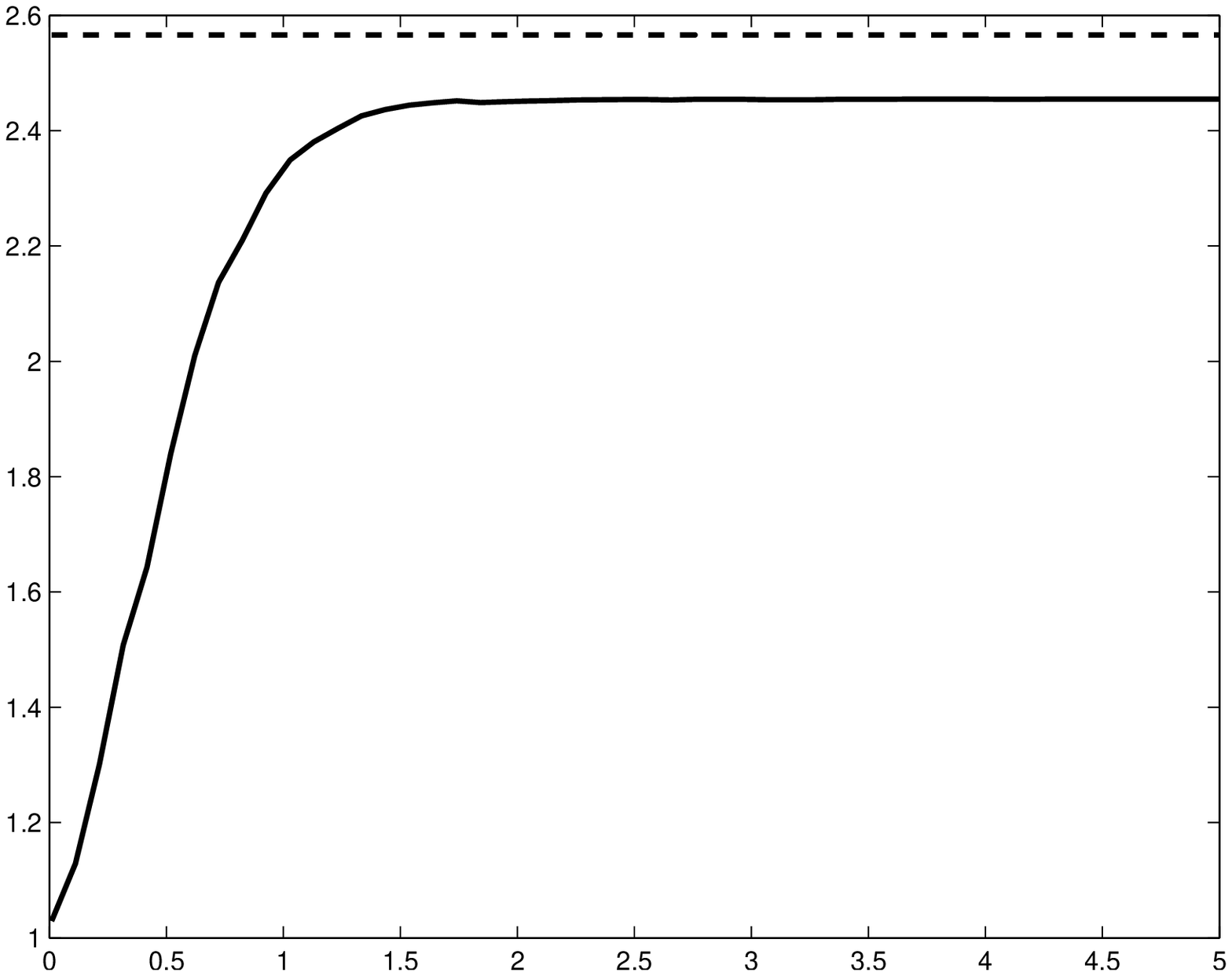} }
\subfigure[Laplace - $j_{1} = 4$]
{ \includegraphics[width=3.5cm]{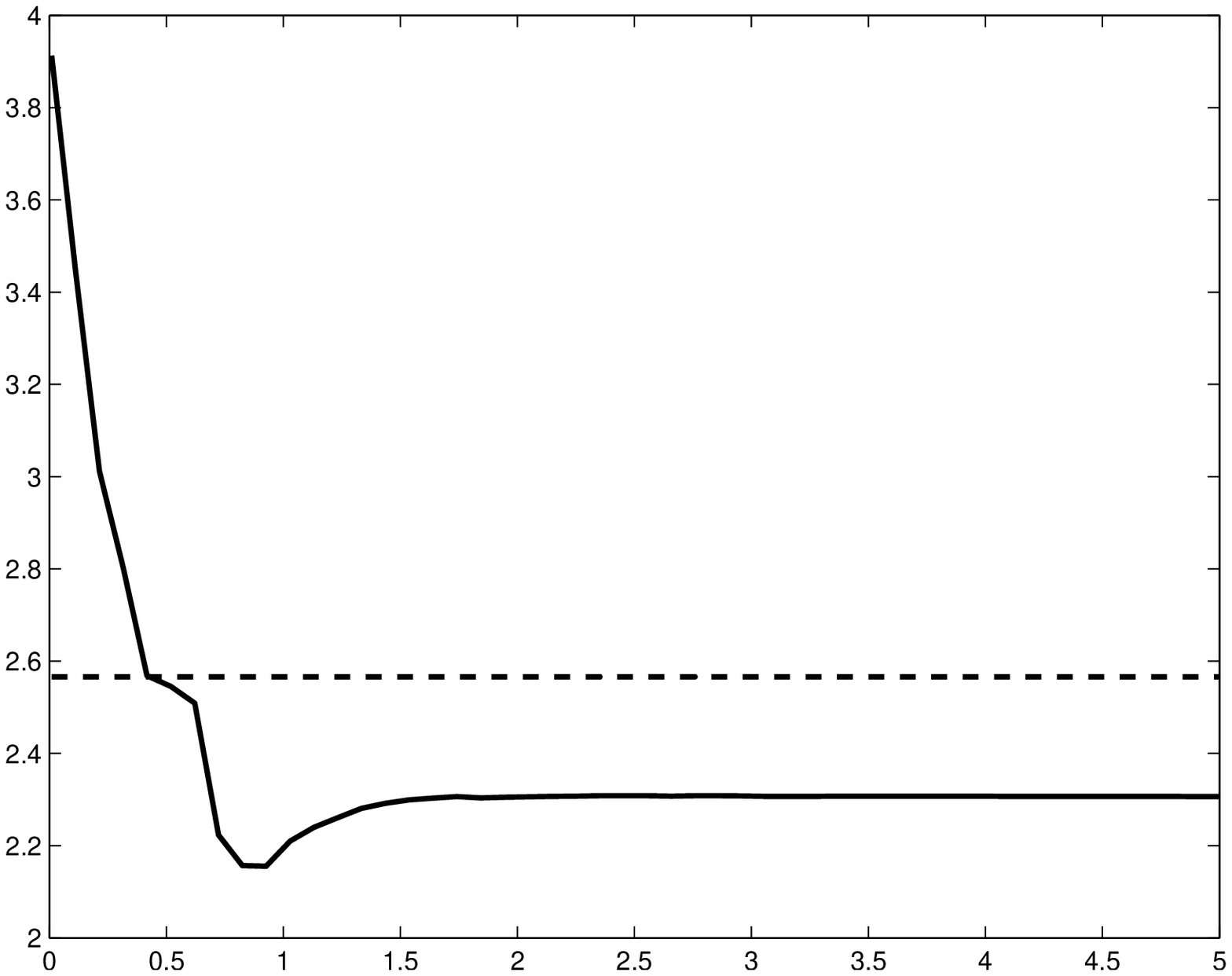} }
\subfigure[Laplace - $j_{1} = 5$]
{ \includegraphics[width=3.5cm]{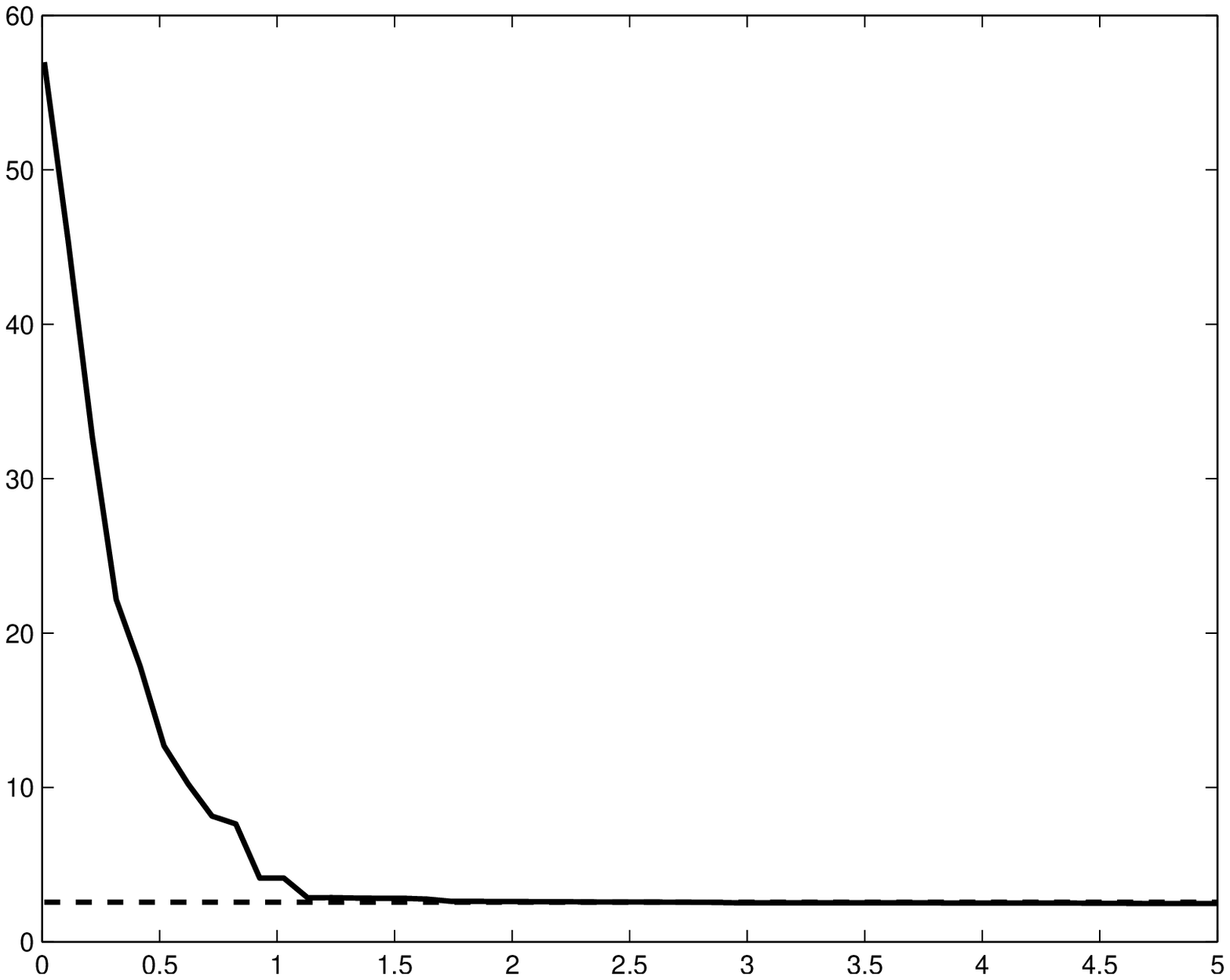} }

\subfigure[MixtGauss - $j_{1} = 3$]
{ \includegraphics[width=3.5cm]{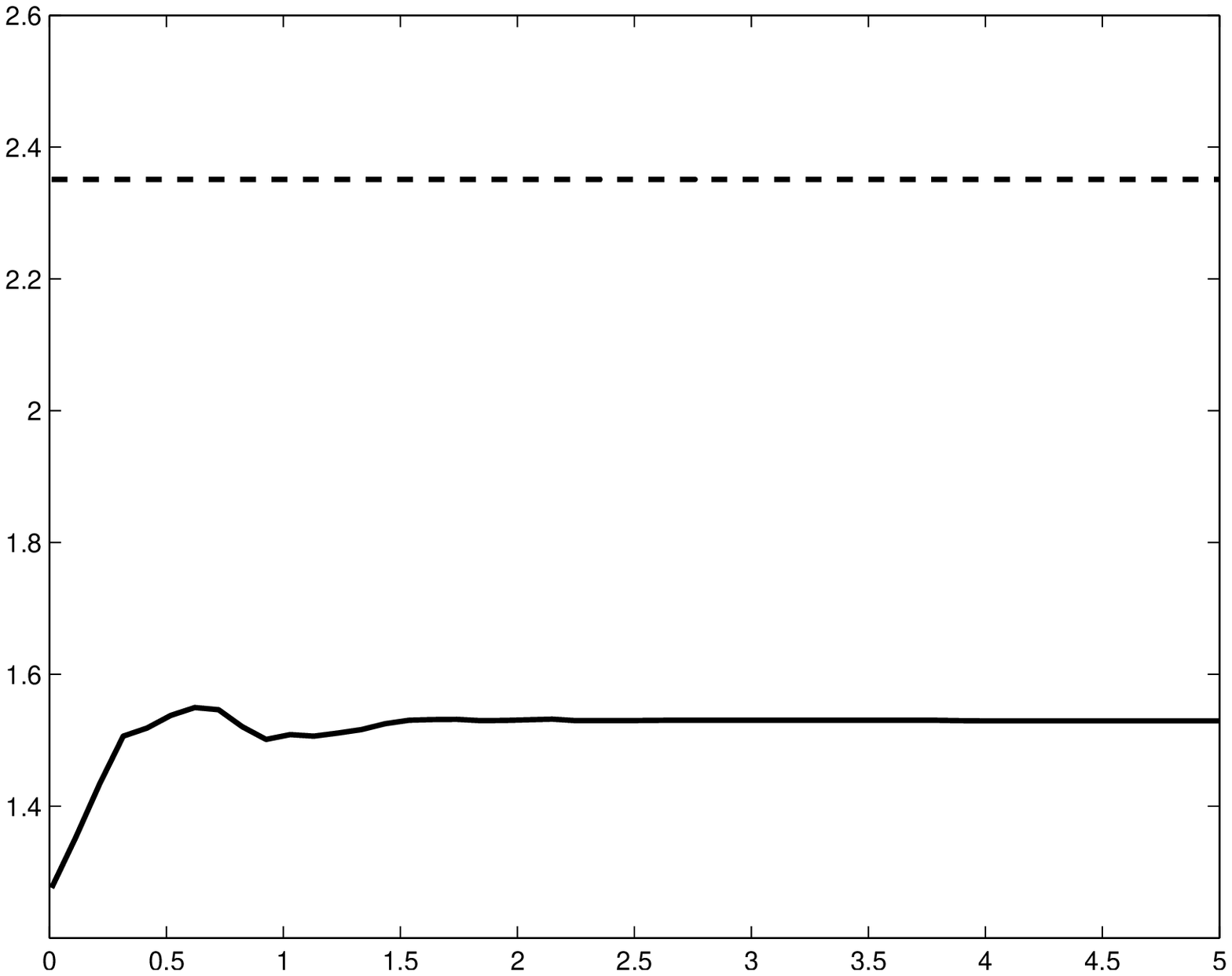} }
\subfigure[MixtGauss - $j_{1} = 4$]
{ \includegraphics[width=3.5cm]{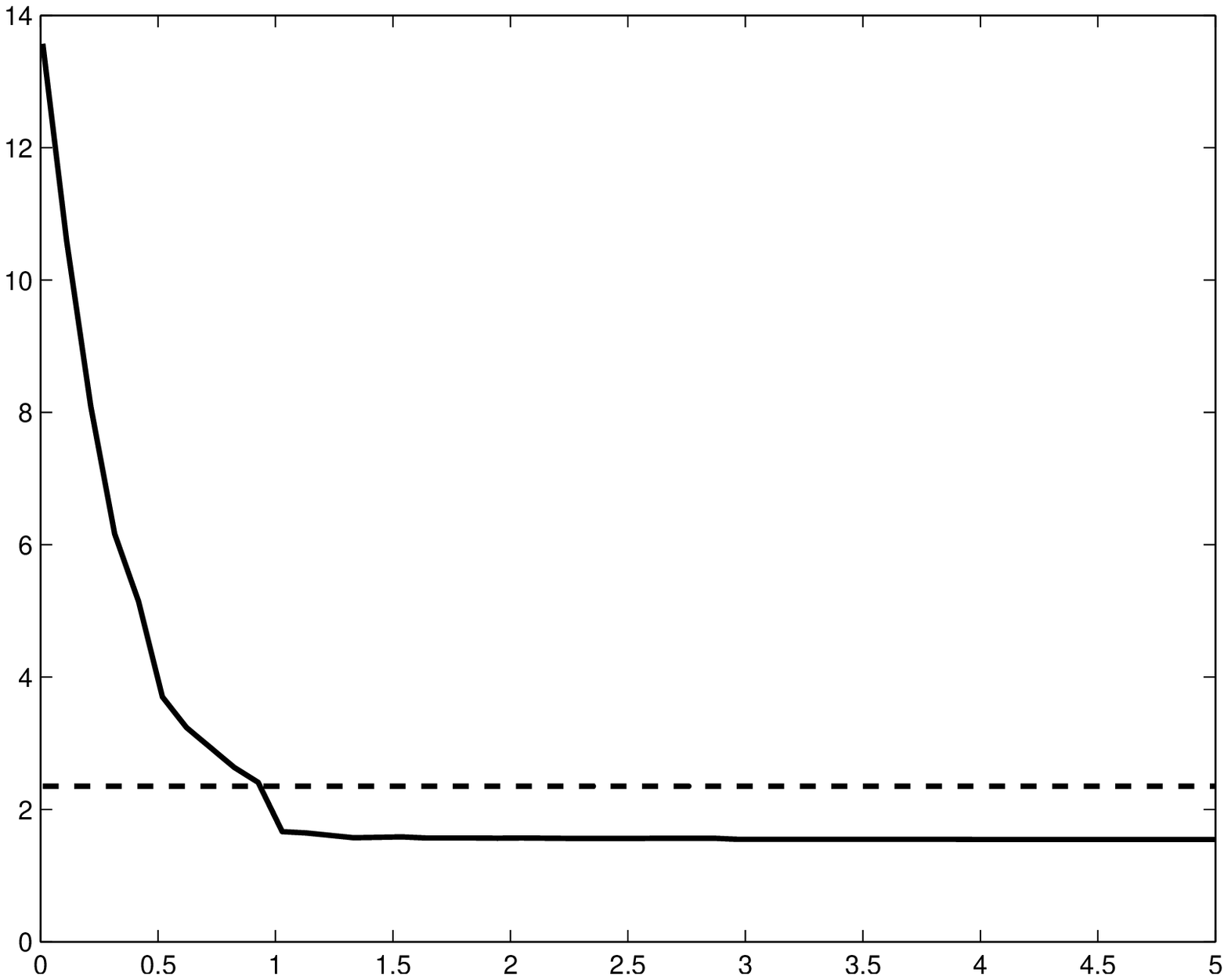} }
\subfigure[MixtGauss - $j_{1} = 5$]
{ \includegraphics[width=3.5cm]{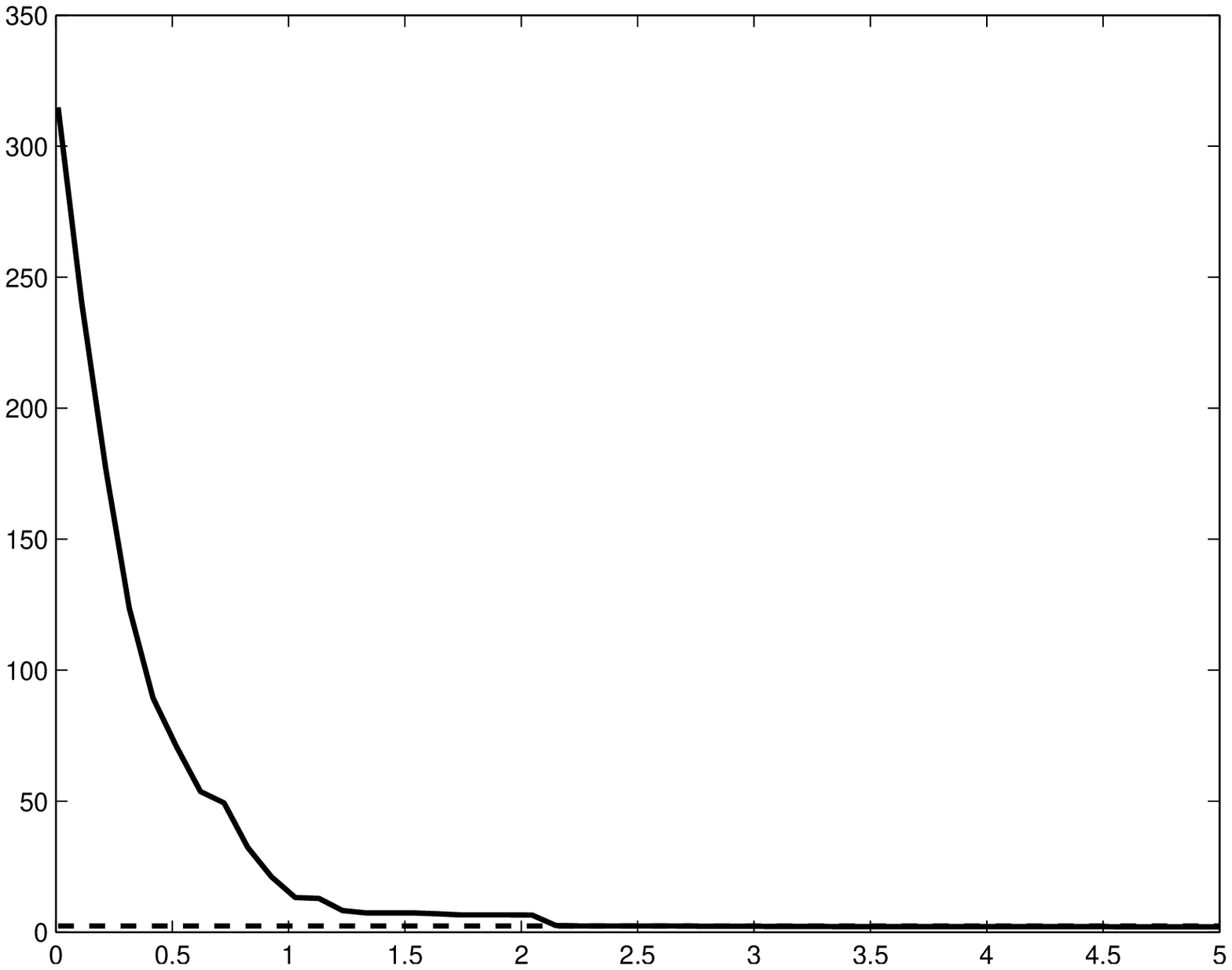} }

\caption{Density deconvolution with $s2n = 3$ and $n = 200$. Evolution of  $R_{n}(j_{1},\delta)$  as a function of $\delta \in [0,5]$  for different values of $j_{1} \geq  j_{0} = 3$ (solid curves). The dotted lines represent the risk of the model selection estimator divided by the risk of the oracle.} \label{fig:riskdecn200s2n3}
\end{figure}

For a given $j_{1}$ and $\delta$, the risk  $R_{n}(j_{1},\delta)$ of the wavelet-based  estimator, see (\ref{eq:riskj1}),  will be compared to the risk of the model section procedure divided by the risk of the oracle $\tilde{f}_{n}$. For each test density, $M=100$ independent samples of size $n = 100,200$ are drawn for $s2n = 3$. Empirical average of $R_{n}(j_{1},\delta)$ and of the risk of the model selection estimator (divided by the oracle risk) over these $M$ repetitions are displayed in Figure \ref{fig:riskdecn100s2n3} ($n=100$) and   Figure \ref{fig:riskdecn200s2n3} ($n=200$) for $\delta \in [0,5]$ and $ j_{0} \leq j_{1} \leq  j_{0} + 2$. For $j_{1} = j_{0}$, wavelet thresholding clearly outperforms the model selection estimator for all values of $\delta$ and all densities, expect for the Uniform distribution with $n=200$ for which it can be seen that model  selection is slightly better if $\delta$ is larger than 0.2. These simulations, also show that the choice $j_{1} = j_{0} = 3$ yields the best results. This observation is consistent with the condition of Theorem \ref{theo:oracledec} that suggests a smaller $j_{1}$ for ill-posed inverse problems than in the direct case. It also confirms that introducing a higher level of resolution clearly deteriorates the quality of the estimator when compared to the oracle risk.

Combining the above remarks, we finally suggest the following choice
$$
j_{1} = j_{0} = \lfloor \log_{2}(\log(n)) \rfloor + 1 \mbox{ and } \delta  = ( 2\nu +1 ) (j_{1}-1-\alpha \log_{2}(\log(n)) )/\log_{2}(n) \mbox{ with } \alpha = 0.5.
$$
For $n=100$, this yields to $j_{1} = 3, \delta \approx 0.6761$ and for $n=200$ to  $j_{1} = 3, \delta \approx 0.5215$. The curves in Figure \ref{fig:riskdecn100s2n3} and Figure \ref{fig:riskdecn200s2n3} indicate that such a choice is reasonable, and leads to satisfactory estimators.

\section*{Appendix}

In what follows $C$ will denote a generic constant whose value may change from line to line. Proofs are given for the case where $h$ satisfies  Assumption \ref{assordi} and   Assumption \ref{ass:h}. The proofs for the case of direct density estimation follow from the same arguments.

\subsection{Technical lemmas}

We start  this technical section by a set of lemmas that will be used in the proof of Theorem \ref{theo:oracledec}. For all the results presented below, it is supposed that $f \in D^{2}([0,1])$ with $\|f\|_{\infty} < + \infty$.

To prove these results we will use the following properties which come directly from the fact that Meyer wavelets are band-limited and that under Assumption \ref{assordi} $ |h_{\ell}| \sim |\ell|^{-\nu} $.
\begin{eqnarray}
|\psi_{j,k}| \leq 2^{-j/2}  \mbox{ and } \# C_{j} \leq 4 \pi 2^{j} \label{eq:boundpsi} \\
|h_{\ell}|^{-2} \leq C 2^{2 j \nu}  \mbox{ for all }  \ell \in C_{j},  \label{eq:boundh}
\end{eqnarray}
where $\# C_{j} $ denotes the cardinality of the set $C_{j}$.

\begin{lemma} \label{lemma:var}
For all integer $j$,
$
\sigma^{2}_{j,k} = \EE (\hat{\beta}_{j,k} - \beta_{j,k} )^{2} \leq V_{j,k} \leq C \|f\|_{\infty} \frac{2^{2j \nu}}{n}
$
and
$
\sigma^{2}_{j_{0},k} = \EE (\hat{c}_{j_{0},k} - c_{j_{0},k} )^{2}  \leq C \|f\|_{\infty} \frac{2^{2j_{0} \nu}}{n}.
$
\end{lemma}
{\bf Proof:} recall that $ V_{j,k}$ is an upper bound for $\sigma^{2}_{j,k}$, and that from equation (\ref{eq:Vdens}), one has that $V_{j,k} =  \frac{1}{n} \int_{\RR} |\tilde{\psi}_{j,k}(y)|^{2} f^{Y}(y)dy$. Then, remark that since $h$ satisfies Assumption \ref{ass:h}, it holds that
\begin{eqnarray*}
V_{j,k}  & =  &  \frac{1}{n} \sum_{m \in \ZZ} \int_{m}^{m+1} |\tilde{\psi}_{j,k}(y)|^{2} f^{Y}(y)dy =  \frac{1}{n} \sum_{m \in \ZZ} \int_{m}^{m+1} |\tilde{\psi}_{j,k}(y)|^{2} \int_{0}^{1} f(u) h(u-y) du dy \\
& \leq &C  \|f\|_{\infty} \frac{1}{n} \sum_{m \in \ZZ} \int_{m}^{m+1} |\tilde{\psi}_{j,k}(y)|^{2} \int_{0}^{1} \frac{1}{|u-y|^{\gamma}} du dy 
\end{eqnarray*}
where $\gamma > 1$ is the real defined in  Assumption \ref{ass:h}. Now, since $\tilde{\psi}_{j,k}(y) $ is a periodic function on $\RR$ with period 1 it follows that
$$
V_{j,k}  \leq C \|f\|_{\infty}  \frac{1}{n} \int_{0}^{1} |\tilde{\psi}_{j,k}(y)|^{2} dy  \sum_{m \in \ZZ} |m|^{-\gamma}
$$
Finally, using    Parseval relation and the bounds   (\ref{eq:boundpsi}) and (\ref{eq:boundh}), it follows that $ \int_{0}^{1} |\tilde{\psi}_{j,k}(y)|^{2} dy = \sum_{\ell \in C_{j}}   |\tilde{\psi}^{j,k}_{\ell}|^{2} = \sum_{\ell \in C_{j}}   \frac{|\psi^{j,k}_{\ell}|^{2} }{|h_{\ell}|^{2}}  \leq C  2^{2 j \nu}$ which finally implies (using the fact that $\gamma > 1$, )
$$
V_{j,k}  \leq C \|f\|_{\infty} \frac{2^{2 j \nu}}{n}
$$
which completes the proof. The argument is the same to bound $\sigma^{2}_{j_{0},k}$.  \hfill $\Box$

\begin{lemma} \label{lemma:moment}
For all $p \geq 2$,
$
\EE (\hat{\beta}_{j,k} - \beta_{j,k} )^{2p} \leq C \max(\|f\|_{\infty}^{p},1) \left( \frac{ 2^{ 2 j \nu p}}{n^{p}} +  \frac{2^{jp(2\nu +1)} }{n^{2p-1}}  \right). 
$
\end{lemma}
{\bf Proof:}  by definition $\hat{\beta}_{j,k} - \beta_{j,k} = \frac{1}{n} \sum_{m=1}^{n} Z_{m}$ where $Z_{m} =  \sum_{\ell \in C_{j}}   \tilde{\psi}^{j,k}_{\ell} \left( e^{-2\pi i \ell Y_{m}} - f_{\ell}h_{\ell} \right)$. Remark that for all $m$, $\EE Z_{m}=0$, and that from Lemma \ref{lemma:var} $\var(Z_{m}) \leq C \|f\|_{\infty}  2^{2 j  \nu}$. Since $f$ and $h$ are densities, $|f_{\ell} h_{\ell}| \leq 1$. Hence using  (\ref{eq:boundpsi}) and (\ref{eq:boundh}), this implies that $|Z_{m}| \leq  2   \sum_{\ell \in C_{j}}  \frac{|\psi^{j,k}_{\ell}|}{|h_{\ell}|} \leq C 2^{j(\nu +1/2)}$. Then the result follows from Rosenthal's inequality (see \citeA{MR0440354}). \hfill $\Box$

\begin{lemma} \label{lemma:ineg1}
For any positive $x$,
$
\PP \left( |\hat{\beta}_{j,k} - \beta_{j,k} | \geq \sqrt{2 V_{jk} x} + \frac{ \eta_{j} }{3n} x \right) \leq 2 \exp(-x)
$, 
where $\eta_{j} = \sum_{\ell  \in C_{j}} |\tilde{\psi}_{jk}^{\ell}|$.
\end{lemma}
{\bf Proof:} note that  $\hat{\beta}_{j,k} - \beta_{j,k} = \frac{1}{n} \sum_{m=1}^{n} (W_{m}- \EE W_{m})$ where $W_{m} =  \sum_{\ell  \in C_{j}}   \tilde{\psi}^{j,k}_{\ell} e^{-2\pi i \ell Y_{m}}$. Then, remark that by definition $V_{jk} = \frac{1}{n^{2}} \sum_{m=1}^{n}  \EE |W_{m}|^{2}$, and that $|W_{m}| \leq \sum_{\ell  \in C_{j}}  |\tilde{\psi}_{jk}^{\ell}|$. Hence, the $Z_{m}$'s are bounded random variables,  and thus the result follows from Bernstein's inequality (see e.g. Proposition 2.9 in \citeA{massart}). \hfill $\Box$

\begin{lemma} \label{lemma:ineg2}
For any positive $x$,
$$
\PP  \left( V_{jk}  \geq  \hat{V}_{jk} + \sqrt{2  \frac{\eta_{j}^{2}}{n^{2}}   \hat{V}_{jk} x}   +  \kappa  \frac{\eta_{j}^{2}}{n^{2}} x \right) \leq \exp \left(-x \right).
$$
where $\kappa = \frac{4}{3} +  \sqrt{\frac{5}{3}}$ and $\eta_{j} = \sum_{\ell  \in C_{j}} |\tilde{\psi}_{jk}^{\ell}|$.
\end{lemma}
{\bf Proof:} the proof is inspired by the proof of Lemma 1 in \citeA{reyriv}. By definition $\hat{V}_{jk} = \frac{1}{n^{2}} \sum_{m=1}^{n} W_{m}$ with $W_{m} =  \sum_{\ell, \ell' \in C_{j}}  \tilde{\psi}^{j,k}_{\ell} \overline{ \tilde{\psi}^{j,k}_{\ell'} }  e^{-2\pi i (\ell-\ell') Y_{m}}$. Then, remark that $|W_{m}| \leq  \sum_{\ell, \ell' \in C_{j}} | \tilde{\psi}^{j,k}_{\ell} | |\tilde{\psi}^{j,k}_{\ell'}|$ which implies that $|W_{m}| \leq \eta_{j}^{2}$ for all $m=1,\ldots,n$. Moreover, one can remark that $\EE |W_{m}|^{2} = \EE |\tilde{\psi}_{j,k}(Y_{m})|^{4} =  \int_{\RR} |\tilde{\psi}_{j,k}(y)|^{4} f^{Y}(y)dy$. Then, since $|\tilde{\psi}_{j,k}(y)|^{2} \leq   \left(\sum_{\ell \in C_{j}} | \tilde{\psi}^{j,k}_{\ell} | \right)^{2}$, it follows that
$$
\EE |W_{m}|^{2} \leq   \eta_{j}^{2} \int_{\RR} |\tilde{\psi}_{j,k}(y)|^{2} f^{Y}(y)dy   \leq \eta_{j}^{2} n V_{jk}.
$$
Hence, by applying Bernstein's inequality (see e.g. Proposition 2.9 in \citeA{massart}), one obtains that for any positive $x$
\begin{equation}
\PP  \left( V_{jk}  \geq \hat{V}_{jk} +  \sqrt{\frac{2 \eta_{j}^{2}}{n^{2}} V_{jk} x } + \frac{\eta_{j}^{2}}{3 n^{2}} x \right) \leq \exp \left(-x \right). \label{eq:bern}
\end{equation}
Now, let $u = \frac{\eta_{j}^{2}}{n^{2}} x$ and let $g(y) = y^{2} - \sqrt{2 u} y - \frac{u}{3}$ for $y \geq 0$. From (\ref{eq:bern}), one has that
$$
\PP  \left(  g(\sqrt{V_{jk}}) \geq  \hat{V}_{jk}  \right) \leq \exp \left(-x \right).
$$
 As $V_{jk} $ and $\hat{V}_{jk} $ are positive, one can check that it is possible to invert the inequality  $g(\sqrt{V_{jk}}) \geq V_{jk}$ to obtain that
$$
\PP  \left(  \sqrt{V_{jk}} \geq g^{-1}( \hat{V}_{jk} ) \right) \leq \exp \left(-x \right),
$$
and that $g^{-1}(y) = \sqrt{y + \frac{5u}{6}} + \sqrt{\frac{u}{2}}$. Hence, one obtains that for any positive $x$
$$
\PP  \left(  V_{jk} \geq \hat{V}_{jk} + \frac{4}{3}u + \sqrt{2u  \hat{V}_{jk}  + \frac{5 u^{2}}{3} } \right) \leq \exp \left(-x \right).
$$
Now, using the fact that $ \sqrt{2u  \hat{V}_{jk}  + \frac{5 u^{2}}{3} } \leq  \sqrt{2u  \hat{V}_{jk} } + \sqrt{\frac{5 u^{2}}{3}}$, it follows that
$$
\hat{V}_{jk} + \frac{4}{3}u + \sqrt{2u  \hat{V}_{jk}  + \frac{5 u^{2}}{3} } \leq \hat{V}_{jk} + \sqrt{2u  \hat{V}_{jk}}   +  \kappa u  ,
$$
which completes the proof. \hfill $\Box$\\

\subsection{Proof of Theorem \ref{theo:oracledec}}

Let us define the following set of integers $\Lambda_{n} = \left\{ (j,k); j_{0} \leq j \leq j_{1}, 0 \leq k \leq 2^{j}-1\right\}$. We first establish the following proposition:
\begin{prop} \label{prop:oracle}
Let $\epsilon > 0$. Then, for any subset of indices $m \subset \Lambda_{n}$
\begin{eqnarray*}
\|\hat{f}_{n} - f \|^{2} & \leq & \sum_{k=0}^{2^{j_{0}}-1} (\hat{c}_{j_{0},k}-c_{j_{0},k})^{2} + \left(\frac{2+\epsilon}{\epsilon}\right)^{2} \sum_{(j,k) \in \Lambda_{n} \backslash m} \beta_{jk}^{2} + (2+\epsilon)  \sum_{(j,k) \in   m} \left(\hat{\beta}_{jk} - \beta_{jk} \right)^{2} \\
& & +  \left(\frac{2+\epsilon}{\epsilon}\right)  \sum_{(j,k) \in   m} \tau_{jk}^{2} + \left(\frac{2+\epsilon}{\epsilon}\right)^{2} \sum_{j = j_{1}+1}^{+ \infty }\sum_{k = 0}^{2^{j}-1} \beta_{jk}^{2} \\
& & +  \left(\frac{2+\epsilon}{\epsilon}\right)  \sum_{(j,k) \in   \hat{m}} \left( (1+\epsilon) \left(\hat{\beta}_{jk} - \beta_{jk} \right)^{2} -  \tau_{jk}^{2} \right)
\end{eqnarray*}
where $\hat{m} = \{(j,k); | \beta_{j,k} |  \geqslant \tau_{j,k} \; , j_{0} \leq j \leq j_{1}, 0 \leq k \leq 2^{j}-1 \} $. 
\end{prop}
{\bf Proof:} the proof  is inspired by model selection techniques (see  \citeA{massart} and also  \citeA{reyriv} ). Let $m \subset \Lambda_{n}$ and define
$$
\hat{f}_{m} = \sum_{(j,k) \in m} \hat{\beta}_{jk} \psi_{jk} \mbox{ and } f_{m} = \sum_{(j,k) \in m}\beta_{jk} \psi_{jk}.
$$
Note that  $\hat{f}_{\hat{m}}$ is given by
$
\hat{f}_{\hat{m}} = \sum_{j = j_{0}}^{j_{1}}\sum_{k = 0}^{2^{j}-1} \hat{\beta}_{jk} \1_{\{| \beta_{j,k} |  \geqslant \tau_{j,k} \}} \psi_{j,k}.
$
Then, define
$$
\gamma_{m} = - \sum_{(j,k) \in m}  \hat{\beta}_{jk}^{2} \mbox{ and } \pen(m) =  \sum_{(j,k) \in m} \tau_{jk}^{2},
$$
and one can check that
$$
\hat{m} = \arg \min_{m \in \Lambda_{n}} \gamma_{m} + \pen(m).
$$
Now let $\tilde{f} =  \sum_{j = j_{0}}^{+ \infty }\sum_{k = 0}^{2^{j}-1} \beta_{jk} \psi_{jk}$ and remark that for any $m \subset \Lambda_{n}$
\begin{equation}
\gamma_{m} = \|\hat{f}_{m}-\tilde{f} \|^{2} -  \sum_{j = j_{0}}^{+ \infty }\sum_{k = 0}^{2^{j}-1} \beta_{jk}^{2} - 2 \sum_{(j,k) \in m } \hat{\beta}_{jk}(\hat{\beta}_{jk}-\beta_{jk}). \label{eq:gamma}
\end{equation}
Then, by definition of $\hat{m}$, equality (\ref{eq:gamma}) implies that for all $m \subset \Lambda_{n}$
$$
 \|\hat{f}_{\hat{m}}-\tilde{f} \|^{2} \leq  \|\hat{f}_{m}-\tilde{f} \|^{2} - 2 \sum_{(j,k) \in m } \hat{\beta}_{jk}(\hat{\beta}_{jk}-\beta_{jk}) + 2 \sum_{(j,k) \in \hat{m} } \hat{\beta}_{jk}(\hat{\beta}_{jk}-\beta_{jk}) + \pen(m) - \pen(\hat{m}).
$$
Using the Pythagorean equality $ \|\hat{f}_{m}-\tilde{f} \|^{2} =  \|f_{m}-\tilde{f} \|^{2} +  \|\hat{f}_{m}-f_{m} \|^{2}  $ one finally obtains that
\begin{eqnarray}
 \|\hat{f}_{\hat{m}}-\tilde{f} \|^{2} & \leq&   \|f_{m}-\tilde{f} \|^{2} -   \|\hat{f}_{m}-f_{m} \|^{2} +  2 \sum_{(j,k) \in \hat{m} } \hat{\beta}_{jk}(\hat{\beta}_{jk}  -\beta_{jk}) - 2 \sum_{(j,k) \in m } \beta_{jk}(\hat{\beta}_{jk}-\beta_{jk}) \nonumber \\ 
 & & + \pen(m) - \pen(\hat{m}). \label{eq:oracle1}
\end{eqnarray}
Let $\epsilon > 0$ and remark that,
\begin{equation}
2 \sum_{(j,k) \in \hat{m} } \hat{\beta}_{jk}(\hat{\beta}_{jk}  -\beta_{jk}) - 2 \sum_{(j,k) \in m } \beta_{jk}(\hat{\beta}_{jk}-\beta_{jk}) \leq 2 \left( \|\hat{f}_{\hat{m}}-\tilde{f} \| + \|f_{m}-\tilde{f} \|  \right)\sqrt{\sum_{(j,k) \in \hat{m} \cup m } (\hat{\beta}_{jk}  -\beta_{jk})^2  } \label{eq:beta}
\end{equation}
Then, remark that $\sum_{(j,k) \in \hat{m} \cup m } (\hat{\beta}_{jk}  -\beta_{jk})^2 \leq  \|\hat{f}_{m}-f_{m} \|^{2}  +  \|\hat{f}_{\hat{m}}-f_{\hat{m}} \|^{2} $ and thus by using twice the inequality $2ab \leq \theta a^{2} + (1/\theta)b^{2}$ with $\theta = 2/(2+\epsilon)$ and $\theta = 2/\epsilon$, and by inserting inequality (\ref{eq:beta}) in (\ref{eq:oracle1}), one obtains that
\begin{eqnarray}
\frac{\epsilon}{2+\epsilon} \|\hat{f}_{\hat{m}}-\tilde{f} \|^{2} & \leq&  \frac{2+\epsilon}{\epsilon} \|f_{m}-\tilde{f} \|^{2} + \epsilon   \|\hat{f}_{m}-f_{m} \|^{2} + \pen(m)  \nonumber \\ 
 & & + (1+\epsilon)\|\hat{f}_{\hat{m}}-f_{\hat{m}} \|^{2}  - \pen(\hat{m}). \label{eq:oracle2}
\end{eqnarray}
which completes the proof since $\hat{f}_{n} =  \sum_{k = 0}^{2^{j_{0}}-1} \hat{c}_{j_{0}k}  \phi_{j_{0},k} + \hat{f}_{\hat{m}}$. \hfill $\Box$\\

Next we prove the following lemma:

\begin{lemma} \label{control2}
For any $\delta > \eta \left(1 + \frac{\nu}{\nu + 1} \right)$, there exists $\epsilon = \epsilon(\delta)  > 0$ such that $\delta > (1 + \epsilon) \left(1 + \frac{\nu}{\nu + 1} \right) $, and a positive constant $C$ not depending on $n$ such that 
$$
\EE \left(  \sum_{(j,k) \in   \hat{m}} (1+\epsilon) \left(\hat{\beta}_{jk} - \beta_{jk} \right)^{2} -  \tau_{jk}^{2} \right) \leq C \max(\|f\|_{\infty},1) \max(\log(n)^{\alpha},1) \frac{ (\log n)^{\alpha (2 \nu +1)}}{n}.
$$
\end{lemma}
{\bf Proof:} let $Z =  \sum_{(j,k) \in   \hat{m}} (1+\epsilon) \left(\hat{\beta}_{jk} - \beta_{jk} \right)^{2} -  \tau_{jk}^{2}$. By definition of $\hat{m}$
\begin{eqnarray*}
\EE Z & = & \EE  \sum_{j = j_{0}}^{j_{1}} \sum_{k=0}^{2^{j}-1}  \left( (1+\epsilon) \left(\hat{\beta}_{jk} - \beta_{jk} \right)^{2} -  \tau_{jk}^{2} \right) \1_{\{ |\hat{\beta}_{jk}| \geq \tau_{jk} \}} \\
& \leq & \EE  \sum_{j = j_{0}}^{j_{1}} \sum_{k=0}^{2^{j}-1}  \left( (1+\epsilon) \left(\hat{\beta}_{jk} - \beta_{jk} \right)^{2} \right) \1_{\{ |\hat{\beta}_{jk}| \geq \tau_{jk} \cap |\hat{\beta}_{jk} - \beta_{jk} | \geq \frac{\tau_{jk}}{\sqrt{1+\epsilon}}  \} } \\
& \leq &   (1+\epsilon) \sum_{j = j_{0}}^{j_{1}} \sum_{k=0}^{2^{j}-1}   \EE \left(   \left(\hat{\beta}_{jk} - \beta_{jk} \right)^{2}  \1_{ \{ |\hat{\beta}_{jk} - \beta_{jk} | \geq \frac{\tau_{jk}}{\sqrt{1+\epsilon}} \} }  \right)
\end{eqnarray*}
Now by applying Holder inequality, one has that for all $p \geq 2$ and $q > 1$ such that $\frac{1}{p} + \frac{1}{q} = 1$
\begin{equation}
\EE Z \leq  (1+\epsilon) \sum_{j = j_{0}}^{j_{1}}   \sum_{k=0}^{2^{j}-1}    \left(  \EE   \left(\hat{\beta}_{jk} - \beta_{jk} \right)^{2p} \right)^{1/p} \left( \PP \left( |\hat{\beta}_{jk} - \beta_{jk} | \geq \frac{\tau_{jk}}{\sqrt{1+\epsilon}}  \right) \right)^{1/q} \label{eq:holder}
\end{equation}
Now by Lemma \ref{lemma:moment} one has that
\begin{equation}
 \left(  \EE   \left(\hat{\beta}_{jk} - \beta_{jk} \right)^{2p} \right)^{1/p} \leq C \max(\|f\|_{\infty},1) \left( \frac{ 2^{ 2 j \nu p}}{n^{p}} +  \frac{2^{jp(2\nu +1)} }{n^{2p-1}}  \right)^{1/p}. \label{eq:c1}
\end{equation}
By definition of $j_{1}$ and $j_{0}$ and since $\eta \leq 1/2$, one has that there exists a constant $C$ such that for all $j_{0} \leq j \leq j_{1}$, $\frac{2^{j}}{\sqrt{n}} \leq C \log(n)^{\alpha}$ which implies (since $p \geq 2$) that $ \frac{2^{jp(2\nu +1)} }{n^{2p-1}} \leq C \frac{2^{2j\nu p} }{n^{p}} \log(n)^{\alpha p}$. By inserting this inequality in (\ref{eq:c1}) one obtains that
\begin{equation}
 \left(  \EE   \left(\hat{\beta}_{jk} - \beta_{jk} \right)^{2p} \right)^{1/p} \leq C \max(\|f\|_{\infty},1) \max(\log(n)^{\alpha},1) \frac{ 2^{ 2 j \nu}}{n}. \label{eq:c1bis}
\end{equation}
Now let $\gamma_{jk} =  \sqrt{2  \delta \log(n) \frac{\eta_{j}^{2}}{n^{2}}   \hat{V}_{jk}}   +  \kappa  \frac{\eta_{j}^{2}}{n^{2}}  \delta \log(n) $, and remark that by definition of the threshold $\tau_{jk}$ and by using Lemmas \ref{lemma:ineg1} and \ref{lemma:ineg2} it follows that
\begin{eqnarray}
\PP \left( |\hat{\beta}_{jk} - \beta_{jk} | \geq \frac{\tau_{jk}}{\sqrt{1+\epsilon}}  \right) & = & \PP \left( |\hat{\beta}_{jk} - \beta_{jk} | \geq \frac{\tau_{jk}}{\sqrt{1+\epsilon}} ; \hat{V}_{jk} + \gamma_{jk} \geq V_{jk} \right) \nonumber \\
& & +  \PP \left( |\hat{\beta}_{jk} - \beta_{jk} | \geq \frac{\tau_{jk}}{\sqrt{1+\epsilon}} ; \hat{V}_{jk} + \gamma_{jk} \leq V_{jk} \right)  \nonumber \\
& \leq &  \PP \left( |\hat{\beta}_{jk} - \beta_{jk} | \geq \sqrt{ \frac{2 V_{jk} \delta \log(n)}{1+\epsilon} } + \frac{\eta_{j}}{3 n} \frac{\delta \log(n)}{1+\epsilon} \right)  \nonumber \\
& & + \PP \left(V_{jk} \geq  \hat{V}_{jk} + \gamma_{jk}  \right) \nonumber \\
& \leq & C (n^{-\frac{\delta}{1+\epsilon}} + n^{-\delta}) \leq C n^{-\frac{\delta}{1+\epsilon}} \label{eq:c2}
\end{eqnarray}
Now inserting (\ref{eq:c1bis}) and  (\ref{eq:c2}) into inequality (\ref{eq:holder}), and using the definition of $j_{1}$ and $j_{0}$ one finally obtains that for any $q > 1$ and $\epsilon > 0$
\begin{eqnarray*}
\EE Z&  \leq & C \max(\|f\|_{\infty},1) \max(\log(n)^{\alpha},1) \sum_{j = j_{0}}^{j_{1}} \sum_{k=0}^{2^{j}-1} 2^{2j \nu} n^{-1 -\frac{\delta}{q(1+\epsilon)}} \\
& \leq & C  \max(\|f\|_{\infty},1) \max(\log(n)^{\alpha},1) n^{-1 -\frac{\delta}{q(1+\epsilon)}} 2^{j_{1}(2 \nu +1)}.
\end{eqnarray*}
By definition of $j_{1}$, $2^{j_{1}(2 \nu +1)} \leq C n^{\eta \frac{2 \nu +1}{\nu +1}} (\log n)^{\alpha (2 \nu +1)} = C n^{\eta\left(1 + \frac{\nu}{\nu+1}\right)}  (\log n)^{\alpha (2 \nu +1)}$ which  implies that 
$$
\EE Z \leq C \max(\|f\|_{\infty},1) \max(\log(n)^{\alpha},1) n^{-1 + \eta\left(1 + \frac{\nu}{\nu+1}\right)-\frac{\delta}{q(1+\epsilon)}}  (\log n)^{\alpha (2 \nu +1)}
$$
By assumption, $\delta > \eta\left(1 + \frac{\nu}{\nu+1}\right)$. Hence there exists $\epsilon > 0$ such that $\delta > (1 + \epsilon)\eta\left(1 + \frac{\nu}{\nu+1}\right)$, and one can then always find some $q > 1$ such that $\delta > q(1+ \epsilon)\eta\left(1 + \frac{\nu}{\nu+1}\right)$. This implies that
$$
\EE Z \leq C \max(\|f\|_{\infty},1) \max(\log(n)^{\alpha},1) n^{-1}  (\log n)^{\alpha (2 \nu +1)}
$$
 which completes the proof.  \hfill $\Box$\\
 
Now, using the inequality $(a+b)^{2} \leq 2a^{2} + 2b^{2}$ one has that
$$
\EE \tau_{jk}^{2} \leq 4 \delta \log(n) \EE \tilde{V}_{jk} + 2 \frac{\eta_{j}^{2}}{9 n^{2}} \delta \log(n),
$$ 
where $\tilde{V}_{jk} =  \hat{V}_{j,k}  + \sqrt{2 \delta \log(n) \hat{V}_{j,k}  \frac{\eta_{j}^{2}}{n^{2}}   } + \delta \log(n) \kappa \frac{\eta_{j}^{2}}{n^{2}}$. Then using the inequality $2 a b \leq a^{2} + b^{2}$ and the fact that $\EE \hat{V}_{j,k}  = V_{jk}$, it follows that 
$$
 \EE \tilde{V}_{jk} \leq 3 V_{jk} +  (1+\kappa) \delta \log(n)  \frac{\eta_{j}^{2}}{n^{2}},
$$
which finally implies that there exists a constant $C$ such that
\begin{equation}
\EE \tau_{jk}^{2} \leq C \left( \log(n) V_{jk} + (\log n)^{2}  \frac{\eta_{j}^{2}}{n^{2}} \right) . \label{eq:taujk}
\end{equation}
Now, using Proposition \ref{prop:oracle}, Lemma \ref{control2} and inequality (\ref{eq:taujk}), it follows that there exists two constants $C(\delta)$ and $C(\delta)'$ not depending on $n$ and $f$, such that for any subset of indices $m \subset \Lambda_{n}$
\begin{eqnarray}
\EE \|\hat{f}_{n} - f \|^{2} & \leq & \sum_{k=0}^{2^{j_{0}}-1} \sigma_{j_{0},k}^{2} +C(\delta) \left( \sum_{(j,k) \in \Lambda_{n} \backslash m} \beta_{jk}^{2} +  \sum_{(j,k) \in   m} (1+\log(n)) V_{jk} +  \sum_{j = j_{1}+1}^{+ \infty }\sum_{k = 0}^{2^{j}-1} \beta_{jk}^{2}  \right) \nonumber \\
& & + C'(\delta) \left(\max(\|f\|_{\infty},1) \max(\log(n)^{\alpha},1) \frac{ (\log n)^{\alpha (2 \nu +1)}}{n} +  (\log n)^{2}   \sum_{(j,k) \in   m}  \frac{\eta_{j}^{2}}{n^{2}} \right), \label{eq:oracle}
\end{eqnarray}
and one can easily check that $\lim_{\delta \to \eta\left(1 + \frac{\nu}{\nu+1}\right)} C(\delta) = \lim_{\delta \to \eta\left(1 + \frac{\nu}{\nu+1}\right)} C'(\delta) = + \infty$
Then, as $\eta_{j} =  \sum_{\ell  \in C_{j}} |\tilde{\psi}_{jk}^{\ell}|$ it follows from (\ref{eq:boundpsi}) and (\ref{eq:boundh}) that $\eta_{j}^{2} \leq C 2^{j (2\nu+1)}$, which implies that for any $m \subset \Lambda_{n}$, $  \sum_{(j,k) \in   m}  \frac{\eta_{j}^{2}}{n^{2}}  \leq  C \frac{1}{n^{2}} \sum_{j = j_{0}}^{j_{1}} \sum_{k=0}^{2^{j}-1}   2^{j (2\nu+1)} \leq C\frac{2^{j_{1}(2 \nu +2)}}{n^{2}} \leq C n^{2(\eta-1)} (\log n)^{\alpha(2 \nu +2)}$ by definition of $j_{1}$. Inserting this inequality in (\ref{eq:oracle}) with the model $m = \left\{  (j,k); |\beta_{jk}|^{2} \geq \log(n) V_{jk} \right\}$ finally yields that there exists two constants $C(\delta)$ and $C(\delta)'$ such that
\begin{eqnarray*}
\EE \|\hat{f}_{n} - f \|^{2} & \leq &  \sum_{k=0}^{2^{j_{0}}-1} \sigma_{j_{0},k}^{2} +C(\delta)\left( \sum_{j = j_{0}}^{j_{1}} \sum_{k=0}^{2^{j}-1} \min( \beta_{jk}^{2},\log(n)V_{jk} ) +  \sum_{j = j_{1}+1}^{+ \infty }\sum_{k = 0}^{2^{j}-1} \beta_{jk}^{2}   \right) \\
& & + C'(\delta) \left( \max(\|f\|_{\infty},1) \max(\log(n)^{\alpha},1)  \frac{(\log n)^{\alpha (2 \nu +1)}}{n} +  \frac{(\log n)^{2 + \alpha (2 \nu +2)  }}{n^{2(1-\eta)}} \right).
\end{eqnarray*}
To finally obtain inequality (\ref{eq:oracle}), remark that $V_{j,k} = \sigma_{j,k}^{2} + \frac{1}{n}\beta_{j,k}^{2}$, and one can easily check that 
$$
\min( \beta_{jk}^{2},\log(n)V_{jk} ) \leq  \min( \beta_{jk}^{2},\log(n)\sigma^{2}_{jk}) + \frac{\log(n)}{n}  \beta_{jk}^{2},
$$
which completes the proof of Theorem \ref{eq:ineqoracledec}.  \hfill $\Box$

\subsection{Proof of Theorem \ref{theo:minimaxdec}}

Given our assumptions, one has that $1 \leq p \leq 2$ which implies that $s^{\ast} = s +1/2-1/p$. First we need the following lemma:

\begin{lemma} \label{lem:bound}
If $f \in B^{s}_{p,q}(A)$  with $1 \leq p \leq 2$ then
\begin{equation}
\sum_{k = 0}^{2^{j}-1} \beta_{jk}^{2} \leq A^{2} 2^{-2js^{\ast}}. \label{eq:john}
\end{equation}
Moreover, if $s > 1/p + 1/2$ with $1 \leq p \leq 2$, then there exists a constant $B > 0$ such that
$$
\sup_{f \in  B^{s}_{p,q}(A)} \|f\|_{\infty} \leq B.
$$
\end{lemma}
{\bf Proof:} since $f \in B^{s}_{p,q}(A)$ one has that
$
 \sum_{k =0}^{2^{j_{0}}-1} |c_{j_{0},k}|^{p} \leq A^{p}
$\
and
$
\sum_{k=0}^{2^{j}-1} |\beta_{j,k}|^{p} \leq A^{p} 2^{-jp(s + 1/2-1/p)}.
$
Since $p \leq 2$ it follows that
\begin{equation} \label{eq:boundalpha}
\left( \sum_{k =0}^{2^{j_{0}}-1} |c_{j_{0},k}|^{2} \right)^{1/2} \leq \left( \sum_{k =0}^{2^{j_{0}}-1} |c_{j_{0},k}|^{p} \right)^{1/p} \leq A,
\end{equation}
and
\begin{equation} \label{eq:boundbeta}
\left( \sum_{k=0}^{2^{j}-1} |\beta_{j,k}|^{2} \right)^{1/2} \leq \left( \sum_{k=0}^{2^{j}-1} |\beta_{j,k}|^{p} \right)^{1/p} \leq A 2^{-j(s + 1/2-1/p)},
\end{equation}
which proves the first part of the Lemma. Then, remark that
\begin{equation} \label{eq:boundf}
\|f\|_{\infty} \leq \| \sum_{k =0}^{2^{j_{0}}-1} c_{j_{0},k} \phi_{j_{0},k} \|_{\infty} +  \sum_{j = j_{0}}^{+ \infty} \|  \sum_{k=0}^{2^{j}-1} \beta_{j,k} \psi_{j,k}  \|_{\infty}
\end{equation}
Now, let $x \in [0,1]$ and remark that by Cauchy-Schwartz inequality 
$$
\left| \sum_{k =0}^{2^{j_{0}}-1} c_{j_{0},k} \phi_{j_{0},k}(x) \right|^{2}  \leq 2^{2j_{0}} \|\phi\|_{\infty}^{2} \left(\sum_{k =0}^{2^{j_{0}}-1} |c_{j_{0},k} |^{2}  \right),
$$
and
$$
\left| \sum_{k=0}^{2^{j}-1} \beta_{j,k} \psi_{j,k}(x) \right|^{2} \leq 2^{2j} \|\psi\|_{\infty}^{2} \left(\sum_{k =0}^{2^{j_{1}}-1} |\beta_{j,k} |^{2}  \right).
$$
Hence the above inequalities, (\ref{eq:boundalpha}), (\ref{eq:boundbeta}), and (\ref{eq:boundf}) imply that
$$
\|f\|_{\infty} \leq A  2^{j_{0}} \|\phi\|_{\infty}  +   A \|\psi\|_{\infty} \sum_{j = j_{0}}^{+ \infty} 2^{-j(s - 1/2-1/p)}
$$
By assumption $s - 1/2-1/p > 0$ which implies that $ \sum_{j = j_{0}}^{+ \infty} 2^{-j(s - 1/2-1/p)} < + \infty$ , and thus the result follows with $B = A\left(2^{j_{0}} \|\phi\|_{\infty} + \|\psi\|_{\infty} \sum_{j = j_{0}}^{+ \infty} 2^{-j(s - 1/2-1/p)}\right)$. \hfill $\Box$\\

Note that although not always stated  Lemma \ref{lem:bound}  will imply that the various bounds given below hold uniformly for $f \in D^{s}_{p,q}(A)$.\\

\noindent Let $R_{n} =  R_{1} + R_{2} + R_{3}$ with
$$
R_{1} = \sum_{k=0}^{2^{j_{0}}-1} \sigma_{j_{0},k}^{2}, \quad R_{2} = \sum_{j = j_{0}}^{j_{1}} \sum_{k=0}^{2^{j}-1} \min( \beta_{jk}^{2}, \log(n) V_{jk} )   \quad \mbox{and}  \quad  R_{3} =  \sum_{j = j_{1}+1}^{+ \infty }\sum_{k = 0}^{2^{j}-1} \beta_{jk}^{2}.
$$
Given our assumptions on $f$ and since $\eta = 1/2$, Theorem \ref{theo:oracledec} and Lemma \ref{lem:bound} imply that there exists two constants $C$ and $C'$  such that for every $n > \exp(1)$ and all $f \in D^{s}_{p,q}(A)$
\begin{equation}
\EE \|\hat{f}_{n} - f \|^{2} \leq C  R_{n}   +  C' \frac{ (\log n)^{2}}{n} \label{eq:risk}
\end{equation}
Hence to study the rate of convergence of $\hat{f}_{n}$, it suffices to study the asymptotic behavior of $R_{n}$. From Lemma \ref{lemma:var} and by definition of $j_{0}$ one has that $R_{1} \leq C \frac{2^{j_{0}(2 \nu +1)}}{n} \leq C  \frac{(\log n )^{(2 \nu +1)}}{n}$ which implies that
\begin{equation}
R_{1} = \opO \left( n^{\frac{-2s}{2s+2\nu+1}} \right)   \mbox{ (in the dense case) or } R_{1} = \opO \left( n^{\frac{-2s^{\ast}}{2s^{\ast}+2\nu}} \right) \mbox{ (in the sparse case) }. \label{eq:R1}
\end{equation}
Then remark that Lemma \ref{lem:bound} implies that $R_{3} =  \opO \left(   2^{-2j_{1}s^{\ast}} \right) =  \opO \left( n^{-\frac{2 s^{\ast}}{2 \nu +2}}(\log n)^{-2\alpha s^{\ast}} \right) $ by definition of $j_{1}$. Now remark that if $p = 2$ then $s^{\ast} = s > 1$ (by assumption) and thus $\frac{2 s^{\ast}}{2 \nu +2} > \frac{2 s}{2s + 2\nu +1}$. If $1 \leq p < 2$ then $s^{\ast} = s + 1/2-1/p$, and one can check that the condition $s > 1/2 + 1/p$ implies that  $\frac{2 s^{\ast}}{2 \nu +2} > \frac{2 s}{2s + 2\nu +1}$ if $\nu(2-p) < ps^{\ast}$, and that $\frac{2 s^{\ast}}{2 \nu +2} > \frac{2 s^{\ast}}{2s^{\ast} + 2\nu}$ if $\nu(2-p) \geq ps^{\ast}$. Hence one obtains that
\begin{eqnarray}
R_{3} & = & \opO \left( n^{-\frac{2s}{2s+2\nu+1}} \right) \mbox{ if } \nu(2-p) < ps^{\ast} ,\label{eq:R3dense} \\
R_{3} & = & \opO \left( n^{-\frac{2 s^{\ast}}{2s^{\ast} + 2\nu}} \right) \mbox{ if } \nu(2-p) \geq ps^{\ast}. \label{eq:R3sparse} \\
\end{eqnarray}
Now it remains to study the term $R_{2}$. \\

For this let us consider first the dense case when $\nu(2-p) < ps^{\ast}$, and decompose $R_{2} = R_{21} + R_{22}$ with
$$
 R_{21} = \sum_{j = j_{0}}^{j_{2}} \sum_{k=0}^{2^{j}-1} \min( \beta_{jk}^{2}, \log(n) V_{jk} )  \mbox{ and }  R_{22} = \sum_{j = j_{2}+1}^{j_{1}} \sum_{k=0}^{2^{j}-1} \min( \beta_{jk}^{2}, \log(n) V_{jk} ) , 
$$
where  $j_{2} = j_{2}(n)$ is the integer such that $2^{j_{2}} > (n/\log(n))^{\frac{1}{2s + 2\nu +1}} \geq 2^{j_{2}-1}$ (note that given our assumptions $j_{2} \leq j_{1}$ for all sufficiently large $n$). Then remark that by Lemma \ref{lemma:var}, $ R_{21} \leq \sum_{j = j_{0}}^{j_{2}} \sum_{k=0}^{2^{j}-1} \log(n) V_{jk} \leq C \sum_{j = j_{0}}^{j_{2}} \log(n) \frac{2^{j(2\nu +1)}}{n} \leq C \log(n) \frac{2^{j_{2}(2\nu +1)}}{n} \leq C (n/\log(n))^{-\frac{2s}{2s+2\nu+1}} $. Hence
\begin{equation}
R_{21} = \opO \left( (n/\log(n))^{-\frac{2s}{2s+2\nu+1}} \right) \label{eq:R21dense}
\end{equation}
Now remark that if $p \geq 2$ then by equation (\ref{eq:john}) and by definition of $j_{1}$ and $j_{2}$  it follows that $R_{22} \leq C \sum_{j = j_{2}+1}^{j_{1}} 2^{-2js} \leq C \left(n ^{-\frac{2s}{2\nu+2}}  - (n/\log(n))^{-\frac{2s}{2s+2\nu+1} } \right)$ which implies that
\begin{equation}
R_{22} = \opO \left( (n/\log(n))^{-\frac{2s}{2s+2\nu+1}} \right) \label{eq:R22dense1}.
\end{equation}
Now if $1 \leq p < 2$ remark that $R_{22}$ can be written as
\begin{eqnarray*}
R_{22} & =  & \sum_{j = j_{2}+1}^{j_{1}} \sum_{k=0}^{2^{j}-1}  \beta_{jk}^{2} \1_{\{ \beta_{jk}^{2} < \log(n) V_{jk}\}} +  \log(n) V_{jk} \1_{ \{\beta_{jk}^{2} > \log(n) V_{jk}\} } \\
& =  & \sum_{j = j_{2}+1}^{j_{1}} \sum_{k=0}^{2^{j}-1}  |\beta_{jk}|^{2-p}  |\beta_{jk}|^{p} \1_{\{ \beta_{jk}^{2} < \log(n) V_{jk}\}} +  (\log(n) V_{jk}) ^{1-p/2} (\log(n) V_{jk})^{p/2}  \1_{ \{\beta_{jk}^{2} > \log(n) V_{jk}\} }\\
& \leq & 2  \sum_{j = j_{2}+1}^{j_{1}} \sum_{k=0}^{2^{j}-1} ( \log(n) V_{jk} )^{1-p/2}    |\beta_{jk}|^{p}. 
\end{eqnarray*}
By Lemma \ref{lemma:var} $V_{jk} \leq C \frac{2^{2j \nu}}{n}$,  and since $f   \in B^{s}_{p,q}(A)$ if follows that there exists a constant $C$ depending only on $p,q,s,A$ such that $\sum_{k=0}^{2^{j}-1}  |\beta_{jk}|^{p} \leq C 2^{-j s ^{\ast}p}$, which implies that
$$
R_{22} \leq C (n/\log(n))^{-1+p/2} \sum_{j = j_{2}+1}^{j_{1}} 2^{j(2 \nu - \nu p-s^{\ast}p)}.
$$
Now in the dense case, one has that $2 \nu - \nu p-s^{\ast}p < 0$ which implies that
\begin{equation}
R_{22} = \opO \left( (n/\log(n))^{-1+p/2}  2^{j_{2}(2 \nu - \nu p-s^{\ast}p)} \right) =  \opO \left( (n/\log(n))^{-\frac{2s}{2s + 2\nu +1}} \right) \label{eq:R22dense2}
\end{equation}
by using the definition of $j_{2}$. Hence combining (\ref{eq:R21dense}), (\ref{eq:R22dense1}) and  (\ref{eq:R22dense2}) it follows that in the dense case for $1 \leq p \leq \infty$
\begin{equation}
R_{2} = \opO \left( (n/\log(n))^{-\frac{2s}{2s+2\nu+1}} \right) \label{eq:R2dense}.
\end{equation}

Now consider the sparse case when $\nu(2-p) \geq ps^{\ast}$, and decompose $R_{2} = R_{21} + R_{22}$  with
$$
 R_{21} = \sum_{j = j_{0}}^{j_{2}} \sum_{k=0}^{2^{j}-1} \min( \beta_{jk}^{2}, V_{jk} )  \mbox{ and }  R_{22} = \sum_{j = j_{2}+1}^{j_{1}} \sum_{k=0}^{2^{j}-1} \min( \beta_{jk}^{2}, V_{jk} ) , 
$$
where  $j_{2} = j_{2}(n)$ is the integer such that $2^{j_{2}} > (n/\log(n))^{\frac{1}{2s^{\ast}+ 2\nu}} \geq 2^{j_{2}-1}$. Note that in the sparse case then necessarily $1 \leq p < 2$, and as previously one can thus remark that $R_{21}$ can be written as
$$
R_{21}  =  \sum_{j = j_{0}}^{j_{2}} \sum_{k=0}^{2^{j}-1}  \beta_{jk}^{2} \1_{\{ \beta_{jk}^{2} < \log(n) V_{jk}\}} +  \log(n) V_{jk} \1_{ \{\beta_{jk}^{2} > \log(n) V_{jk}\} }  \leq 2  \sum_{j = j_{0}}^{j_{2}} \sum_{k=0}^{2^{j}-1}  (\log(n) V_{jk})^{1-p/2}    |\beta_{jk}|^{p}. 
$$
% and since $ 2 \nu - \nu p-s^{\ast}p \geq 0$
Now using again the fact that  $V_{jk} \leq C \frac{2^{2j \nu}}{n}$ (by Lemma \ref{lemma:var}), that  $\sum_{k=0}^{2^{j}-1}  |\beta_{jk}|^{p} \leq C 2^{-j s ^{\ast}p}$ (since $f   \in B^{s}_{p,q}(A)$),  and by  definition of $j_{2}$ and $j_{0}$ it follows that
\begin{equation}
R_{21}  \leq C (n/\log(n))^{-1 + p/2} \left( 2^{j_{2} (2 \nu -p\nu -p s^{\ast}) } - 2^{j_{0} (2 \nu -p\nu -p s^{\ast}) } \right) = \opO \left( (n/\log(n))^{-\frac{2s^{\ast}}{2s^{\ast} + 2 \nu}}\right). \label{eq:R21sparse}
\end{equation}
Then, remark that by equation (\ref{eq:john}), $R_{22} \leq   \sum_{j = j_{2}+1}^{j_{1}} \sum_{k=0}^{2^{j}-1}	\beta_{jk}^{2} \leq  \sum_{j = j_{2}+1}^{j_{1}}  2^{-2js^{\ast}}$. Hence, $R_{22} = \opO \left( 2^{-2j_{2}s^{\ast}} \right)$, and thus by definition of $j_{2}$
\begin{equation}
R_{22}  = \opO \left( (n/\log(n))^{-\frac{2s^{\ast}}{2s^{\ast} + 2 \nu}}\right). \label{eq:R22sparse}
\end{equation}
Finally combining (\ref{eq:R21sparse}) and (\ref{eq:R22sparse}) it follows that in the dense case
\begin{equation}
R_{2} = \opO \left( (n/\log(n))^{-\frac{2s^{\ast}}{2s^{\ast} + 2 \nu}}\right) \label{eq:R2sparse}.
\end{equation}
Then, combining (\ref{eq:R1}),  (\ref{eq:R3dense}),  (\ref{eq:R3sparse}), (\ref{eq:R2dense}) and (\ref{eq:R2sparse}) implies that $R_{n} = \left( (n/\log(n))^{-\frac{2s}{2s + 2 \nu + 1}}\right)$  in the dense case, and $R_{n} = \left( (n/\log(n))^{-\frac{2s^{\ast}}{2s^{\ast} + 2 \nu}}\right)$ in the sparse case. Using inequality (\ref{eq:risk}) this completes the proof of Theorem \ref{theo:minimaxdec}.  \hfill $\Box$

\bibliography{DensDecOracle}
\bibliographystyle{apacite}

\end{document}